\documentclass[12pt]{elsart}
\usepackage[a4paper,text={16.0cm,22.0cm},centering]{geometry}
\usepackage{graphicx}
\usepackage{amsfonts}
\usepackage{amsmath}
\usepackage{amssymb}
\usepackage{lscape}
\usepackage[all]{xy}
\usepackage{mathrsfs}
\usepackage{psfrag}
\usepackage{bm}
\usepackage{color}

\usepackage[textmath,displaymath]{preview}

\def\lbk{\lbrack\!\lbrack}
\def\lbkbig{\left\lbrack\!\left\lbrack}
\def\rbk{\rbrack\!\rbrack}
\def\rbkbig{\right\rbrack\!\right\rbrack}

\def\i{\mathrm i}
\def\re{\mathrm e}
\def\d{\mathrm d}

\def\Z{\mathbb Z}
\def\R{\mathbb R}
\def\P{\mathbb P}
\def\C{\mathbb C}

\def\sech{\mathrm {sech}}
\def\ess{\mathrm {ess}}

\def\Det2{\mathrm{Det}^2}
\def\fr{\mbox{$\frac{1}{2}$}}
\def\bw{\mbox{$\bigwedge^2(\mathbb{R}^{4})$}}
\def\ri{{\rm i}}
\def\rd{{\rm d}}
\def\bwedge{\mbox{$\bigwedge$}}
\def\unit{{\rm U}(2)}
\def\vol{{\textsf{vol}}}
\def\symp{{\bm\omega}}
\def\rp{{\mathbb R}\mathbb{P}}

\def\cwedgeknn{\hbox{$\bigwedge^{k}({\mathbb R}^{2n})$}}

\newcommand{\derp}[3]{\frac
  {\ifthenelse{\equal{#3}{1}}{\partial #1}{\partial^{#3} #1}}
  {\ifthenelse{\equal{#3}{1}}{\partial #2}{\partial #2^{#3}}}
}

\renewcommand{\theequation}{\arabic{section}.\arabic{equation}}

\begin{document}

\begin{frontmatter}

\title{Computational aspects of the Maslov index of solitary waves}

\author{Fr\'ed\'eric Chardard$^a$, Fr\'ed\'eric Dias$^a$, Thomas J. Bridges$^b$}

\address[a]{Centre de Math\'ematiques
et de Leurs Applications, Ecole Normale Sup\'erieure de Cachan,
61 avenue du Pr\'esident Wilson, 94235 Cachan cedex, France}
\address[b]{Department of Mathematics, University of Surrey,
Guildford GU2 7XH, England}

\begin{abstract}
When solitary waves are characterized as homoclinic orbits of a finite-dimensional
Hamiltonian system, they have an integer-valued topological invariant, the Maslov index.
We are interested in developing a robust numerical algorithm to compute
the Maslov index, to understand its properties, and to study the
implications for the stability of solitary waves.
The algorithms reported here 
are developed in the exterior algebra representation, which leads to a robust and
fast algorithm with some novel properties.  We use two different
representations for the Maslov index, one based on an intersection index and one
based on approximating the homoclinic orbit by a sequence of periodic orbits.
New results on the Maslov index for solitary wave
solutions of reaction-diffusion equations,
the fifth-order Korteweg-De Vries equation,
and the longwave-shortwave resonance equations are presented.  
Part 1 considers the case of four-dimensional
phase space, and Part 2 considers the case of $2n-$dimensional phase space with $n>2$.
\end{abstract}

\end{frontmatter}

\section{Introduction}
\label{intro}
\setcounter{equation}{0}

Hamiltonian evolution equations in one space dimension, such as the nonlinear
Schr\"odinger (NLS) equation, fifth-order Korteweg-De Vries (KdV) equation, longwave-shortwave
resonance (LW-SW) equations, have the property that their steady part is a finite-dimensional
Hamiltonian system.  For such systems, solitary wave solutions
can be characterized as homoclinic orbits of the
Hamiltonian ordinary differential equation (ODE).  The spectral problem
associated with the linearization about a given
homoclinic orbit, in the time-dependent equations, then leads to a parameter-dependent
family of linear Hamiltonian systems.  The advantage of these Hamiltonian
structures is that the linear and nonlinear
Hamiltonian systems have global geometric properties that aid in proving existence of
the basic solitary wave and in understanding its stability as a solution of the
time-dependent equation. Our interest in this paper is in a particular
geometric invariant -- the Maslov index of homoclinic orbits.

The study of the stability of solitary waves using the
Maslov index was pioneered in the papers by \textsc{Jones}~\cite{jones}
and \textsc{Bose \& Jones}~\cite{bose}.  
The linear stability of
standing wave solutions of a spatially-dependent NLS equation is studied in \cite{jones}.
The linearization about a steady solution results in a linear $\lambda-$dependent
Hamiltonian system with two degrees of freedom of the form (\ref{h1}) 
and $\lambda$ a spectral parameter.   Geometric methods are then used to
determine the Maslov index, and it is used to prove an instability result.
Gradient parabolic PDEs of the form
\begin{equation}\label{grad-pde}
\begin{array}{rcl}
u_t &=& d_1 u_{xx} + f_u(u,v) \\[2mm]
v_t &=& d_2 v_{xx} + f_v(u,v) \end{array}
\end{equation}
are considered in \cite{bose}, where $d_1$ and $d_2$ are positive parameters, $f(u,v)$ is
a given smooth function with gradient $(f_u,f_v)$.  Linearizing
about a steady solution $(\widehat u(x),\widehat v(x))$, and
introducing a spectral parameter leads to a coupled pair of linear
second-order ODEs which can be put into the standard form (\ref{h1}),
with the asymptotic
property (\ref{A_infty}) and $\lambda$ the spectral parameter.  Since
the PDE is a gradient system it is sufficient to restrict the spectral
parameter to be real. Singular perturbation methods are then used
to determine the Maslov index, which in turn is related to stability.
A key feature of this work is the analysis of the induced system on
the exterior algebra space $\bwedge^2(\R^4)$.

Many of the most interesting solitary waves
are only known numerically and therefore
a numerical approach to the Maslov index is of interest.
It is the aim of this paper to develop a numerical framework 
for computing the Maslov index of homoclinic orbits.   Once the solitary
wave solution is known, analytically or numerically, it is the linearization
about that solitary wave which encodes the Maslov index.  Therefore, the
starting point for developing the theory is the following class of parameter-dependent
Hamiltonian systems
\begin{equation}\label{h1}
{\bf J}{\bf u}_x = {\bf B}(x,\lambda){\bf u}\,,\quad {\bf u}\in\R^{4}\,,
\quad x\in\R\,,\quad \lambda\in\R\,,
\end{equation}
where ${\bf J}$ is the standard symplectic operator on $\R^{4}$
\begin{equation}\label{J-def}
{\bf J} = \left[\begin{matrix}
\ {\bf 0} & - {\bf I}\ \\ \ {\bf I} & \hfill{\bf 0}\ \end{matrix}\right]\,,
\end{equation}
and ${\bf B}(x,\lambda)$ is a symmetric matrix depending smoothly on
$x$ and $\lambda$. Let
\begin{equation}\label{A-def}
{\bf A}(x,\lambda) = {\bf J}^{-1}{\bf B}(x,\lambda)\,.
\end{equation}
The fact that ${\bf A}(x,\lambda)$ is obtained from the linearization about a solitary
wave suggests the following asymptotic property.  
It is assumed throughout the paper that
\begin{equation}\label{A_infty}
{\bf A}_\infty(\lambda)=\lim_{x\to\pm\infty}{\bf A}(x,\lambda)\,,
\end{equation}
and that ${\bf A}_\infty(\lambda)$ is strictly hyperbolic for an
open set of $\lambda$ values that includes $0$.

The theory applies to linear Hamiltonian systems on a phase space of dimension $2n$ with
$n$ any finite natural number.  In Part 1 attention is restricted to the case of $4-$dimensional
phase space which simplifies formulae, and general aspects of the case $n>2$ are given in
Part 2 \cite{cdb-part2}.

The Maslov index is a winding number associated with paths of solutions
of (\ref{h1}), in particular, paths of Lagrangian planes.
A Lagrangian plane in $\R^4$ is
a $2-$dimensional subspace of $\R^{4}$, say ${\rm span}\{{\bf z}_1,{\bf z}_2\}$, satisfying
$\langle{\bf J}{\bf z}_1,{\bf z}_2\rangle=0$,
where $\langle\cdot,\cdot\rangle$ is a standard inner product on $\R^{4}$.

Suppose $\lambda$ is fixed and on the interval $a\leq x\leq b$
consider a path of Lagrangian planes
\[
[a,b] \mapsto\ {\bf Z}(x,\lambda) = [{\bf z}_1(x,\lambda)\,|\,{\bf z}_2(x,\lambda)]
\in\R^{4\times 2}\,,
\]
satisfying ${\bf Z}_x={\bf A}(x,\lambda){\bf Z}$ for $a\leq x\leq b$.
The Maslov index of this path is a count of the number of times
this path of Lagrangian planes has a non-trivial intersection with
a fixed reference Lagrangian plane.
A precise definition is given in \S\ref{sec2}.

A byproduct of the present theory is some new observations
about the properties of the Maslov index, which in turn are
useful in computation.
In the numerics, the exterior
algebra formulation is also advantageous.  We give formulas for different
representations of the Maslov index for Lagrangian planes on
$\bw$ (and for any $n$ in Part 2 \cite{cdb-part2}),
and present a general algorithm that works -- in principle -- for
any dimension $n$.  However, the dimension of
$\bigwedge^{n}(\R^{2n})$ increases rapidly with $n$ and
so the algorithm is most effective for low dimensional systems.
The algorithm is constructed so that the manifold of Lagrangian
planes is attracting. Numerical results are presented in this paper
for the cases of $\R^2$ and $\R^4$.

In order to develop a numerical framework for the Maslov index,
one of the first difficulties is defining the Maslov index.  Although
it is easy to give a rough definition, making it precise depends greatly
on the context and a surprising number of special difficulties and cases
arise.

We will appeal to two constructions of the Maslov index.  The first is
based on an intersection index between the Lagrangian
path and a reference plane.
This definition was used in Maslov's original work, and was developed
further by \textsc{Arnold}~\cite{arnold1} and \textsc{Duistermaat}~\cite{duistermaat}.
It is this definition that is used by \textsc{Jones}~\cite{jones}
and \textsc{Bose \& Jones}~\cite{bose}, taking the Lagrangian path to
be a path of unstable subspaces and taking the reference plane to be
the stable subspace at infinity.  

Independently, \textsc{Chen \& Hu}~\cite{chen_hu}
give two constructions of the Maslov index of a homoclinic orbit.   Their first definition
is based on an intersection index and is equivalent to the definition in
\cite{jones,bose}.  Their second definition is based on a Fredholm index of
(\ref{h1}) viewed as an operator in a function space on the real line.  
However, this latter definition, although equivalent to the definition 
based on an intersection index, is not convenient in the numerics since
it requires to compute a spectral projector of an infinite dimensional
operator and then to determine a Fredholm index.

In this paper, the definition of the Maslov index based on an intersection index
is extended by introducing an explicit and computable formula
for the crossing form.  This theory is developed in \S\ref{sec-paths}.

Our second method for computing the Maslov index is to approximate
the homoclinic orbit by a sequence of periodic orbits, apply the
Maslov index for periodic orbits, and then take limits.  
There does not appear to be any loss of generality in using periodic
approximates.  \textsc{Vanderbauwhede \& Fiedler}~\cite{vf92} prove
that homoclinic orbits in Hamiltonian systems (as well as reversible systems)
can be approximated as the limit of a sequence of periodic orbits.

The Maslov index for periodic
orbits has been widely developed because of its interest in
semi-classical quantization (e.g.\ \cite{gutzwiller,crl,lr87,robbins91,pb03,m-g03}
and references therein).
In \cite{F.Chardard2006}
a new numerical scheme is developed to compute the Maslov index 
of \emph{hyperbolic} periodic orbits, and 
\textsc{Chardard}~\cite{F.Chardard2007} proves under suitable
hypotheses that if the periodic orbit is asymptotic to a homoclinic orbit,
the Maslov index converges to the Maslov index of the limiting
homoclinic orbit.  
This approach ties in with 
existing schemes for computing the basic solitary wave,
where the solitary wave is approximated by a periodic
orbit and then a spectral method is used for computation.

The computational framework for the Maslov index is illustrated by
application to four examples.  The first is a tutorial example on $\R^2$, where
the details can be given explicitly.  It is a scalar-reaction diffusion equation
with an explicit localized solution.  The second example is a coupled reaction-diffusion
equation which also has an explicit solution.
The third example is solitary wave solutions of KdV5.
The fourth example is the LW-SW resonance
equation which arises in fluid mechanics and consists of a NLS equation coupled
to a KdV equation.  This latter example has two new interesting features: 
the spectral problem is on a six-dimensional phase space,
and for appropriate parameter values it has a Maslov index which is a non-monotone
function of $\lambda$.

\section{Linear Hamiltonian systems and Lagrangian subspaces}
\label{sec2}
\setcounter{equation}{0}

A Lagrangian subspace can be represented by a
\emph{Lagrangian frame}:
a $4\times 2$ matrix of rank $2$
\begin{equation}\label{Z-def}
{\bf Z} = \begin{pmatrix}{\bf X}\\ {\bf Y}\end{pmatrix}\,,
\end{equation}
where ${\bf X}$ and ${\bf Y}$ are $2\times 2$ matrices satisfying
\begin{equation}\label{Z-sym}
{\bf Y}^T{\bf X} = {\bf X}^T{\bf Y}\,.
\end{equation}
There is a correspondence between
elements of the unitary group $\unit$ and Lagrangian subspaces.
When ${\bf X}^T{\bf X} + {\bf Y}^T{\bf Y}={\bf I}$ then
${\bf X}\pm\ri{\bf Y}$ are unitary matrices.  The determinant
of a unitary matrix lies on the unit circle.  This property
suggests defining the angle
\begin{equation}\label{angle-def}
{\rm e}^{\ri\kappa} = \frac{{\rm det}[{\bf X} - \ri{\bf Y}]}{
{\rm det}[{\bf X} + \ri{\bf Y}]}\,.
\end{equation}
Along a path of Lagrangian subspaces, this angle
will change, and the winding of this angle is the basis
of the Maslov index.

Let ${\bf Z}(x)$, $a \leq x \leq b$ be any smooth path of Lagrangian subspaces.
If the path is a loop: ${\bf Z}(b) = {\bf Z}(a)$, then there is an integer
associated with the path: the number of times the induced 
path on the unit circle, represented by ${\rm e}^{\ri \kappa(x)}$, encircles the origin.  
Define the angle associated with the path by
\begin{equation}\label{angle-def-1}
{\rm e}^{\ri\kappa(x)} = \frac{{\rm det}[{\bf X}(x) - \ri{\bf Y}(x)]}{
{\rm det}[{\bf X}(x) + \ri{\bf Y}(x)]}\,.
\end{equation}
Then the Maslov index of the path is the integer
\begin{equation}\label{maslov-angle}
\textsf{Maslov}({\bf Z},\kappa) := \frac{ \kappa(b) - \kappa(a) }{2\pi}\,,
\end{equation}
where $\kappa$ here is viewed as the lift from $S^1$ to $\R$.

Geometrically $\operatorname{Span}(\mathbf Z)\to \textsf{Maslov}(\mathbf Z,\kappa)$ induces an isomorphism
between the homotopy group of $\Lambda(n)$ at a point $A$ and the group
of integers $(\Z,+)$.

The Maslov index can also be defined for an arbitrary path of
Lagrangian subspaces by introducing the idea of an intersection form.
This approach to the Maslov index goes back to \textsc{Maslov} and
\textsc{Arnold}~\cite{arnold1}.  The key to the definition in this
case is the use of an \emph{intersection form} or \emph{crossing form}.
Here we will follow the construction in
\textsc{Robbins}~\cite{robbins91,robbins92} and
\textsc{Robbin \& Salamon}~\cite{rs93}.  Modulo a choice of orientation
these definitions are equivalent.

Let $V$ be a fixed reference plane, represented
by a Lagrangian frame.  For example a typical choice for the reference
plane is
\[
V = \begin{pmatrix} {\bf I} \\ {\bf 0} \end{pmatrix}\,.
\]
In the case of homoclinic orbits a natural choice for the reference
plane is the stable or unstable manifold at some value of $x$.

Suppose, for some value
of $x$, denoted $x_0$, there is a simple intersection between
the reference plane $V$ and the path:
that is, ${\bf Z}(x_0)\cap V$ is
one dimensional.  The intersection index at $x_0$ is determined by
the sign of $\Gamma$, the crossing form, defined by
\begin{equation}\label{Gamma-def}
\Gamma({\bf Z},V,x_0) =  \langle {\bf J}{\bf Z}'(x_0)\beta,
{\bf Z}(x_0)\beta\rangle\,\vol\,.
\end{equation}
Here and throughout the paper $\vol$ is taken to be the standard
volume form on $\R^4$,
\[
\vol = {\bf e}_1\wedge{\bf e}_2\wedge{\bf e}_3\wedge{\bf e}_4\,.
\]
In the formula (\ref{Gamma-def}),
\[
{\bf Z}'(x_0):= \frac{d\ }{dx}\bigg|_{x=x_0}{\bf Z}(x),
\] 
and  $\beta\in\rp^{1}$, where  $\rp^{1}$ is the one-dimensional real
projective space\footnote{\[
\rp^{1} = \{\ \beta\in\R^2\ :\ \beta\neq0\,,\quad
\beta \sim c\,\beta\,,\quad
c\in\R\setminus\{0\}\ \}\,.
\]
}.
The parameter $\beta$ determines the linear combination
of the columns of ${\bf Z}(x_0)$
which span the intersection subspace
\[
{\bf Z}(x_0)\cap V = {\rm span}\{\xi \}\,,\quad
\xi:=\beta_1{\bf z}_1+\beta_2{\bf z}_2\,,
\]
where ${\bf z}_1,{\bf z}_2$ are the columns of ${\bf Z}(x_0)$.
At each simple intersection between $V$ and the path ${\bf Z}(x)$ the
sign of the intersection form is $\pm 1$.  
Adding the intersection indices over the path
gives the Maslov index
\begin{equation}\label{Maslov-def}
\textsf{Maslov}({\bf Z},V) =
\sum_{a<x_0<b} {\rm sign}\,\Gamma({\bf Z},V,x_0)\,.
\end{equation}
This formula assumes that ${\bf Z}(a)\cap V = \{0\}$ and
${\bf Z}(b)\cap V = \{0\}$ and that there is only a finite
number of $x_0$ where ${\bf Z}(x_0)\cap V \neq \{0\}$.
When the endpoints have non-trivial intersection the formula can be modified
to contribute a half-integer for each end intersection (see page 831 of \cite{rs93}).
Non-regular intersections are not generic in the one-parameter family, and so
can be eliminated by perturbing $\lambda$.

Geometrically, two paths $\operatorname{image}(\mathbf Z),
\operatorname{image}(\mathbf W)$ in $\Lambda(n)$ with the same endpoints
and such that $\textsf{Maslov}({\bf Z},V) $ and $\textsf{Maslov}({\bf W},V) $
are well-defined, 
are in  homotopy if and only if there exists $V$ such that
$\textsf{Maslov}({\bf Z},V)=\textsf{Maslov}({\bf W},V)$.

\section{The Evans function associated to the set of linear Hamiltonian systems}
\label{sec-evans-def}
\setcounter{equation}{0}

In order to compare the number of eigenvalues of (\ref{h1}) with the
Maslov index, we will use the Evans function to determine eigenvalues
based on the setup in \textsc{Alexander, Gardner \& Jones}~\cite{agj},
adapted to the symplectic setting in
\textsc{Bridges \& Derks}~\cite{bd-analyticity}, restricted to the
case of $\R^4$.

Consider the linear system of ODEs,
\begin{equation}\label{E1}
{\bf u}_x = {\bf A}(x,\lambda){\bf u}\,, \quad {\bf u}\in\R^{4}\,,
\end{equation}
where ${\bf A}(x,\lambda)={\bf J}^{-1}{\bf B}(x,\lambda)$ and
${\bf B}(x,\lambda)$ is symmetric and depends smoothly
on $x$ and $\lambda$.  In general $\lambda$ can be
complex but in this paper it will be restricted to be real.
Assume that ${\bf A}(x,\lambda)$ tends exponentially fast to a matrix
$\mathbf A_{\infty}(\lambda)$ when $x\to \pm\infty$.

Define the stable and unstable subspaces of ${\bf A}_\infty(\lambda)$ by
\[
{\rm E}^s_{\infty}(\lambda) := \{ {\bf u}\in\R^{4}\ :\ \lim_{x\to+\infty}{\rm e}^{{\bf A}_\infty(\lambda)x}{\bf u}=0 \}
\]
and
\[
{\rm E}^u_{\infty}(\lambda) := \{ {\bf u}\in\R^{4}\ :\ \lim_{x\to-\infty}{\rm e}^{{\bf A}_\infty(\lambda)x}{\bf u}=0 \}
\]
${\rm E}^s(\lambda)$ (${\rm E}^u_{\infty}(\lambda)$)
is the direct sum of the generalized eigenspaces associated with
the eigenvalues of ${\bf A}_\infty(\lambda)$ with negative (positive) real part.
The matrix ${\bf A}_\infty(\lambda)$ is said to be hyperbolic if $\R^{4}=
{\rm E}_\infty^u(\lambda)\oplus {\rm E}_\infty^s(\lambda)$; equivalently
if ${\bf A}_\infty(\lambda)$ has no purely imaginary eigenvalues.  Purely imaginary eigenvalues
of ${\bf A}_\infty(\lambda)$ are associated with the essential spectrum.
The essential spectrum is
\begin{equation}\label{sigma-ess}
\begin{array}{rcl}
\sigma_\ess&=&\{\lambda\in\C\mid \mathbf A_{\infty}(\lambda)\ \mbox{is not hyperbolic}\}\\[2mm]
&=&\{\ \lambda\in\C\ :\ {\rm det}[ {\bf A}_\infty(\lambda)-\ri\kappa{\bf I}]=0\quad
\mbox{for some $\kappa\in\R$}\ \}\,.
\end{array}
\end{equation}
We will assume throughout that $\lambda\notin\sigma_\ess$.
Then the Hamiltonian symmetry of ${\bf A}_\infty(\lambda)$ gives that
${\rm dim}\,{\rm E}_\infty^u(\lambda) = {\rm dim}\,{\rm E}_\infty^s(\lambda) = 2$.

Let $\bwedge^2(\R^{4})$ be the vector space of $2-$vectors in $\R^{4}$.  There is an induced system from (\ref{E1})
\begin{equation}\label{E2}
{\bf U}_x = {\bf A}^{(2)}(x,\lambda){\bf U}\,,\quad {\bf U}\in\bwedge^2(\R^{4})\,.
\end{equation}
Let $\sigma_+(\lambda)$ be the sum of the eigenvalues of $\mathbf A_\infty(\lambda)$
with positive real part, and let $\sigma_-(\lambda)$ be the sum of the eigenvalues
with negative real part.
Then there are solutions ${\bf U}^+(x,\lambda)$ and ${\bf U}^-(x,\lambda)$
of (\ref{E2}) with maximal decay as $x$ goes to $-\infty$ and $+\infty$ respectively
satisfying
\begin{equation}\label{Eplus-asymp}
{\rm lim}_{x\to-\infty}{\rm e}^{-\sigma_+(\lambda)x}{\bf U}^+(x,\lambda) = 
\zeta^+(\lambda)\in\bwedge^2(\R^{4})\,,
\end{equation}
and
\begin{equation}\label{Eminus-asymp}
{\rm lim}_{x\to+\infty}{\rm e}^{-\sigma_-(\lambda)x}{\bf U}^-(x,\lambda) = 
\zeta^-(\lambda)\in\bwedge^2(\R^{4})\,,
\end{equation}
where $\zeta^{\pm}(\lambda)$ are eigenvectors
\begin{equation}\label{eig-eqn}
{\bf A}^{(2)}_\infty(\lambda)\zeta^\pm(\lambda) = \sigma_\pm(\lambda)\zeta^\pm(\lambda)\,.
\end{equation}

The eigenvalues $\sigma_{\pm}(\lambda)$ are analytic functions of $\lambda$
and so the eigenvectors $\zeta^{\pm}(\lambda)$ can be chosen to be analytic
as well.

A value $\lambda\in\R\setminus\sigma_{ess}$ is called an eigenvalue if the stable solutions
${\bf U}^-(x,\lambda)$ and unstable solutions ${\bf U}^+(x,\lambda)$ have nontrivial
intersection.  Eigenvalues are detected by the \emph{Evans function} \cite{agj}
which is defined by
\begin{equation}\label{evans-def}
D(\lambda)\,\vol = {\bf U}^-(x,\lambda)\wedge{\bf U}^+(x,\lambda)\in\bwedge^{4}(\R^{4})\,.
\end{equation}
The Evans function is independent of $x$ and is an analytic function
of $\lambda$ \cite{agj}. 
Analyticity assures that the zeros of $D(\lambda)$ are isolated.  
In the definition
(\ref{evans-def}), the property of (\ref{h1}),
\begin{equation}\label{trace-prop}
{\rm Trace}({\bf A}(x,\lambda))=0\,,
\end{equation}
has been used.  This property
follows since ${\bf A}={\bf J}^{-1}{\bf B}$ with ${\bf J}$ skew-symmetric and
${\bf B}$ symmetric.

Because of the
Hamiltonian structure, the Evans function is invariant under exponential scaling
of the following form.  Let
\[
\widehat{\bf U}^{\pm}(x,\lambda) = {\rm e}^{-\sigma_{\pm}(\lambda)x}{\bf U}^{\pm}(x,\lambda)\,,
\]
Then the scaled functions satisfy
\[
\widehat{\bf U}^{\pm}_x = [{\bf A}^{(2)}(x,\lambda)-\sigma_{\pm}(\lambda){\bf I}]\widehat{\bf U}^\pm\,,
\]
but the Evans function becomes
\[
D(\lambda)\,\vol = {\rm e}^{(\sigma_-(\lambda)+\sigma_+(\lambda))x}
\widehat{\bf U}^-(x,\lambda)\wedge\widehat{\bf U}^+(x,\lambda)
=\widehat{\bf U}^-(x,\lambda)\wedge\widehat{\bf U}^+(x,\lambda)\,,
\]
since $\sigma_-(\lambda)+\sigma_+(\lambda) = {\rm Trace}({\bf A}_\infty(\lambda)) = 0$.

An explicit expression for the entries of ${\bf A}^{(2)}$ as a function
of the entries of ${\bf A}$, when $n=2$, is given in Appendix \ref{A2-kernel}.
It is natural to ask whether there is an induced symplectic structure on
$\bw$.  For example, can the induced system be written in the form
\[
{\bf J}^{(2)}{\bf U}_x = {\bf B}^{(2)}(x,\lambda){\bf U}\,,
\]
where ${\bf J}^{(2)}$ is the induced matrix from ${\bf J}$ on
$\bw$.  However, this is not the case.  The most significant obstacle is
the fact that ${\bf J}^{(2)}$ is not invertible.  The precise relation
between ${\bf J}^{(2)}$, ${\bf B}^{(2)}$ and ${\bf A}^{(2)}$ is given
in Appendix \ref{app_J2B2}.

\section{${\bf U}^\pm(x,\lambda)$ represent paths of Lagrangian planes}
\label{subsec-u-lagr}
\setcounter{equation}{0}

The paths of stable and unstable subspaces ${\bf U}^\pm(x,\lambda)$ (or their
scaled versions $\widehat{\bf U}^\pm(x,\lambda)$) are paths of Lagrangian
subspaces.

Let $\Phi(x,s,\lambda)$ be a fundamental solution matrix for
(\ref{h1}), that is,
\[
{\bf J}\Phi_x ={\bf B}(x,\lambda)\Phi\,,\quad \Phi(s,s,\lambda) = {\bf I}\,,
\]
and $\Phi(x,s,\lambda)\Phi(s,t,\lambda) = \Phi(x,t,\lambda)$.
Define the stable and unstable subspaces for each $x_0\in\R$ \cite{rs95},
\[
{\rm E}^s(x_0,\lambda) = \{ {\bf u}\in\R^4\ :\ \lim_{x\to+\infty}\Phi(x,x_0){\bf u} = 0 \}
= {\rm span}\left\{ {\rm col}({\bf U}^-(x_0,\lambda))\right\}\,,
\]
and
\[
{\rm E}^u(x_0,\lambda) = \{ {\bf u}\in\R^4\ :\ \lim_{x\to-\infty}\Phi(x,x_0){\bf u} = 0 \}
= {\rm span}\left\{ {\rm col}({\bf U}^+(x_0,\lambda))\right\}\,.
\]
Both subspaces define invariant vector bundles over $\R$.  This means that
\begin{equation}\label{invariant-vector-bundle}
{\rm E}^s(x,\lambda) = \Phi(x,s,\lambda) {\rm E}^s(s,\lambda)\quad \mbox{and}\quad 
{\rm E}^u(x,\lambda) = \Phi(x,s,\lambda) {\rm E}^u(s,\lambda)\,.
\end{equation}
Moreover
\[
\lim_{x\to+\infty}{\rm E}^s(x,\lambda) = {\rm E}_\infty^s(\lambda)\quad\mbox{and}\quad
\lim_{x\to-\infty}{\rm E}^u(x,\lambda) = {\rm E}_\infty^u(\lambda)\,.
\]

If ${\rm E}^u(x,\lambda)$ is Lagrangian for some $x$ then it is Lagrangian for all $x$.  This
observation is implicit in (\ref{invariant-vector-bundle}) 
but a direct proof can be given as follows.
When ${\bf u},{\bf v}\in{\rm E}^u(x,\lambda)$ are solutions of (\ref{h1}), 
\[
\frac{d\ }{dx}\langle {\bf J}{\bf u},{\bf v}\rangle
= \langle {\bf J}{\bf u}_x,{\bf v}\rangle +
 \langle {\bf J}{\bf u},{\bf v}_x\rangle=
 \langle {\bf B}(x,\lambda){\bf u},{\bf v}\rangle -
 \langle {\bf u},{\bf B}(x,\lambda){\bf v}\rangle= 0\,,
\]
using symmetry of ${\bf B}(x,\lambda)$.
Hence the value of $\langle {\bf J}{\bf u},{\bf v}\rangle$ is an invariant of 
(\ref{h1}) for any pair of vectors
${\bf u},{\bf v}\in\R^{4}$:
\[
\langle \mathbf J \mathbf u(x,\lambda), \mathbf v(x,\lambda) \rangle
= \langle \mathbf J \mathbf u(x_0,\lambda), \mathbf v(x_0,\lambda) \rangle\,,\quad
\forall x \,.
\]
But ${\bf u},{\bf v}\in{\rm E}^u(x,\lambda)$ and so
$\lim_{x\to-\infty} \mathbf u(x,\lambda)=
\lim_{x\to-\infty} \mathbf v(x,\lambda)=0$, therefore:
\[
\lim_{x\to-\infty} \langle\mathbf J \mathbf u(x,\lambda), \mathbf v(x,\lambda)\rangle=0\quad
\Rightarrow\quad \langle
\mathbf J \mathbf u(x_0,\lambda), \mathbf v(x_0,\lambda)\rangle = 0\,.
\]
Therefore, ${\rm E}^u(x,\lambda)$ is a Lagrangian subspace for any $x$. A similar proof confirms the result for ${\rm E}^s(x,\lambda)$.  Another proof is to use Montaldi's Theorem
\cite{montaldi} on Lagrangian planes in $\R^4$ and a sketch is given in Appendix
\ref{montaldi-lagrangian}.

\section{An example in $\R^2$}
\label{sec-example-r2}
\setcounter{equation}{0}

Before proceeding to the full definition and properties of the Maslov index for
paths of Lagrangian subspaces which are also solutions of (\ref{h1}) it will be useful to
consider the simplest possible context for the Maslov index, linear systems on $\R^2$.
To illustrate the role of $\lambda$, a stability problem for a reaction-diffusion equation is used.
It is a simplified version of the
class of nonlinear parabolic PDEs studied in \cite{bose}.

\begin{figure}
\begin{center}
\includegraphics[width=8cm]{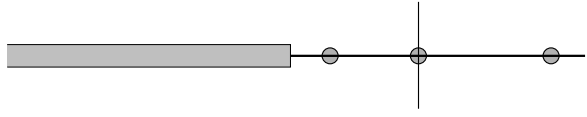}
\caption{Elements of the spectrum of $\mathscr{L}$ in the complex plane. 
The big dots are the elements of the point spectrum $\mathscr{L}$ and the dashed
interval represents the essential spectrum. The origin $\lambda=0$ is located at
the intersection of the two axes.}
\label{fig-spectrum}
\end{center}
\end{figure}

Consider the nonlinear parabolic PDE
\begin{equation}\label{pde-1}
\frac{\partial \phi}{\partial t} = \frac{\partial^2\phi}{\partial x^2} - \phi
+ \phi^2 \,,\quad x\in\R\,,
\end{equation}
for the scalar-valued function $\phi(x,t)$.  There is a basic 
steady solitary wave solution
\begin{equation}\label{basic-state}
\widehat\phi(x) = \mbox{$\frac{3}{2}$}\,{\rm sech}^2\left(\fr x\right)\,,
\end{equation}
which satisfies $\widehat\phi_{xx}-\widehat\phi +\widehat\phi^2=0$.
Linearizing (\ref{pde-1}) about the basic state $\widehat\phi$ and looking for
solutions proportional to ${\rm e}^{\lambda t}$ leads to the spectral problem
\begin{equation}\label{L-2}
\mathscr{L}\phi = \lambda\phi\,,\quad \mbox{with}\quad
\mathscr{L}\phi := \frac{d^2\phi}{dx^2}-\phi + 2\widehat\phi(x)\phi\,.
\end{equation}
The basic state (\ref{basic-state}) is said to be (spectrally) unstable
if any part of the spectrum of $\mathscr{L}$ is positive.
The spectrum of $\mathscr{L}$ can be explicitly constructed.
It consists of a branch of essential spectra and a point spectrum
\[
\sigma(\mathscr{L}) = \sigma_{ess}(\mathscr{L}) \cup \sigma_p(\mathscr{L})\,,
\]
with $\sigma_{ess}(\mathscr{L}) = \{ \lambda\in\R\ :\ \lambda\leq -1 \}$
and
$\sigma_p(\mathscr{L}) = \left\{ -\frac{3}{4}\,,\,0\,,\,\frac{5}{4}\right\}$.
The spectrum is illustrated in Figure \ref{fig-spectrum}.

The point spectrum can be verified by constructing the Evans function.
First reformulate
(\ref{L-2}) as a first-order system.  Let
\[
{\bf u}(x,\lambda) = \begin{pmatrix} \phi(x,\lambda)\\
\phi_x(x,\lambda) \end{pmatrix}\,,
\]
then
\begin{equation}\label{J-B-2}
{\bf J}{\bf u}_x = {\bf B}(x,\lambda){\bf u}\,,\quad{\bf u}\in\R^2\,,
\quad \lambda\in\R\,,
\end{equation}
with
\[
{\bf B}(x,\lambda) = \left[\begin{matrix}
\ \lambda+1-3{\rm sech}^2(\fr x) & \ 0 \\ \ 0\  & \ 1\ \end{matrix}\right]\,.
\]
The eigenvalues of ${\bf A}_{\infty}(\lambda)$ are real and hyperbolic when
$\lambda+1>0$.  In this formulation the stable (${\bf u}^-$)
and unstable (${\bf u}^+$)
subspaces are represented by
\[
{\bf u}^{\pm}(x,\lambda) = {\rm e}^{\pm\gamma s }\begin{pmatrix}
h^\pm\\
\fr(h^\pm_ s \pm\gamma h^\pm)\end{pmatrix}\,,
\]
where $ s =\fr x$, $\gamma=2\sqrt{\lambda+1}$,
\[
h^{\pm}( s ,\lambda) = \pm a_0 + a_1 {\rm tanh}( s )\pm a_2{\rm tanh}^2( s ) +
a_3\,{\rm tanh}^3( s )\,,
\]
and
\begin{equation}\label{a-def}
a_0 = \frac{\gamma}{15}(4-\gamma^2)\,a_3\,,\quad
a_1 = \frac{1}{5}(2\gamma^2-3)\,a_3\,,\quad
a_2 = -\gamma\, a_3\,,
\end{equation}
and $a_3$ is an arbitrary nonzero real number.
The Evans function is then
\[
D(\lambda) = {\bf u}^-(x,\lambda)\wedge{\bf u}^+(x,\lambda)\,.
\]
Evaluating at $x=0$, a straightforward calculation leads to
\[
D(\lambda) 
= -2\sqrt{\lambda+1}\,\left(\frac{2a_3}{15}\right)^2
\lambda\,(4\lambda+3)(4\lambda-5)\,.
\]
The zeros of $D(\lambda)$ are the eigenvalues, confirming the point spectrum
$\left\{-\frac{3}{4},0,+\frac{5}{4}\right\}$.

For linear Hamiltonian systems on $\R^2$ Lagrangian subspaces are just
one-dimensional subspaces.  The path of unstable subspaces
${\bf u}^+(x,\lambda)$ is used to define the Maslov index.  The natural one-dimensional 
subspace to choose for the reference space is ${\rm E}_\infty^s(\lambda)$,
\[
{\rm E}_\infty^s(\lambda)  ={\rm span}\left\{ \begin{pmatrix} 2 \\ -\gamma \end{pmatrix}
\right\}\,.
\]

Then, assume simple intersections between ${\rm E}_\infty^s(\lambda)$ and ${\bf u}^+(x,\lambda)$ -- which
can be confirmed a posteriori for the example (\ref{J-B-2}) -- and assume that
\[
\lim_{x\to\pm\infty} {\bf u}^+(x,\lambda)\bigcap {\rm E}_\infty^s(\lambda) = \{ 0\}\,.
\]
This latter assumption is equivalent to assuming that $\lambda$ is not an eigenvalue.
The Maslov index for this case is
\[
\textsf{Maslov}({\bf u}^+,{\rm E}_\infty^s) = \sum_{x_0} {\rm sign}\, \langle
{\bf Ju}^+_x,{\bf u}^+\rangle\,\vol\,,
\]
with $x_0$ the points at which ${\bf u}^+(x,\lambda)\cap{\rm E}^s_{\infty}$ is non-trivial, and the volume form
can be taken to be $\vol={\bf e}_1\wedge{\bf e}_2$.
This expression for the Maslov index is the one-dimensional version of
(\ref{Maslov-def}).

\begin{figure}
\begin{center}
\includegraphics[width=7cm]{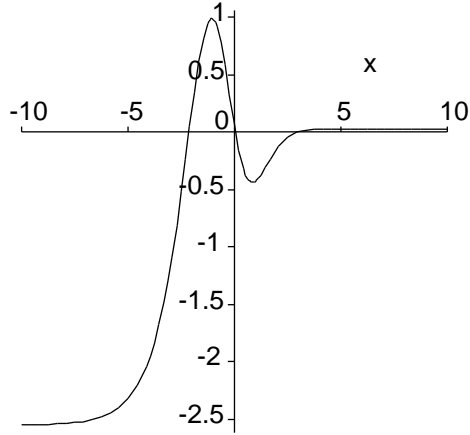}
\caption{Plot of $\xi^s(\lambda)\wedge{\bf u}^+(x,\lambda)$ as a function of $x$ for the case $\lambda=-0.8$.}
\label{int-lambda-58}
\end{center}
\end{figure}

The path of unstable subspaces is
\begin{equation}\label{W-u}
{\bf u}^{+}(x,\lambda) = \fr {\rm e}^{\gamma s }\begin{pmatrix}
2h^+\\
h^+_ s +\gamma h^+\end{pmatrix}\,.
\end{equation}
The intersection form in this case is
\[
\begin{array}{rcl}
\Gamma({\bf u}^+,{\rm E}_\infty^s,x_0) &=& \langle
{\bf Ju}^+_x,{\bf u}^+\rangle\bigg|_{x=x_0}\,\vol\,,\\[2mm]
&=& (-u_1^+\dot u_2^++u_2^+\dot u_1^+)\bigg|_{x=x_0}\,\vol\\[2mm]
&=& \left[
( (u_2^+)^2 - \lambda - 1 + 12\,{\rm sech}^2 s )(u_1^+)^2 \right]\bigg|_{x=x_0}\vol\,.
\end{array}
\]
However, at a point $x_0$ where ${\bf u}^+$ intersects ${\rm E}_\lambda^s$,
$u_2^+ = -\fr \gamma u_1^+$ and so
\[
\Gamma({\bf u}^+,{\rm E}_\infty^s,x_0)= 12\,{\rm sech}^2\fr x_0 \,(u_1^+)^2\,\vol\,.
\]
Hence $\Gamma({\bf u}^+,{\rm E_\infty}^s)>0$ at each intersection, and the Maslov
index is just the sum of the intersections.
An intersection occurs when
\[
\xi^s\wedge{\bf u}^+=0\,,\quad \mbox{where}\quad
 \xi^s = \begin{pmatrix} 2 \\ -\gamma\end{pmatrix}\,.
\]
Now
\[
\xi^s\wedge{\bf u}^+= 2u_2^+ + \gamma u_1^+\,.
\]
The factor ${\rm e}^{\gamma s}$ is not important and so can be divided out,
giving
\[
\xi^s\wedge{\bf u}^+ \sim \frac{dh^+}{d s} + 2\gamma h^+\,.
\]
This function has $0,1,2$ or $3$ zeros depending on the value of $\lambda$.
Each zero corresponds to an intersection between the unstable subspace
with ${\rm E}_\infty^s(\lambda)$.  The function $\xi^s\wedge{\bf u}^+$ is illustrated
in Figure \ref{int-lambda-58} for the case $\lambda=-0.8$ where $\xi^s\wedge{\bf u}^+ $
has three zeros indicating
three intersections.  A summary of the Maslov index in each region is tabulated below.  
\vspace{.25cm}

\begin{center}
\begin{tabular}{|c|c|c|c|c|}
\hline
\phantom{$\bigg|$}
$\lambda$ & $-1<\lambda<-\frac{3}{4}$  & $-\frac{3}{4}<\lambda<0$ & $0<\lambda<\frac{5}{4}$ & $\lambda>\frac{5}{4}$ \\
\hline
\phantom{$\bigg|$}$\textsf{Maslov}({\bf u}^+,{\rm E}_\infty^s)$ & $3$ & $2$ & $1$ & $0$ \\
\hline
\end{tabular}
\end{center}
\vspace{.25cm}

Note that the Maslov index jumps by one at each eigenvalue.
Let $\lambda_0$ be any fixed real value of
$\lambda$ such that $\lambda_0>-1$ and $\lambda_0$ is
not an eigenvalue, then the value of the
Maslov index equals the number of eigenvalues of
$\mathscr{L}$ in the set $\lambda>\lambda_0$. 

This link between the Maslov index and the number of eigenvalues
is not a coincidence specific to this operator but a quite general
fact for systems such that
$\partial_\lambda \mathbf B(x,\lambda)$ is a semi-definite matrix  (see \cite{chardard-thesis} for example). 
In the case where $\partial_\lambda \mathbf B(x,\lambda)$ is a not semi-definite, we give a counterexample
in section \ref{sec-10}.

\subsection{The Maslov angle in $\bwedge^1(\R^2)$}
\label{subsec-2d-angle}

Another way to count intersections between the path ${\bf u}^+$ and some reference
plane is to use the Maslov angle (\ref{angle-def-1}).
In this case the angle $\kappa(x,\lambda)$ is
just the angle determined by a polar representation of ${\bf u}^+$
\[
{\rm e}^{\ri\kappa(x,\lambda)}
:= \frac{u_1^+(x,\lambda) - \ri u_2^+(x,\lambda)}{u_1^+(x,\lambda) + \ri u_2^+(x,\lambda)}\,.
\]
As $x\to\pm\infty$
\[
\lim_{x\to\pm\infty}{\rm e}^{\ri\kappa(x,\lambda)} = \frac{2-\ri\gamma}{2+\ri\gamma}\,.
\]
The Maslov index is then the count of the number of times that $\kappa$ crosses
some reference angle, such as the angle associated with the stable subspace.
This version of the Maslov index
is equivalent to the definition based on intersection index.

\section{The Maslov index for paths}
\label{sec-paths}
\setcounter{equation}{0}

In this section we look at some of the properties of paths of Lagrangian subspaces that
are also solutions of (\ref{h1}), and bring in the exterior algebra representation.

Here and throughout, let ${\bf e}_1,\ldots,{\bf e}_{4}$
be the canonical basis for $\R^{4}$. The symplectic form 
in standard form is then
\begin{equation}\label{omega-def}
\symp = {\bf e}_1\wedge{\bf e}_{3} +  {\bf e}_2\wedge{\bf e}_{4}\,.
\end{equation}

The symplectic form is related to the symplectic operator ${\bf J}$ in
(\ref{J-def}) by
\[
\langle  {\bf J}{\bf u} , {\bf v}\rangle = 
\lbk \symp,{\bf u}\wedge{\bf v}\rbk_2 \,,\quad \forall\ {\bf u},{\bf v}
\in\R^{2n}\,.
\]

Here and throughout $\langle\cdot,\cdot\rangle$ is the standard
inner product on $\R^{4}$ and $\lbk\cdot,\cdot\rbk_k$ is the induced
inner product on $\bwedge^k(\R^{4})$. We will sometimes write $\langle\cdot,
\cdot\rangle_d$ on $\R^d$ when the associated dimension is not clear.
The equivalence between the induced inner product $\lbk\cdot,\cdot\rbk_k$
and $\langle \cdot,\cdot\rangle_d$ when $d={\rm dim}\left(
\bigwedge^k(\R^{2n}\right)$
is established in Appendix \ref{app-ip}.  The basic properties of symplectic
exterior algebra can be found in \cite{klr}.

The set of all Lagrangian subspaces of $\R^{4}$ associated with 
the standard symplectic operator $\symp$ will be denoted by $\Lambda(2)$.
$\Lambda(2)$ is a manifold of dimension $3$.  It is 
a submanifold of ${\rm G}_2(\R^{4})$, the Grassmannian of all
$2-$dimensional subspaces of $\R^{4}$ \cite{klr}.

A path of unstable subspaces will be represented by
${\bf U}^+(x,\lambda)\in\bwedge^2(\R^{4})$, which satisfies the equation (\ref{E2})
and (\ref{Eplus-asymp}).  When $\lambda$ is not an eigenvalue then 
${\bf U}^+(x,\lambda)\to{\rm E}_\infty^u(\lambda)$ as $x\to\pm\infty$ and so in
$\mathbb{P}({\rm E}^u(x,\lambda))$, projective space based on ${\rm E}^u(x,\lambda)$,
a loop of Lagrangian subspaces is obtained.  The Maslov index is then taken to be
the sum of the weighted intersections of ${\bf U}^+(x,\lambda)$ with ${\rm E}_\infty^s(\lambda)$.

${\rm E}_\infty^s(\lambda)$ is represented by the $2-$form $\zeta^-(\lambda)$ defined in
(\ref{eig-eqn}).  An analytic basis can always be constructed for ${\rm E}_\infty^s(\lambda)$
\cite{bd-analyticity}.  Denote this basis by
\[
{\rm E}_\infty^s(\lambda) = {\rm span}\{\xi_1^s(\lambda),\xi_2^s(\lambda)\}\,.
\]
 
Consider the $3-$form
\[
{\bf U}^+(x,\lambda)\wedge( \alpha_1\xi_1^s(\lambda) + \alpha_2\xi_2^s(\lambda) ) \,.
\]
If for any fixed $x$ and $\lambda$,
\[
{\bf U}^+(x,\lambda)\wedge( \alpha_1\xi_1^s(\lambda) + \alpha_2\xi_2^s(\lambda) ) = 0 
\quad\Rightarrow\quad \alpha = 0 \,,
\]
then we say that ${\rm E}_\infty^s(\lambda)$ is transverse to ${\rm E}^u(x,\lambda) := {\rm image}(
{\bf U}^+(x,\lambda))$ at that value of $(x,\lambda)$.

We say that ${\rm E}_\infty^s(\lambda)$ and ${\rm E}^u(x,\lambda)$ have a simple intersection
(or regular intersection) if the intersection is one dimensional;
there exists $\alpha\in\rp^{1}$ such that
\begin{equation}\label{reg-inters}
{\bf U}(x,\lambda)\wedge\xi=0\quad\mbox{with}\quad \xi =\alpha_1\xi_1^s(\lambda) 
+ \alpha_2\xi_2^s(\lambda) \,.
\end{equation}
We will assume that non-trivial intersections are regular.  It is proved in \cite{arnold1} that
non-regular intersections can be eliminated by perturbation.  This property can also be
proved using the homotopy equivalence property of the Maslov index, and a nice proof of this
is given in \S3.4 of \cite{m-g03}.

Now, let ${\bf u}^+_1(x,\lambda),{\bf u}^+_2(x,\lambda)$ be a basis for ${\rm E}^u(x,\lambda)$
such that
\[
{\bf U}^+(x,\lambda) = {\bf u}^+_1(x,\lambda)\wedge {\bf u}^+_2(x,\lambda)\,.
\]
Then (\ref{reg-inters}) implies that there exists $(\beta_1,\beta_2)\in\R^2\setminus\{0\}$ 
such that
\begin{equation}\label{beta-def}
\xi = \beta_1{\bf u}^+_1(x,\lambda) + \beta_2{\bf u}^+_2(x,\lambda) := {\bf Z}^+(x,\lambda)
\beta\,.
\end{equation}
The Maslov index is a count of how many times the path ${\bf U}^+(x,\lambda)$
crosses the reference space ${\rm E}^s_\infty(\lambda)$, weighted by the intersection form.
Starting with (\ref{Gamma-def}), and using the equation (\ref{h1}) and
the equivalence (\ref{beta-def}) it is
\begin{equation}\label{g-b}
\Gamma(\Lambda,V,x_0) = \langle {\bf B}(x_0,\lambda)\xi,\xi\rangle\,\vol\,.
\end{equation}
Note that the $x-$dependent path ${\bf U}^+(x,\lambda)$ drops out.  For fixed $\lambda$,
once we have found a point $x_0$ corresponding to a regular intersection, the intersection form
can be evaluated using eigenvectors spanning ${\rm E}_\infty^s(\lambda)$ and the known symmetric
matrix ${\bf B}(x_0,\lambda)$.  
However there is the hidden calculation of determining $\alpha\in\rp^1$ and an algorithm for
this is developed in \S\ref{sec-numer-algor-1}.

Suppose $\lambda$ is not an eigenvalue.  Then the Maslov index of the path ${\bf U}^+$
relative to ${\rm E}_\infty^s$ is
\begin{equation}\label{maslov-def}
\textsf{Maslov}({\bf U}^+,{\rm E}_\infty^s) = \sum_{x_0} {\rm sign}\, \langle \xi,{\bf B}(x_0,\lambda)\xi\rangle\,\vol\,,
\end{equation}
where the sum is over all points $x_0$ of intersection in the interval $-\infty<x_0<+\infty$.

This expression will serve as a definition for the Maslov 
index of a solitary wave
at $\lambda$:

\begin{equation}
  \label{eq:def}
\mathsf{Maslov}(\lambda)=\textsf{Maslov}({\bf U}^+(\cdot,\lambda),\mathrm E^s_{\infty}(\lambda)).
\end{equation}
\vspace{.1cm}

\begin{prop} Suppose $\lambda$ is not an eigenvalue.  Then under the given hypotheses on
(\ref{h1}) the Maslov index of a solitary wave is finite.
\end{prop}

\noindent{\bf Proof.}
Introduce a metric $\textsf{dist}(\cdot,\cdot)$
on the manifold of $2-$dimensional subspaces of $\R^{4}$.  For example this
can be the standard metric on the Grassmannian ${\rm G}_2(\R^{4})$ \cite{wong}.

Let $\lambda\in\R\setminus\sigma_{ess}$ which is not an eigenvalue.
Since ${\rm E}^s_\infty(\lambda)$ and ${\rm E}^u_\infty(\lambda)$
are transverse,
by a suitable scaling of the stable subspace ${\rm E}^s_\infty(\lambda)$ we can take
\[
\textsf{dist}(({\rm E}^s_\infty(\lambda)),
{\rm E}^u_\infty(\lambda)) > 1 \,,
\] 
where  $({\rm E}^s_\infty(\lambda))$ is the closed set
of planes which are not transverse to $({\rm E}^s_\infty(\lambda))$.
To simplify notation, the argument in $\textsf{dist}(\cdot,\cdot)$ should
be interpreted as the representation of ${\rm E}^{s,u}$ on 
the Lagrangian Grassmannian.

When $\lambda$ is not an eigenvalue we have that
${\rm E}^u(x,\lambda) \to {\rm E}^u_\infty(\lambda)$ as $x\to\pm\infty$.  Therefore,
given any $\epsilon>0$, there exists $\delta>0$ such that 
\[
\textsf{dist}({\rm E}^u(x,\lambda),{\rm E}^u_\infty(\lambda)) < \varepsilon\quad\mbox{for}\quad
|x|>\frac{1}{\delta}\,.
\]
Now use the triangle inequality 
\[
1<\textsf{dist}(({\rm E}^s_\infty(\lambda)),{\rm E}^u_\infty(\lambda))\leq
\textsf{dist}(({\rm E}^s_\infty(\lambda)),{\rm E}^u(x,\lambda)) +
\textsf{dist}({\rm E}^u(x,\lambda),{\rm E}^u_\infty(\lambda))\,,
\]
or
\[
\textsf{dist}({\rm E}^u(x,\lambda),({\rm E}^s_\infty(\lambda))) > 1-\varepsilon\quad
\mbox{for}\quad |x|>\frac{1}{\delta}\,.
\]
Hence there exists $x_*>0$ such that for $|x|>x_*$ 
${\rm E}^u(x,\lambda)$ and ${\rm E}^s_\infty(\lambda)$ are transverse.
Intersections are therefore limited to the finite interval $-x_*<x<x_*$.
Since intersections are generically isolated their number is finite.
$\hfill\blacksquare$

\section{The Maslov index on $\bwedge^2(\R^4)$}
\label{sec-r4}
\setcounter{equation}{0}

The vector space $\bigwedge^2(\R^4)$ is six-dimensional, and the 
orthonormal basis induced from the basis of $\R^4$ is
\begin{equation}\label{E-basis}
\begin{array}{rcl}
{\bf E}_1 &=& {\bf e}_1\wedge{\bf e}_2\,,\quad {\bf E}_2 = {\bf e}_1\wedge{\bf e}_3\,,\quad
{\bf E}_3 = {\bf e}_1\wedge{\bf e}_4\,,\\[2mm]
{\bf E}_4 &=& {\bf e}_2\wedge{\bf e}_3\,,\quad {\bf E}_5 = {\bf e}_2\wedge{\bf e}_4\,,\quad
{\bf E}_6 = {\bf e}_3\wedge{\bf e}_4\,.
\end{array}
\end{equation}
Any ${\bf U}\in\bwedge^2(\R^4)$ can be represented in the form
\begin{equation}\label{U-24}
{\bf U} = \sum_{j=1}^6 U_j\,{\bf E}_j\,.
\end{equation}
The Grassmannian $G_2(\R^4)$ is a subset of $\bwedge^2(\R^4)$ defined by
\begin{equation}\label{g24}
0={\bf U}\wedge{\bf U} = I_1\,\vol\,,\quad I_1 := U_1U_6-U_2U_5+U_3U_4\,.
\end{equation}
The Lagrangian-Grassmannian is the subset defined by
\begin{equation}\label{Lambda2}
0 = \symp\wedge{\bf U} = I_2\,\vol\,,\quad I_2 := U_2+U_5 \,.
\end{equation}
The Lagrangian-Grassmannian $\Lambda(2)$ is the three dimensional submanifold of
$\P\left(\bwedge^2(\R^4)\right)$ defined by $I_1=I_2=0$.

Let $V\in\bwedge^2(\R^4)$ be a fixed Lagrangian plane.  Then 
\[
\Lambda^1(2) = \{\ {\bf U}\in \bwedge^2(\R^4)\cap\Lambda(2)\ : \ {\bf U}\wedge V=0\ \}\,,
\]
is a codimension one submanifold of $\Lambda(2)$ \cite{arnold1}.
We have a sequence of manifolds
\vspace{.25cm}

\begin{center}
\begin{tabular}{|c|c|c|c|c|c|}
\hline
\phantom{$\bigg|$}
\textsf{Manifold} & $\bwedge^2(\R^4)$ & $\rp^5$ & ${\rm G}_2(\R^4)$ & $\Lambda(2)$ & $\Lambda^1(2)$ \\
\hline
\phantom{$\bigg|$}\textsf{Dimension} & $6$ & $5$ & $4$ & $3$ & $2$ \\
\hline
\end{tabular}
\end{center} 
\vspace{.25cm}

In this table $\rp^5$ represents $\P\left(\bwedge^2(\R^4)\right)$. 

Consider the class of linear Hamiltonian systems (\ref{h1})
with ${\bf B}(x,\lambda)$ satisfying the asymptotic properties (\ref{A-def})-(\ref{A_infty}).

\begin{prop}\label{prop-5-im}
$\Lambda(2)$ is an invariant manifold of (\ref{E2}). \end{prop} 

\noindent{\bf Proof.} 
\[
\begin{array}{rcl}
\frac{d\ }{dx}{\bf U}\wedge{\bf U} &=& {\bf U}_x\wedge{\bf U} + {\bf U}\wedge{\bf U}_x \\
&=& {\bf A}^{(2)}{\bf U}\wedge{\bf U} + {\bf U}\wedge{\bf A}^{(2)}{\bf U} \\
&=& {\rm Trace}({\bf A})\,{\bf U}\wedge{\bf U}\\
&=& 0 
\end{array}
\]
since ${\rm Trace}({\bf A})=0$, also using the property \cite{allen}
\[
{\bf A}^{(2)}{\bf U}\wedge{\bf U} + {\bf U}\wedge{\bf A}^{(2)}{\bf U} = {\rm Trace}({\bf A}){\bf U}\wedge{\bf U}\,.
\]
This proves that ${\bf U}\wedge{\bf U}$ is a constant along solutions.  Similarly,
\[
\begin{array}{rcl}
\frac{d\ }{dx}\symp\wedge{\bf U} &=& \symp\wedge{\bf U}_x \\
&=& \symp\wedge{\bf A}^{(2)}{\bf U} \\
&=& \symp\wedge{\bf A}^{(2)}{\bf U} +
{\bf A}^{(2)}\symp\wedge{\bf U} - {\bf A}^{(2)}\symp\wedge{\bf U}\\
&=& {\rm Trace}({\bf A})\symp\wedge{\bf U} - {\bf A}^{(2)}\symp\wedge{\bf U}\\
&=& - {\bf A}^{(2)}\symp\wedge{\bf U}\\
&=& 0\,,
\end{array}
\]
since $\symp$ is in the kernel of ${\bf A}^{(2)}$, a property which 
is proved in Appendix \ref{A2-kernel}.  This proves that
${\bf U}\wedge{\bf U}$ and $\symp\wedge{\bf U}$ are constant along solutions.  Hence the
special case ${\bf U}\wedge{\bf U}=\symp\wedge{\bf U}=0$ completes the proof.
$\hfill\blacksquare$
\vspace{.15cm}

Let $ V={\rm span}\{\xi_1,\xi_2\}$ be a fixed Lagrangian plane; that
is $\xi_1$ and $\xi_2$ are linearly independent and $\langle{\bf J}
\xi_1,\xi_2\rangle=0$.  The reference subspace $V$ is represented by the form
\[
{\bf V} = \xi_1\wedge\xi_2\,.
\]
The intersection between a path of Lagrangian subspaces ${\bf U}(x,\lambda)$ and 
$\mathbf V$ can be described as follows.  For each fixed $\lambda$ define
\[
\mathscr{A}(x) = \left\{\alpha\in\R^2 \ :\ 
{\bf U}(x,\lambda)\wedge( \alpha_1\xi_1 + \alpha_2\xi_2 ) =0 \right\}\,.
\]
Then there are three cases
\begin{itemize}
\item If $\mathscr A(x_0)= \{0\}$ then ${\bf U}$ and $V$ are transverse at $x=x_0$.
\item If $\mathscr A(x_0)$ is one dimensional then ${\bf U}$ and $V$ intersect in a one-dimensional subspace 
at $x=x_0$ (the case of regular crossing).
\item If $\mathscr A(x_0)=\R^2$ then ${\bf U}$ and $V$ intersect in a two-dimensional subspace at $x=x_0$
(this case is sometimes referred to as an intersection between ${\bf U}$ and the \emph{vertex} of $V$).
\end{itemize}

The Maslov index for a path of Lagrangian subspaces is 
given by (\ref{maslov-def}).  However, in the case of $\bwedge^2(\R^4)$
a new representation of the intersection form 
can be obtained.  Suppose that a regular crossing occurs
\[
{\bf U}\cap{\bf V} = {\rm span}\{\xi\}\,,
\]
at $x=x_0$, then the crossing form is
\begin{equation}\label{ifw}
\Gamma({\bf U},V,x_0) = \symp\wedge \xi \wedge {\bf A}\xi \,,
\end{equation}
To verify this formula, note that
\[
\symp\wedge{\bf a}\wedge{\bf J}^{-1}{\bf c} = \langle{\bf a},{\bf c}\rangle\,\vol\,,\quad
\mbox{for any}\ {\bf a},{\bf c}\in\R^4\,.
\]
Hence
\[
\symp\wedge\xi\wedge{\bf A}\xi = \symp\wedge\xi\wedge{\bf J}^{-1}{\bf B}\xi =
\langle\xi,{\bf B}\xi\rangle\,\vol\,,
\]
recovering the expression in (\ref{g-b}).  There is an interesting
geometric interpretation of (\ref{ifw}).  At a regular intersection
the two-plane $\xi\wedge{\bf A}\xi$ is \emph{not} a Lagrangian plane.
It is an element of ${\rm G}_2(\R^4)$ but not an element of $\Lambda(2)$.
Since $\Lambda(2)$ is a codimension one submanifold of ${\rm G}_2(\R^4)$,
the sign of $\symp\wedge\xi\wedge{\bf A}\xi$ determines which side
of $\Lambda(2)$ in ${\rm G}_2(\R^4)$ it lies.
See Appendix \ref{app-plus-minus} for further discussion of this case.

Hence the Maslov index of a path ${\bf U}(x,\lambda)$ relative to $V$ is
\[
\textsf{Maslov}({\bf U},V) = \sum_{x_0} {\rm sign} (\symp\wedge\xi\wedge{\bf A}\xi)\,,
\]
where the sum is over all interior intersections.

\section{The Maslov angle on $\bwedge^2(\R^4)$}
\label{sec-angle-wedge}
\setcounter{equation}{0}

In $\R^2$, the Maslov angle is just the angle associated with the
polar representation of a vector in $\R^2$ as shown in \S\ref{subsec-2d-angle}.
For a Lagrangian frame of the form (\ref{Z-def})
the Maslov angle is defined as in (\ref{angle-def}) and (\ref{angle-def-1}).
In this section a new formula for the Maslov angle is given for the exterior
algebra representation of a Lagrangian plane.  Here, the result for 4D phase
space is given and in Part 2 \cite{cdb-part2}, the general result for
$2n-$dimensional phase space is given.

A Lagrangian frame can be partitioned into two $4\times 2$ blocks as in
(\ref{Z-def}) and it can also be represented in terms of its columns:

\begin{equation}
  \label{eq:Z_eq}
 {\bf Z} =
\begin{pmatrix}
  \mathbf X\\
  \mathbf Y
\end{pmatrix}
=[\, {\bf z}_1\,|\,{\bf z}_2\,]\,,\quad\mbox{with}\quad
\langle{\bf J}{\bf z}_1,{\bf z}_2\rangle = 0\,. 
\end{equation}
Denote the exterior algebra representation of the Lagrangian plane by

\begin{equation}
  \label{eq:Z_eq_wedge}
{\bf U} =   {\bf z}_1\wedge {\bf z}_2\,.  
\end{equation}

\begin{prop}
\label{prop-C-def}
There exists a $2-$form ${\bf C}$,
\[
{\bf C} = {\bf C}_1 + \ri {\bf C}_2\,,\quad\mbox{with}\quad
{\bf C}_1,{\bf C}_2\in\bwedge^2(\R^4)\,,
\]
such that
\begin{equation}\label{angle-formula}
{\rm det}[ {\bf X}-\ri{\bf Y}]\vol = {\bf C}\wedge{\bf U}\,.
\end{equation}
\end{prop}

It follows from this proposition that there exists a scalar complex-valued
function
$K$ such that:

 \begin{equation}
   \label{eq:K_def}
   {\bf C}\wedge{\bf U}=K(\mathbf U)\vol.
 \end{equation}  

A formula for the Maslov angle $\kappa$ 
(the real number defined up to an addition by a multiple 
of $2\pi$ by $e^{\i\kappa}=\det( ({\bf X}-\ri{\bf Y})({\bf X}+\ri{\bf Y})^{-1})$.) is then immediate:

\begin{prop}
\[
{\rm e}^{\ri\kappa} = {K(\mathbf U)}/{\overline{K(\mathbf U)}}\,.
\]
\end{prop}

It remains to prove Proposition \ref{prop-C-def}.
The proof is by explicit construction.  Let
\[
{\bf c}_j = {\bf e}_j - \ri{\bf J}{\bf e}_{j}\,,\quad j=1,2\,.
\]
Then 
\begin{equation}\label{xy-c}
{\bf X} - \ri  {\bf Y}
=\begin{pmatrix} \phantom{-\ri}{\bf I} \\ -\ri{\bf I}\end{pmatrix}^T
\begin{pmatrix} {\bf X} \\ {\bf Y}\end{pmatrix} 
= [{\bf c}_1\,|\,{\bf c}_2]^T
[\, {\bf z}_1\,|\,{\bf z}_2\,]
=\begin{pmatrix}\langle{\bf c}_1,{\bf z}_1\rangle & \langle{\bf c}_1,{\bf z}_2\rangle \\
\langle{\bf c}_2,{\bf z}_1\rangle & \langle{\bf c}_2,{\bf z}_2\rangle \end{pmatrix}
\,,
\end{equation}
and so, using the induced inner
product\footnote{A real inner product is used throughout the paper.
Complexification is used so rarely, a Hermitian inner product is not
necessary.  One just needs to keep track of the complex conjugations.}
on $\bw$ (see Appendix \ref{app-ip})
\[
{\rm det}[ {\bf X} - \ri  {\bf Y} ]\vol =
{\rm det}\left[ \begin{matrix}\langle{\bf c}_1,{\bf z}_1\rangle & 
\langle{\bf c}_1,{\bf z}_2\rangle \\
\langle{\bf c}_2,{\bf z}_1\rangle & \langle{\bf c}_2,{\bf z}_2\rangle
\end{matrix}\right]\vol
= \lbk {\bf c}_1\wedge{\bf c}_2,{\bf U}\rbk_2\vol \,.
\]
This gives a formula for $K$,
\[
K(\mathbf U) =_{def} \lbk {\bf c}_1\wedge{\bf c}_2,{\bf U}\rbk_2\,.
\]
It is not necessary to give an expression for ${\bf C}$ since in computation
it is $K$ that is needed.  However, for completeness it is given.  
Let ${\bf C}$ be an $2-$form satisfying
\begin{equation}\label{C-def}
{\bf c}_1\wedge{\bf c}_2 \wedge \overline{\bf C} = \lbk
{\bf c}_1\wedge{\bf c}_2,\overline{{\bf c}_1\wedge{\bf c}_2}\rbk_2\,\vol\,.
\end{equation}
Then 
\[
{\rm det}[ {\bf X} - \ri  {\bf Y} ]\vol = {\bf C}\wedge{\bf U}\,.
\]
The $2-$form ${\bf C}$ is in fact the Hodge star of
${\bf c}_1\wedge{\bf c}_2$ although the details of
that characterization are not needed. 

Now compute the formula in coordinates on $\bw$.
On $\R^4$ with the standard basis,
\[
\begin{array}{rcl}
{\bf c}_1\wedge{\bf c}_2 &=& ({\bf e}_1 - \ri{\bf J}{\bf e}_{1})\wedge({\bf e}_2 - \ri{\bf J}{\bf e}_{2}) \\
&=& ({\bf e}_1 - \ri{\bf e}_{3})\wedge({\bf e}_2 - \ri{\bf e}_{4}) \\
&=& {\bf e}_1\wedge{\bf e}_2 -\ri {\bf e}_1\wedge{\bf e}_4 +\ri{\bf e}_2\wedge{\bf e}_3 -
{\bf e}_3\wedge{\bf e}_4 \,,
\end{array}
\]
and so, if ${\bf U} = \sum_{j=1}^6U_j{\bf E}_j$, with ${\bf E}_1,\ldots,{\bf E}_6$
the standard basis on $\bwedge^2(\R^4)$,
\[
K(\mathbf U) = \lbk{\bf c}_1\wedge{\bf c}_2,{\bf U}\rbk_2 = U_1 - \ri U_3 + \ri U_4 - U_6 \,,
\]
and so the expression for the Maslov angle is
\begin{equation}\label{kappa-wedge2}
{\rm e}^{i\kappa} = \frac{ U_1-U_6 -\ri U_3 + \ri U_4}{ U_1-U_6 +\ri U_3 - \ri U_4}\,.
\end{equation}
This expression is equivalent to the formula derived in equation (22) of \cite{bose}.

The two-form ${\bf C}$ in (\ref{angle-formula}) is computed to be
\[
{\bf C} = -{\bf e}_1\wedge{\bf e}_2 + {\bf e}_3\wedge{\bf e}_4 +\ri({\bf e}_1\wedge{\bf e}_4-
{\bf e}_2\wedge{\bf e}_3 )\,.
\]

\section{Further decomposition of the Maslov angle using the eigenvalues
of a unitary matrix}
\label{app-eigs-Q}
\setcounter{equation}{0}

Let ${\bf Z}\in\R^{4\times2}$ be a Lagrangian frame on $\R^4$
of the form (\ref{eq:Z_eq}) and $\mathbf U$ defined as in (\ref{eq:Z_eq_wedge}).  Then the matrix
\[
{\bf Q} = ({\bf X}-\ri{\bf Y})({\bf X}+\ri{\bf Y})^{-1}\,,
\]
is a unitary and symmetric (but not Hermitian) matrix.

The Maslov angle for a closed path  (\ref{angle-def}) is defined using the determinant
of ${\bf Q}$.  However, ${\bf Q}$ has $2$ eigenvalues
of unit modulus.  Denote these eigenvalues by ${\rm e}^{\ri\kappa_j}$,
$j=1,2$ with $\kappa_j$ real.  Then
\[
{\rm e}^{\ri\kappa} = {\rm e}^{\ri\kappa_1}\,{\rm e}^{\ri\kappa_2}\quad
\Rightarrow\quad \kappa=\kappa_1+\kappa_2\ \mbox{(mod $2\pi$)}\,.
\]
These eigenvalues are independent of the choice of $\mathbf Z\in \R^{4\times 2}$ as
a representative of a Lagrangian space:  choosing another
representation leads to similar matrix.

These eigenvalues can also be used to give another formula for the
sign of each intersection.  Fix the reference angle to be $0$ (mod $2\pi$).
There is a $1$-dimensional intersection at $x_0$ with a reference plane if
and only if there exists $\re^{\i\kappa_r(x_0)}=1$ with $r=1$ \emph{or} $r=2$.
There is a $2$-dimensional intersection at $x_0$ with the reference plane if
and only if $\re^{\i\kappa_r(x_0)}=1$ for $r=1$ \emph{and} $r=2$.
If $\re^{\i\kappa_r(x_0)}\neq1$ for $r=1$ and $r=2$ then the intersection
is \emph{transverse}.

When the intersection is regular, the sign of the intersection is given by:
\[
\lim_{x\to x_0^+} \#\{r\in S |  \kappa_r(x_0)\in\ (0,\pi)+2\pi \Z\}
- \#\{r\in S |  \kappa_r(x_0)\in\ (-\pi,0)+2\pi \Z\}\,.
\]
Thus, it is possible to determine the Maslov index, defined
with intersections by simply tracking the crossings of the
angles $\kappa_i$ with $2\pi\Z$.

\subsection{The angles $\kappa_j$ in the exterior algebra representation}

The two angles $\kappa_1,\kappa_2$ satisfy
\begin{equation}\label{Q-eigs}
{\rm det}( \mu{\bf I}-{\bf Q}) = \mu^2 - {\rm Trace}({\bf Q})\mu +
{\rm det}({\bf Q})=0\,,
\end{equation}
with $\mu$ of unit modulus and $\mu_r={\rm e}^{\i\kappa_r}$, $r=1,2$.
Both the trace and determinant can be expressed in terms of the exterior
algebra representation.

For the determinant, as shown in \S\ref{sec-angle-wedge},

\[
{\rm det}[{\bf X}-\ri{\bf Y}] = K(\mathbf U) := \lbk{\bf c}_1\wedge{\bf c}_2,
{\bf U}\rbk_2= U_1 - \ri U_3 + \ri U_4 - U_6
\,,
\]
for $\bwedge^2(\R^4)\ni{\bf U}= \sum_{j=1}^6U_j{\bf E}_j$.  Hence
\begin{equation}\label{det-Q}
{\rm det}({\bf Q}) = K(\mathbf U)/\overline{K(\mathbf U)}\,.
\end{equation}
It remains to express the Trace of {\bf Q} in terms of the exterior algebra
representation.

\begin{prop}
\begin{equation}\label{trace-Q}
{\rm Trace}({\bf Q}) = \frac{2}{\overline{K(\mathbf U)}}(U_1+U_6)\,.
\end{equation}
\end{prop}

To prove this proposition,
use (\ref{xy-c}) to relate the columns of ${\bf Z}$ to the ${\bf X}-{\bf Y}$
decomposition
\[
{\bf X}+\ri{\bf Y} = \left[\begin{matrix}
\langle\overline{{\bf c}_1},{\bf z}_1\rangle & 
\langle\overline{{\bf c}_1},{\bf z}_2\rangle \\
\langle\overline{{\bf c}_2},{\bf z}_1\rangle & 
\langle\overline{{\bf c}_2},{\bf z}_2\rangle \end{matrix}\right]\,.
\]
Hence
\[
{\bf Q} = \frac{1}{\overline{K}}
\begin{pmatrix}
\langle{\bf c}_1,{\bf z}_1\rangle & 
\langle{\bf c}_1,{\bf z}_2\rangle \\
\langle{\bf c}_2,{\bf z}_1\rangle & 
\langle{\bf c}_2,{\bf z}_2\rangle \end{pmatrix}
\left[\begin{matrix}
\langle\overline{{\bf c}_2},{\bf z}_2\rangle & 
-\langle\overline{{\bf c}_1},{\bf z}_2\rangle \\
-\langle\overline{{\bf c}_2},{\bf z}_1\rangle & 
\langle\overline{{\bf c}_1},{\bf z}_1\rangle \end{matrix}\right]\,,
\]
and so
\[
\begin{array}{rcl}
{\rm Trace}({\bf Q}) &=& \frac{1}{\overline{K}}\left(
\langle{\bf c}_1,{\bf z}_1\rangle\,\langle\overline{{\bf c}_2},{\bf z}_2\rangle
-\langle{\bf c}_1,{\bf z}_2\rangle\,\langle\overline{{\bf c}_2},{\bf z}_1\rangle
-\langle{\bf c}_2,{\bf z}_1\rangle\,\langle\overline{{\bf c}_1},{\bf z}_2\rangle
+\langle{\bf c}_2,{\bf z}_2\rangle\,\langle\overline{{\bf c}_1},{\bf z}_1\rangle
\right)\,,\\[3mm]
&=& \frac{1}{\overline{K}}\left(
{\rm det}\left[\begin{matrix}
\langle{\bf c}_1,{\bf z}_1\rangle & \langle{\bf c}_1,{\bf z}_2\rangle\\
\langle\overline{{\bf c}_2},{\bf z}_1\rangle &
\langle\overline{{\bf c}_2},{\bf z}_2\rangle \end{matrix}\right] +
{\rm det}\left[\begin{matrix}
\langle\overline{{\bf c}_1},{\bf z}_1\rangle & \langle\overline{{\bf c}_1},{\bf z}_2\rangle\\
\langle\overline{{\bf c}_2},{\bf z}_1\rangle & \langle\overline{{\bf c}_2},{\bf z}_2\rangle \end{matrix}\right] \right) \,,\\[7mm]
&=& \frac{1}{\overline{K}}\left(
\lbk {\bf c}_1\wedge\overline{{\bf c}_2},{\bf z}_1\wedge{\bf z}_2\rbk_2 +
\lbk \overline{{\bf c}_1}\wedge{\bf c}_2,{\bf z}_1\wedge{\bf z}_2\rbk_2 
\right)\\[3mm]
&=& \frac{1}{\overline{K}}\left(
\lbk {\bf c}_1\wedge\overline{{\bf c}_2}+ \overline{{\bf c}_1}\wedge{\bf c}_2,{\bf z}_1\wedge{\bf z}_2\rbk_2 \right)\\.
&=& \frac{2}{\overline{K}}\left(
\lbk {\bf e}_1\wedge{\bf e}_2 +{\bf e}_3\wedge{\bf e}_4,{\bf z}_1\wedge{\bf z}_2\rbk_2 \right)\\
&=& \frac{2}{\overline{K}}(U_1+U_6)\,,
\end{array}
\]
using
\[
{\rm Re}\left({\bf c}_1\wedge\overline{{\bf c}_2}\right) =
{\rm Re}\left( ({\bf e}_1 - \ri{\bf e}_3)\wedge({\bf e}_2+\ri{\bf e}_4)\right) =
{\bf e}_1\wedge{\bf e}_2 +{\bf e}_3\wedge{\bf e}_4\,,
\]
proving (\ref{trace-Q}).

Given a path ${\bf U}\in\bwedge^2(\R^4)$ the eigenvalues of ${\bf Q}$ can
be computed by substituting (\ref{det-Q}) and (\ref{trace-Q}) in (\ref{Q-eigs})
leading to
\[
\mu_{1,2} := {\rm e}^{\ri\kappa_{1,2}} = \frac{ U_1+U_6 \pm \sqrt{
4U_1U_6+2U_3U_4 - U_3^2-U_4^2}}{U_1+\ri U_3-\ri U_4 - U_6}\,.
\]
Using the properties of a Lagrangian plane, $U_2+U_5=0$
and $U_1U_6-U_2U_5+U_3U_4=0$, this formula reduces to
\[
\mu_{1,2} := {\rm e}^{\ri\kappa_{1,2}} = \frac{ U_1+U_6 \pm \i\sqrt{
4U_5^2 + (U_3+U_4)^2}}{U_1+\ri U_3-\ri U_4 - U_6}\,.
\]
The formula\footnote{In the most general case,
the eigenvalues $e^{\i\kappa_j}$ are the roots
of the following polynomial:
\[
P(\lambda)=\det(({\bf X}-\i{\bf Y})- \mu({\bf X}+\i{\bf Y}))\,,\quad \mu\in S^1\,.
\]
The coefficients of $P$ are antisymmetric multi-linear functions of $\mathbf Z$. 
As a consequence, they can always be expressed as a linear combination
of the minors of $\mathbf Z$ and hence, of the components of $\mathbf U$.
Therefore, it is possible to compute the eigenvalues   $e^{\i\kappa_j}$ 
from $\mathbf U$.} for the decomposition of the 
Maslov angles in  dimension 3 is given in Part 2 \cite{cdb-part2}.

\section{$\lambda$ dependence of the Maslov index}
\label{lambda-properties}
\setcounter{equation}{0}

Let ${\bf B}(x,\lambda)$ be as defined in (\ref{h1}) with the asymptotic
property (\ref{A_infty}).  When the $\lambda-$dependence of ${\bf B}(x,\lambda)$
takes a simple form one can say more about the $\lambda-$dependence of the
Maslov index.  For example, in \S\ref{sec-example-r2} the matrix ${\bf B}(x,\lambda)$
in (\ref{J-B-2}) satisfies
\[
\frac{\partial\ }{\partial\lambda}{\bf B}(x,\lambda) = \left[\begin{matrix}
1 & 0 \\ 0 & 0 \end{matrix}\right]\,,
\]
that is, it is positive semi-definite.  In the example in \S\ref{sec-tutorial-II}
the matrix ${\bf B}(x,\lambda)={\bf J}{\bf A}(x,\lambda)$ with ${\bf A}(x,\lambda)$
defined in (\ref{ode4}) has the property
\[
\frac{\partial\ }{\partial\lambda}{\bf B}(x,\lambda) =
\left[\begin{matrix} -1 & 0 & 0 & 0 \\
 0 & -1 & 0 & 0 \\ 0 & 0 & 0 & 0\\ 0 & 0 & 0 & 0 \end{matrix}\right]\,,
\]
which is negative semi-definite.

When $\partial_\lambda{\bf B}(x,\lambda)$ is semi-definite, the Maslov index is
a monotone function of $\lambda$.  For the case of gradient systems (as in
\S\ref{sec-example-r2} and \S\ref{sec-tutorial-II}) this property is
proved in Lemmas 3.3 and 3.7 of \textsc{Bose \& Jones}~\cite{bose}.
Related results are proved by \textsc{Arnold}~\cite{arnold2}
and generalizations of these results are proved by
\textsc{Chardard}~\cite{chardard-thesis}.
These results are summarized in

\begin{lem}
\label{lem:eigenvalues}
Assume that:
  \begin{itemize}
  \item  $\mathbf B(x,\lambda)$ is a smooth function with respect to $x$ and analytic with respect
   to $\lambda$.
  \item  There exists  $\mathbf B_{\infty}(\lambda)$, $\gamma>0$ 
    and $F>0$ such that 
    $\forall x,\lambda \quad \|\mathbf B(x,\lambda)-
    \mathbf B_{\infty}(\lambda)\|\leq Fe^{-\gamma|x|} $.
  \item  
The open set $\mathbb{X}=\R-\sigma_{ess}$ of real numbers is not empty.
\item $\partial_{\lambda} \mathbf B(x,\lambda)$  is semi-definite symmetric matrix.
\end{itemize}

If $[\lambda_1,\lambda_2]\cap \sigma_{ess}=\emptyset$ and $\lambda_1,\lambda_2\notin \sigma$, 
then $\mathsf{Maslov}(\lambda_2)-\mathsf{Maslov}(\lambda_1)$ is equal to the number of eigenvalues with multiplicity
in $[\lambda_1,\lambda_2]$.
\end{lem}

The first three assumptions are the usual hypotheses made to
prove the analyticity of the Evans function and the theorems 
linking eigenvalues and the zeros of the Evans function.
Using this lemma, it is possible to define the following invariant 
for the homoclinic orbit:

\begin{defn}
$\textsf{Maslov}^{\rm homoclinic}$ is defined as $\lim_{\lambda\to 0^+} \textsf{Maslov} (\lambda)$.
\end{defn}

If we make the hypotheses of Lemma \ref{lem:eigenvalues},
 $\textsf{Maslov}^{\rm homoclinic}$ is only dependent on $\mathbf A(x,0)$, and hence on the linearization
of the ODE satisfied by the homoclinic orbit.

In fact, it is possible to define $\mathsf{Maslov}^{\rm homoclinic}$ without any reference\footnote{For example, in the case
of a transverse orbit (i.e. $\dim(\mathrm E^s(x,0)\cap \mathrm E^u(x,0) =1$), we have:
$$\begin{array}{ll}
  \textsf{Maslov}^{\rm homoclinic}= 
\lim_{x\to+\infty} \lim_{\varepsilon\to 0^+,\varepsilon>0} \max( & \textsf{Maslov}(
\mathbf Z(]-\infty,x-\varepsilon],0), \mathrm E^s(x,0)),\\
& \textsf{Maslov}(
\mathbf Z(]-\infty,x+\varepsilon],0),\mathrm E^s(x,0))).
\end{array}$$}
to a parameter $\lambda$ (see \cite{bose,chen_hu,chardard-thesis}). 
In these references, the quantity $\textsf{Maslov}^{\rm homoclinic}$ is defined
like a Maslov index for the path $\mathbf Z(.,0)$ with respect 
to $\mathrm E^s(0)$. Unfortunately,  $\mathbf Z(.,0)$
does not necessarily admit a right end point, and when it
has one, the intersection with  $\mathrm E^s(0)$ is not trivial and
some extra work is needed.

Anyway, these geometrical constructions are not convenient for
numerical computations and we will not use them in the sequel.

\subsection{The Maslov index for large values of $|\lambda|$}
\label{subsec-largelambda}

When $\lambda\to+\infty$ (or $\lambda\to-\infty$ if the essential spectrum extends
to minus infinity) we expect the Maslov index to converge to some finite value.
This property is similar
to the property of the Evans function for large $\lambda$.  The hypotheses
are based on the  analogous result of \cite{pw92}, adapted to the setting
of the Maslov index.   

\begin{hypo}
\label{hyp_asymptotic}
Suppose
\begin{itemize}
  \item that there exists $\lambda_0\in\R$ such that $\sigma_{ess}$ is empty for all
$\lambda>\lambda_0$;
  \item 
  For large enough $\lambda$, ${\bf A}_{\infty}(\lambda)$ has no purely imaginary eigenvalues;
  \item
Let ${\bf K}(\lambda)$ be a symplectic $4\times 4$ matrix depending 
analytically on $\lambda$
whose first $2$ columns are a basis for ${\rm E}^u_\infty(\lambda)$
and  whose last $2$ columns are a basis for ${\rm E}^s_\infty(\lambda)$
Define
\[
{\bf F}(x,\lambda)= {\bf K}^{-1}(\lambda)({\bf A}(x,\lambda)-{\bf A}_{\infty}(\lambda)){\bf K}(\lambda)\,.
\]
\item Suppose that, for large enough $\lambda$:
\[
\begin{array}{rcl}
\quad && \int_{\R} |{\bf F}(x,\lambda)|\d x\quad\mbox{is bounded, uniformly in $\lambda$}\\
\quad &&\int_{|x|>x_0} |{\bf F}(x,\lambda)|\d x\quad\mbox{tends to $0$ when $x_0\to\infty$,
uniformly in $\lambda$}\\
\quad && \int_{\R} |{\bf F}^{(2)}(x,\lambda){\bf e}_1|\d x\quad\mbox{tends to $0$}\,.
\end{array}
\]
\end{itemize}
\end{hypo}

{\bf Remark.} If there exists $\lambda_0$ such that $\sigma_{ess}$ is empty for all
$\lambda<\lambda_0$ the above hypotheses can be modified accordingly.
\vspace{.15cm}

\begin{prop}
\label{asymptotic}
Assume that hypothesis \ref{hyp_asymptotic} is met by $\mathbf A(x,\lambda)$, then
\[
\lim_{\lambda\to +\infty} D(\lambda)=1\,,\quad \mbox{and}\quad \lim_{\lambda\to +\infty} 
\textsf{Maslov}(\cdot;\lambda)=0\,.
\]
\end{prop}

The proof is obtained by following the argument in Proposition 1.17
in \textsc{Pego \& Weinstein}~\cite{pw92}
and the Appendix of \textsc{Bridges \& Derks}~\cite{bd-analyticity}.
Let $\mathcal K(\lambda)$ be the matrix
whose entries are the $2\times 2$ minors of ${\bf K}(\lambda)$.
In terms of bialternate 
product, $\mathcal K(\lambda)={\bf K}(\lambda)\odot {\bf K}(\lambda)$.
 One proves that $\mathcal K(\lambda)^{-1}Y(.,\lambda)$ converges,
uniformly in $x$, to the 
constant vector ${\bf e}_1$ when $\lambda\to-\infty$.
Then, for large enough $\lambda$,
$\mathcal K(\lambda)^{-1}Y(.,\lambda)$ has a null Maslov index 
and so does $Y(.,\lambda)$.  See \cite{chardard-thesis} for a detailed proof.

\subsection{Defining a Maslov index  at $\lambda=0$ when the basic state
is approximated by hyperbolic periodic solutions}

In this paragraph, we assume that $\mathbf B(x,0)=D^2H(\widehat{\phi(x)})$ 
where $H:\R^{4}\to \R$ is the Hamiltonian function and $\widehat\phi$ is the
basic homoclinic solution of the autonomous system   
\begin{equation}
  \label{eq:aut_sys}
\mathbf J{\bf u}_x= \nabla H ({\bf u})\,.
\end{equation}
This hypothesis is natural in the solitary wave context and is indeed
satisfied by all examples presented in this article.

Suppose that  $\partial_{\lambda} \mathbf B$ is semi-definite near  $\lambda=0$
and that the dimension of the space of square integrable solutions 
of $\mathbf J {\bf u}_x = \mathbf B(x,0){\bf u}$ is one. 
Furthermore, suppose $\widehat\phi$ is approximated by 
hyperbolic $\frac{2\pi} k$-periodic solutions $\widehat\phi_k$, 
for which the Maslov index is well-defined \cite{F.Chardard2006,F.Chardard2007}.

If the Maslov index at  $\lambda=0$  is defined as the 
limit of the Maslov index of the periodic solutions $\widehat\phi_k$, when 
$k\to 0$, then, as shown in \cite{F.Chardard2007} under natural hypotheses,
the Maslov index at $0$ is the value of limit of the Maslov indices of the
periodic orbits when $\lambda$ 
is close to $0$ and it has the sign of  $f'(k)$ near $k=0$, where
$f(k):=H(\widehat\phi_k)$.

\section{A coupled reaction-diffusion equation with explicit Maslov index}
\label{sec-tutorial-II}
\setcounter{equation}{0}

Consider the system of reaction-diffusion equations
\begin{equation}\label{rd1}
\begin{array}{rcl}
\displaystyle\frac{\partial u}{\partial t} &=& \displaystyle
\frac{\partial^2u}{\partial x^2} 
-4u + 6u^2 - c(u-v) \\[4mm]
\displaystyle\frac{\partial v}{\partial t} &=& \displaystyle
\frac{\partial^2v}{\partial x^2} 
-4v + 6v^2 + c(u-v)\,,
\end{array}
\end{equation}
where $c$, the coupling
constant, is a non-zero real parameter, restricted to the values $c>-2$.
When $c>-2$, it is straightforward to show
that the trivial solution $u=v=0$ is stable in the time dependent problem,
and the trivial solution of the steady equation is hyperbolic.

This system has the exact steady solitary-wave solution
\[
u=v := \widehat u(x) = {\rm sech}^2(x)\,.
\]
Linearizing (\ref{rd1}) about the basic state $\widehat u$
and taking perturbations of the form 
\[
{\rm e}^{\lambda t}(u(x,\lambda),v(x,\lambda))\,,
\]
leads to the coupled ODE eigenvalue problem
\begin{equation}
  \label{eq:sec-tutorial-II-spec}
\begin{array}{rcl}
u_{xx} &=& (\lambda + 4 + c  - 12\widehat u(x) )\,u - c\, v \\[2mm]
v_{xx} &=& -c\, u + (\lambda + 4 + c - 12\widehat u(x) )\,v\,.
\end{array}  
\end{equation}

This eigenvalue problem can be written in the standard form
\begin{equation}\label{ode4}
{\bf u}_x = {\bf A}(x,\lambda){\bf u}\,,\quad {\bf u}\in\R^4\,,
\end{equation}
with ${\bf u}=(u,v,u_x,v_x)$ and 
\[
{\bf A}(x,\lambda) = \begin{pmatrix} 0 & 0 & 1 & 0 \\
0 & 0 & 0 & 1 \\ f(x,\lambda) & -c & 0 & 0 \\ -c & f(x,\lambda) & 0 & 0
\end{pmatrix}\,,\quad\mbox{with}\quad
f(x,\lambda) = \lambda + 4 + c - 12\,{\rm sech}^2(x)\,.
\]
The system (\ref{ode4}) is Hamiltonian: ${\bf JA}$ is symmetric.

The spectral problem (\ref{ode4}) can be solved explicitly.  Write
the second-order problem in the form,
\begin{equation}\label{eqn-2}
\begin{pmatrix}u\\ v\end{pmatrix}_{xx} = f(x,\lambda)\begin{pmatrix}u\\ v\end{pmatrix} -c \left[\begin{matrix} 0 & 1 \\ 1 & 0 \end{matrix}\right]
\begin{pmatrix}u\\ v\end{pmatrix}\,.
\end{equation}
Let
\[
\begin{array}{rcl}
u &=& \widetilde u - \widetilde v \\[2mm]
v &=& \widetilde u + \widetilde v \,.
\end{array}
\]
Then substitution into (\ref{eqn-2}) leads to the decoupled system
\begin{equation}\label{decoupl}
\begin{array}{rcl}
\widetilde u_{xx} + 12\,{\rm sech}^2(x)\,\widetilde u &=& (\lambda + 4)\widetilde u\\[2mm]
\widetilde v_{xx} + 12\,{\rm sech}^2(x)\,\widetilde v &=& (\lambda + 4+2c)\widetilde v \,.
\end{array}
\end{equation}
These two systems have explicit solutions (cf. Appendix \ref{app_spectral}),
and using these results one finds that there are exactly six eigenvalues
for the spectral problem (\ref{eq:sec-tutorial-II-spec}):
\[
\begin{array}{rcl}
\lambda_1 &=& -3 - 2c \,,\quad \lambda_2 = -3  \,,\quad
\lambda_3 = -2c  \,,\\[2mm]
\lambda_4 &=& 0 \,,\quad \lambda_5 = 5 - 2c  \,,\quad
\lambda_6 = 5\,. 
\end{array}
\]
The essential spectrum is
\[
\sigma_{\rm ess} = \{\lambda\in\R\ :\ \lambda\leq -4\}\cup
\{\lambda\in\R\ :\ \lambda\leq -4-2c\}\,.
\]

When $c=0$ then there are three
double eigenvalues: 
\[
\lambda_1=\lambda_2 = -3\,,\quad \lambda_3=\lambda_4 = 0\,,\quad
\lambda_5=\lambda_6 = 5\,.
\]
For $c$ nonzero and small the eigenvalues $\lambda_1$, $\lambda_3$
and $\lambda_5$ are perturbed to the left (when $c>0$) and
to the right (when $c<0$).  Hence positive coupling is stabilizing and
negative coupling is destabilizing.
\begin{figure}
\begin{center}
\includegraphics[width=6cm]{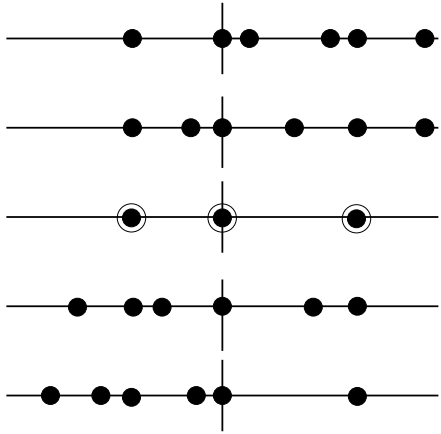}
\caption{Plots of the spectrum in the $\lambda$-plane for $c=-2$ (starting
at the top), $c=-1$, $c=0$,
$c=1$ and $c=3$. For each plot, the intersection of the axes
is the origin $\lambda=0$.}
\label{fig-cfig}
\end{center}
\end{figure}

When $c=-2$ (the lower bound on $c$ for stability of the zero state
and hyperbolicity) there are four positive eigenvalues, which is
the maximum number of positive eigenvalues. 

At $c=-\frac{3}{2}$ one
of the positive eigenvalues passes through zero, leaving $3$ positive
eigenvalues.  Then when $c=0$ another eigenvalue arrives at zero
leaving two positive eigenvalues.  Then when $c=\frac{5}{2}$ a
third eigenvalue passes through zero.  For all $c>\frac{5}{2}$
there is one positive eigenvalue, one zero eigenvalue and four
negative eigenvalues.  The configuration of the eigenvalues as a function
of $c$ is shown in Figure \ref{fig-cfig}.

When the system decouples into two subsystems, the Maslov index is the sum
of the Maslov indices of the subsystems
\begin{equation}\label{sum_formula}
\textsf{Maslov}^{2D}\oplus\textsf{Maslov}^{2D} = \textsf{Maslov}^{4D}\,.
\end{equation}
If we choose a Lagrangian frame
 of the form  $\begin{pmatrix}
  X\\
  Y
\end{pmatrix}=
\begin{pmatrix}
  X_1  & 0\\
  0   & X_2\\
  Y_1  & 0\\
  0   & Y_2\\
\end{pmatrix}
$ to represent the unstable space,
we have indeed $\det(X+\i Y)=\det(X_1+\i Y_1)\det(X_2+\i Y_2)$. From this property,
\eqref{sum_formula} follows easily.

This property is obvious in the present example, since the system decouples.
The Maslov index for the first 2D system in (\ref{decoupl}) is tabulated below.
\vspace{.25cm}

\begin{center}
\begin{tabular}{|c|c|c|c|c|}
\hline
\phantom{$\bigg|$}
$\lambda$ & $-4<\lambda<-3$  & $-3<\lambda<0$ & $0<\lambda<5$ & $\lambda>5$ \\
\hline
\phantom{$\bigg|$}$\textsf{Maslov}_1^{2D}$ & $3$ & $2$ & $1$ & $0$ \\
\hline
\end{tabular}
\end{center}

\vspace{.25cm}

The Maslov index for the second 2D system in (\ref{decoupl}) is tabulated here.

\begin{center}
\begin{tabular}{|c|c|c|c|c|}
\hline
\phantom{$\bigg|$}
$\lambda$ & $-4-2c<\lambda<-3-2c$  & $-3-2c<\lambda<-2c$ & $-2c<\lambda<5-2c$ & $\lambda>5-2c$ \\
\hline
\phantom{$\bigg|$}$\textsf{Maslov}_2^{2D}$ & $3$ & $2$ & $1$ & $0$ \\
\hline
\end{tabular}
\end{center}
\vspace{.5cm}

The Maslov index for the full 4D system for any $\lambda$
is then obtained by fixing $c$ and then applying the sum formula (\ref{sum_formula}).
For example, fix $c=-1$ and compute $\textsf{Maslov}^{\rm homoclinic}$,
\[
\textsf{Maslov}^{\rm homoclinic}\bigg|_{c=-1} = \lim_{\lambda\to0^+} \left[
 \textsf{Maslov}_1^{2D}\oplus \textsf{Maslov}_2^{2D}\right]\bigg|_{c=-1} = 3\,.
\]

A summary of the Maslov index of the homoclinic orbit as a function of
$c$ is given in the following table.
\vspace{.15cm}

\begin{center}
\begin{tabular}{|c|c|c|c|c|}
\hline
\phantom{$\bigg|$}
$\lambda$ & $-2<c<-\frac{3}{2}$  & $-\frac{3}{2}<c<0$ & $0<c<\frac{5}{2}$ & $c>\frac{5}{2}$ \\
\hline
\phantom{$\bigg|$}$\textsf{Maslov}^{\rm homoclinic}$ & $4$ & $3$ & $2$ & $1$ \\
\hline
\end{tabular}
\end{center}

\section{Numerical implementation -- approximation by periodic orbits}
\label{sec-periodic}
\setcounter{equation}{0}

When the solitary wave is approximated by a hyperbolic periodic orbit,
the Maslov index of the periodic orbit is computed using the Maslov angle
(\ref{angle-def-1})-(\ref{maslov-angle}) in the exterior algebra
representation (e.g. equation (\ref{kappa-wedge2})).
An algorithm for this case has been proposed in  \cite{F.Chardard2006},
and a proof of convergence of the Maslov index in the limit as the
periodic solution converges to the solitary wave is given in
\cite{F.Chardard2007}. A first possibility would be to compute the
Maslov index for a sufficiently good periodic approximation
of the solitary wave and to compute the its Maslov index.

However, it is possible to avoid the computation of 
periodic approximants and to adapt the algorithm presented 
in \cite{F.Chardard2006} to the homoclinic case.
Fix $\lambda$ and a solitary
wave.  The steps in this new algorithm are as follows. 
\begin{enumerate}
\item Choose a large enough interval $[-L,L]$.
\item Compute the eigenvalue with largest real part of $\mathbf A^{(2)}_\infty(\lambda)$,
denoted by $\sigma_+(\lambda)$,
and its associated eigenvector $\zeta^+(\lambda)$.  
\item Integrate equation 
\begin{equation}\label{key}
{\bf U}^+_x = [{\bf A}^{(2)}(x,\lambda) -\sigma_+(\lambda){\bf I}]{\bf U}^+\,,
\end{equation}
on $[-L,L]$, taking 
$\zeta^+(\lambda)$ as initial condition at $x=-L$, using any
standard numerical integration scheme.  The justification for
the arbitrariness in choice of numerical scheme is given in
Appendix \ref{sec_attractivity}. 
\item ${\bf U}^+(L,\lambda)$ and ${\bf U}^+(-L,\lambda)$ are nearly collinear, and an 
approximation to the
Evans function is determined from their constant of proportionality
\[
{\bf U}^+(L,\lambda) = D(\lambda){\bf U}^+(-L,\lambda) + \textsf{Error}\,,
\]
where the error is generally of the order of machine precision.
\item Compute ${\rm e}^{\ri\kappa(x)}$ using equation (\ref{kappa-wedge2})
or analogous representation.
\item Compute a lift of $\kappa(x)$ and
choose the stepsize $\Delta x$ so that $|\kappa(x+\Delta x)-\kappa(x)|<\pi$.
\item Compute the Maslov index using (\ref{maslov-angle}).
\end{enumerate}
\vspace{.15cm}

There are a number of sources of error in the algorithm.
There is an approximation error due to the fact that the solitary wave is
approximated by a periodic orbit.  
Two parameters have to be chosen: $L$ and the step size $\Delta x$.  The choice
of step size is a familiar source of error.
The consistency error of the numerical integration scheme will be of
the form $C\,\Delta x^p$, for some natural number $p$, at each step, where
$C$ is a constant
depending on the derivatives of $\mathbf A$.
The choice of numerical scheme will also impose some stability condition.

Since the Maslov index is an integer, the proposed scheme will give 
the Maslov index if the relative error on
${\bf U}^+(\cdot,\lambda)$ is small enough.
However, if $\sup_\R \|{\bf U}^+(\cdot,\lambda)\|$ is very small (for example when
the Evans function is small),
the relative error may be too big and lead to a miscomputed Maslov index.
This is the case when $\lambda$ is an eigenvalue or near an eigenvalue, since
there is an integer-valued jump in the Maslov index at eigenvalues.

\section{Numerical implementation -- intersection index algorithm}
\label{sec-numer-algor-1}
\setcounter{equation}{0}

The numerical algorithm based on the intersection index is similar to the 
algorithm in \S\ref{sec-periodic} except that the computation of the angle $\kappa(x)$
is replaced by the computation of the angles $\kappa_1$ and $\kappa_2$.
The previous algorithm can be modified as follows:

\begin{enumerate}
\item Choose a large interval $-L\leq x \leq L$. 
Initialize $\textsf{Maslov}$ to 0.
\item 
Construct a symplectic matrix ${\bf K}(\lambda)$ such that
\[
{\bf K}(\lambda)\begin{pmatrix} I\\
0
\end{pmatrix}\quad\mbox{and}\quad {\bf K}(\lambda)\begin{pmatrix} 
0\\
I
\end{pmatrix}\,,
\]
represent the stable and unstable spaces of $\mathbf A_{\infty}(\lambda)$.
${\bf K}(\lambda)$ defines a symplectic change of coordinates in which 
the coordinates of stable and unstable spaces in the exterior algebra are 
respectively
\[
{\bf U}_0=\begin{pmatrix} 
1\\
0\\
\vdots\\
0
\end{pmatrix}\quad\mbox{and}\quad {\bf V}_0=\begin{pmatrix} 
0\\
\vdots\\
0\\
1
\end{pmatrix}\,.
\]
Let $\mathcal K(\lambda)$ be the matrix
whose entries are the $2\times 2$ minors of ${\bf K}(\lambda)$.
In terms of bialternate 
product, $\mathcal K(\lambda)={\bf K}(\lambda)\odot {\bf K}(\lambda)$.

\item Compute the eigenvalue with largest real part, $\sigma_+(\lambda)$,
of $\mathbf A^{(2)}_\infty(\lambda)$.
\item Define $\widetilde{{\bf U}^+}(x,\lambda) = {\mathcal K}(\lambda)^{-1}
{\bf U}^+(x,\lambda)$ and integrate the equation for $\widetilde{{\bf U}^+}$,
\begin{equation}\label{key-1}
\frac{d\ }{dx}\widetilde{{\bf U}^+} = 
[\mathcal K(\lambda)^{-1}{\bf A}^{(2)}(x,\lambda)\mathcal K(\lambda) -\sigma_+(\lambda){\bf I}])\widetilde{{\bf U}^+}\,,
\end{equation}
on $[-L,L]$, taking ${\bf U}_0$ as initial condition for 
$\widetilde{{\bf U}^+}$ at $x=-L$, using any standard numerical integration scheme.

\item Compute the angles $(\kappa_1$ and $\kappa_2$
corresponding to $\widetilde{{\bf U}^+}$ over $[-L,L]$. 
If an angle $\kappa_i$ crosses $2\pi\Z$ between $x$ and $x+\Delta x$,
update the value of the Maslov index to:

\[
\textsf{Maslov}\mapsto\textsf{Maslov} + {\rm sign}\,
\left(\kappa_i(x+\Delta x)-\kappa_i(x)\right)\,.
\]

\item Return $\widetilde{{\bf U}^+}(L,\lambda)\wedge {\bf U}_0$
as an approximation to the Evans function.
\item At $x=+L$, return the value of the Maslov index. 
\end{enumerate}

\section{The Maslov index of solitary wave solutions of KdV5}
\setcounter{equation}{0}
\label{sec-multi-pulse}

In this section we study the Maslov index as a function of
$\lambda$ for the ODE eigenvalue problem
\begin{equation}\label{4d-example}
\phi_{xxxx} - P\phi_{xx} + a(x) \phi = \lambda\phi\,,
\end{equation}
where $\phi(x,\lambda)$ is scalar valued,
$P$ is a real parameter and $a(x)$ is a localized function which
satisfies $a(x)\to a_{\infty}$ as $x\to\pm\infty$, with exponential decay 
of $a(x)$ at infinity.
For definiteness it is assumed that $a_\infty>0$.
The ODE (\ref{4d-example}) can be put in the form (\ref{h1}) with
\begin{equation}\label{u-B-def}
{\bf u} = \begin{pmatrix} \phi \\ \phi_{xx} \\ \phi_{xxx}-P \phi_x \\ \phi_x \end{pmatrix}
\quad\mbox{and}\quad 
{\bf B}(x,\lambda) = \left[\begin{matrix} 
a(x)-\lambda & \hfill 0 & \ 0 & \ 0 \\ 0 & -1 & \ 0 & \ 0 \\
0 & \hfill0 & \ 0 & \ 1 \\ 0 & \hfill0 & \ 1 & \ P \end{matrix} \right]\,.
\end{equation}
The spectrum of the system at infinity ${\bf A}_\infty(\lambda)=
{\bf J}^{-1}{\bf B}_\infty(\lambda)$ has the characteristic polynomial
\begin{equation}\label{kdv5-poly}
{\rm det}[ {\bf A}_\infty(\lambda) - \mu{\bf I}] =
\mu^4 - P\mu^2 + a_\infty-\lambda\,,
\end{equation}
With $a_\infty>0$ and $\lambda=0$ the four roots are hyperbolic for all
$P$ such that
\[
P+2\sqrt{a_\infty}>0\,,
\]
which is assumed to be satisfied henceforth.  When $\lambda\neq0$ the essential
spectrum will form the boundary of the hyperbolic region.
The essential spectrum is
\[
\sigma_\ess  =  \{  \lambda\in\R\ :\ \lambda=a_\infty + Ps^2 + s^4\,,\  s\in\mathbb{R} \}\,.
\]
When
\[
\lambda < \lambda^{\rm edge} = a_\infty - \frac{1}{8}P(P-|P|)\,,
\]
the spectrum of ${\bf A}_\infty(\lambda)$ is hyperbolic.  
Hence, all the hypotheses
for the existence of the Evans function and the Maslov index are satisfied.  We will
apply this theory to determine the Maslov index of a class of homoclinic orbits.

The eigenvalue problem (\ref{h1}) appears in the linearization about
a solitary wave solution of the fifth-order Korteweg-de Vries equation (KdV5). 
KdV5 appears for example as a model equation in plasma physics,
and in the study of capillary--gravity water waves
\cite{kawahara,kodama,diasiooss,champ,bd-sima,diaskuz,bdg02}. 

To see the role of (\ref{h1}) in the linearization of KdV5, consider
the following form of the fifth-order KdV equation
relative to a moving frame, moving at speed $c$, 
\begin{equation}\label{newKdV5}
\frac{\partial \phi}{\partial t} - c\frac{\partial\phi}{\partial x}
+ \frac{\partial\ }{\partial x}\left(
\phi^{q+1}\right) + P\frac{\partial^3\phi}{\partial x^3} -
\frac{\partial^5\phi}{\partial x^5}=0\,,\quad q\geq1\,.
\end{equation}
A further scaling can be introduced so that $c=1$, but including $c$ is useful for comparing with
results in the literature on KdV5. Effectively, $q$ is a third parameter, but its
value is restricted to natural numbers.

Steady solutions of (\ref{newKdV5}), that decay to zero as $x\to\pm\infty$
satisfy the $4^{th}-$order ordinary
differential equation
\begin{equation}\label{steady_kdv}
\phi_{xxxx}-P\phi_{xx}+c\phi-\phi^{q+1}=0. 
\end{equation}
The system (\ref{newKdV5}) linearized about a solitary wave $\widehat\phi(x)$ solution of (\ref{steady_kdv})
takes the form
\[
\frac{\partial\phi}{\partial t} = \frac{\partial\ }{\partial x}\left(
{\mathscr L}\phi \right)\,,
\]
with
\begin{equation}\label{sa_operator}
\mathscr{L}\phi := \phi_{xxxx} - P\,\phi_{xx} + c\phi
- (q+1)\widehat\phi(x)^q\,\phi \,.
\end{equation}
In this case, $a_\infty=c$.
There are two spectral problems:
\begin{equation}\label{two-spectral-def}
\mathscr{L}\phi = \lambda\phi \quad\mbox{and}\quad {\bf L}\phi = \widehat\lambda \phi\,,
\quad {\bf L}\phi := \frac{d\ }{dx}\mathscr{L}\,.
\end{equation}
The operator $\mathscr{L}$ is self-adjoint (in a suitably-chosen Hilbert space) and
so $\lambda\in\R$ whereas ${\bf L}$ is not self-adjoint and $\widehat\lambda$ --
which is the stability exponent -- can in general be complex.  The relationship between
these two eigenvalue problems is discussed in \S\ref{stability-kdv5}.  
First the Maslov index of the spectral problem $\mathscr{L}\phi = \lambda\phi$, 
which can be put in the form (\ref{h1}), is studied.

\begin{figure}
\centering
\includegraphics[width=17cm]{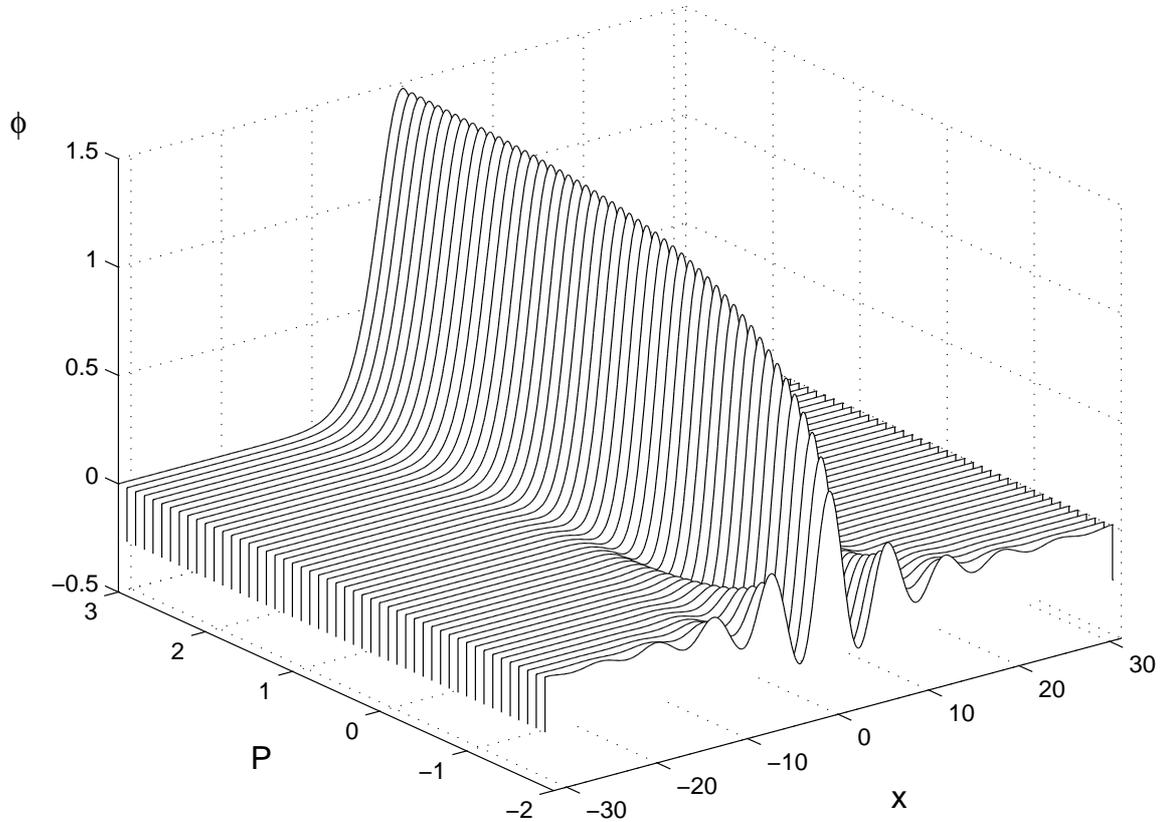}
\caption{ 
Numerically computed solitary waves for the KdV 5  equation for the case
$q=1$, $c=1$ ($a_\infty=1$) and $-2<P<3$.
For each value of $P$, the value of the corresponding 
unimodal homoclinic solution $\phi$ is plotted as a function of $x$.}
\label{solp1}
\end{figure}

The ODE (\ref{steady_kdv}) has been extensively studied and many solitary wave solutions
have been found; a classification is given in \cite{bct}.  
There are some special cases where
explicit solitary wave solutions can be constructed.  An example is the explicit solution
$\widehat\phi(x)=\frac{35}{24}\sech^4\left(\frac{x}{2\sqrt{6}}\right)$
which exists when $q=1$, $c=1$ and $P=\frac {13}{6}$.  However, the interesting solutions 
of (\ref{steady_kdv}) need
to be computed numerically.  They can be computed using a spectral method (approximate
the solitary wave by a periodic function of large wavelength and then use Fourier series
to represent it), or in the case of symmetric solitary waves a shooting algorithm can be
used.  We used both methods to compute solitary waves.  
An example of the family of
one-mode solitary waves as a function of $P$, computed using a spectral method,
is shown in Figure \ref{solp1}. 

Although these solitary waves are solutions of the model ODE, they are representative
of solutions of the full water-wave problem.  
\textsc{Dias, Menasce \& Vanden-Broeck}~\cite{dmv} have found large-amplitude branches
of these solutions in the full water-wave equations.

Symmetric solutions are computed numerically using a shooting method:
the starting point is an element of the tangent space of
the unstable manifold and the ending point is a symmetric point.
The time step typically used is $\frac 1 {1000}$ and the integrator chosen was
the fourth-order Runge-Kutta method. For this integrator and this time
step, the error on the Evans function at $0$ (which is then equal to $0$)
was of order $10^{-12}$. This is quite close to the minimum error reachable
by using double precision.

The ODE (\ref{steady_kdv}) can be characterized as
a Hamiltonian system on $\R^4$ and the Hamiltonian in the original coordinates is
\begin{equation}\label{kdv5-energy}
E(\phi)= \fr\phi_{xx}^2 +\fr P\phi_x^2 -\fr c\phi^2 + \frac{1}{q+2}\phi^{q+2}-\phi_x\phi_{xxx}\,.
\end{equation}
The function $E(\phi)$ is constant along solutions of (\ref{steady_kdv}) (i.e. $\frac{dE}{dx}=0$).
Physically, for equations like KdV5, this quantity is associated with momentum flux.
For simplicity, we
will just refer to it as energy.  The
energy of the periodic approximants gives some information about the nature of
the limiting homoclinic orbit. 
  
\begin{figure}
\centering
\includegraphics[width=10cm]{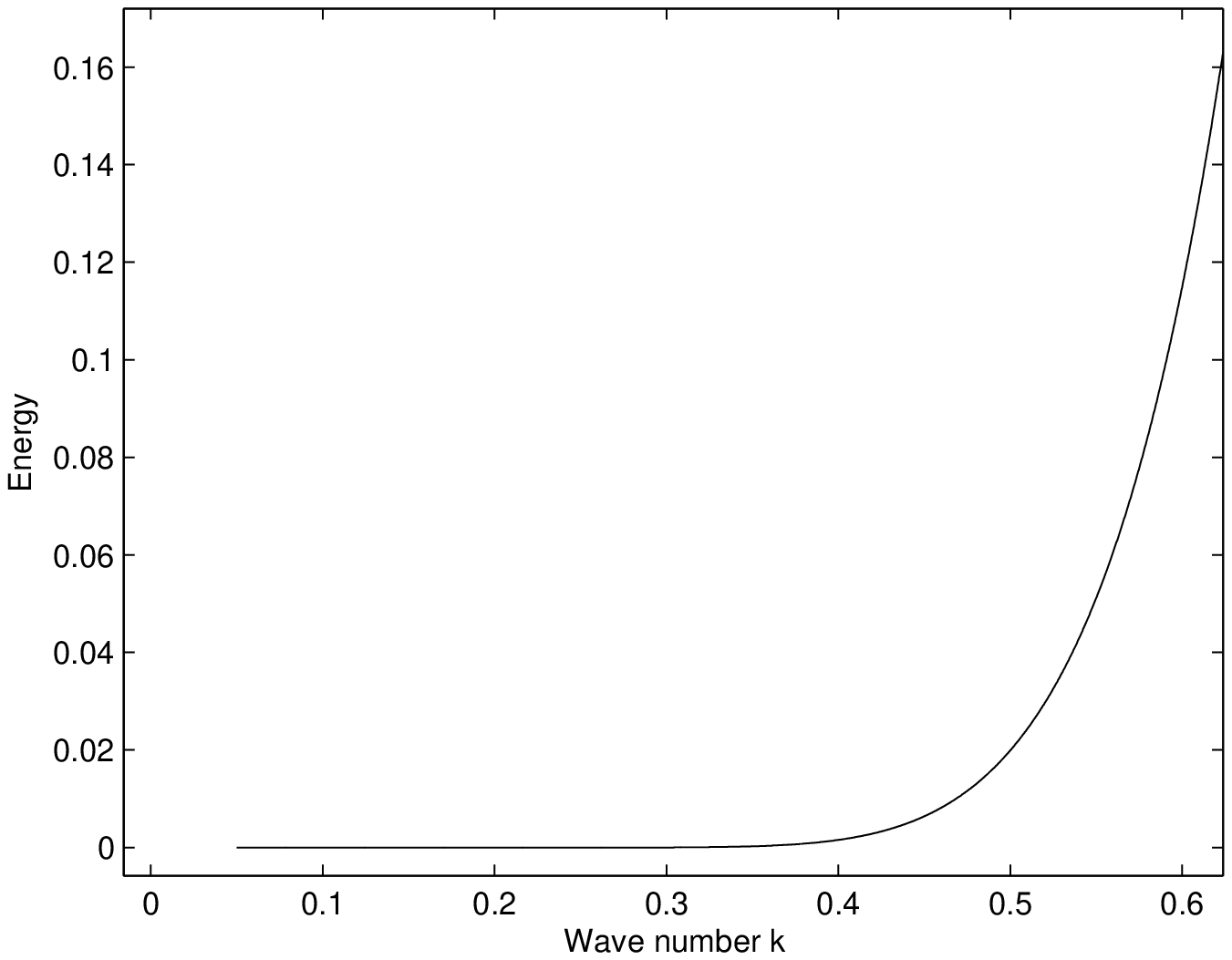}
\caption{Energy of the $\frac {2\pi} k$-periodic solutions as function of $k$
for $q=1$, $c=1$ and $P=\frac{13}6$.}
\label{fig_energy_1}
\end{figure}

The energy can be plotted as a function of wavenumber $k$ along a branch of periodic solutions
as the wavelength tends to infinity ($k\to0$) as a function of $q$ and $P$. 
First the case $P=\frac{13}{6}$, $c=1$
and $q=1$ is considered and it is shown in Figure \ref{fig_energy_1}.

\begin{figure}
\centering
\includegraphics[width=10cm]{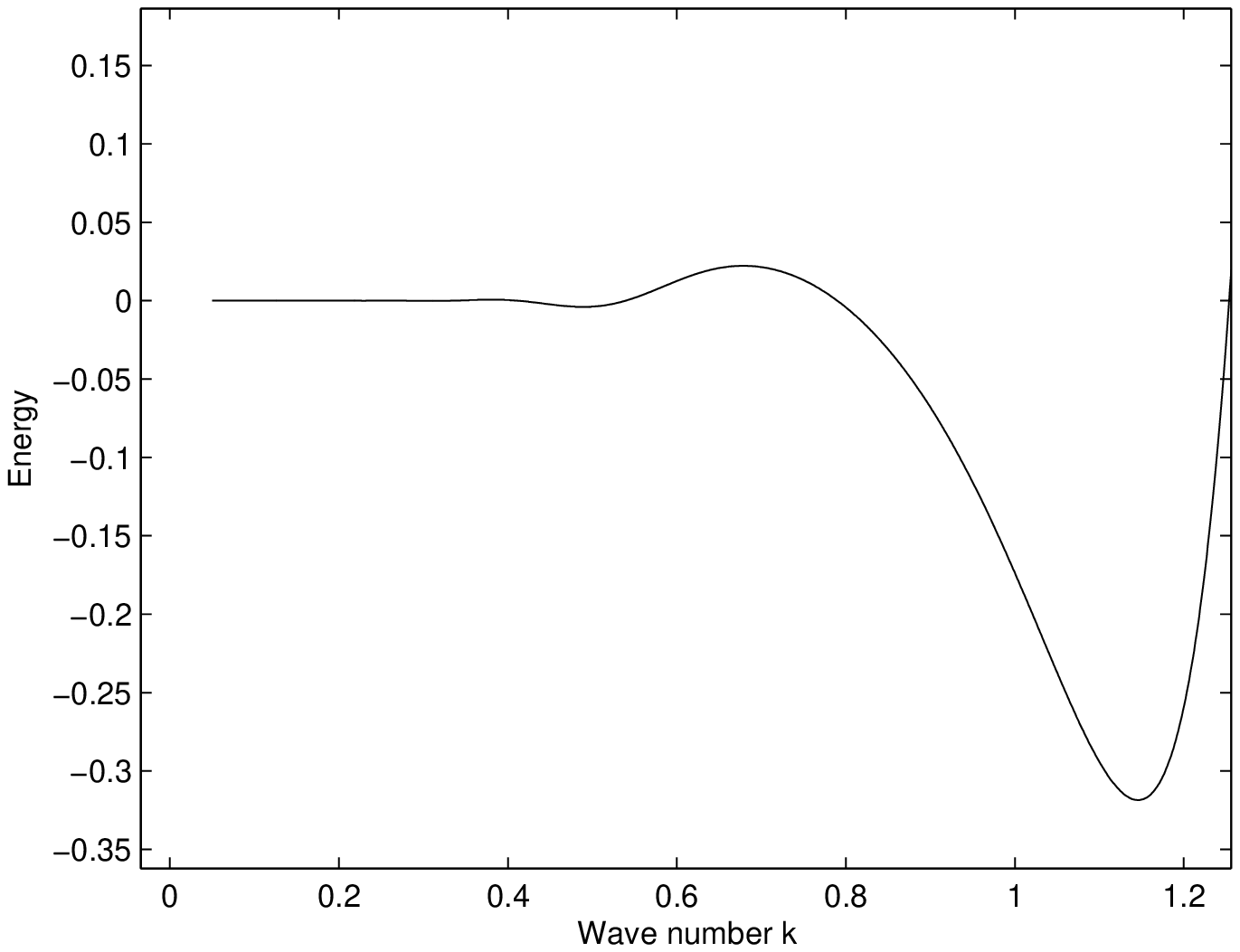}
\caption{Energy of the $\frac{2\pi}k$-periodic solutions as a function of $k$
for $q=1$, $c=1$ and $P=-1$.}
\label{fig_energy_2}
\end{figure}

\begin{figure}
\centering
\includegraphics[width=10cm]{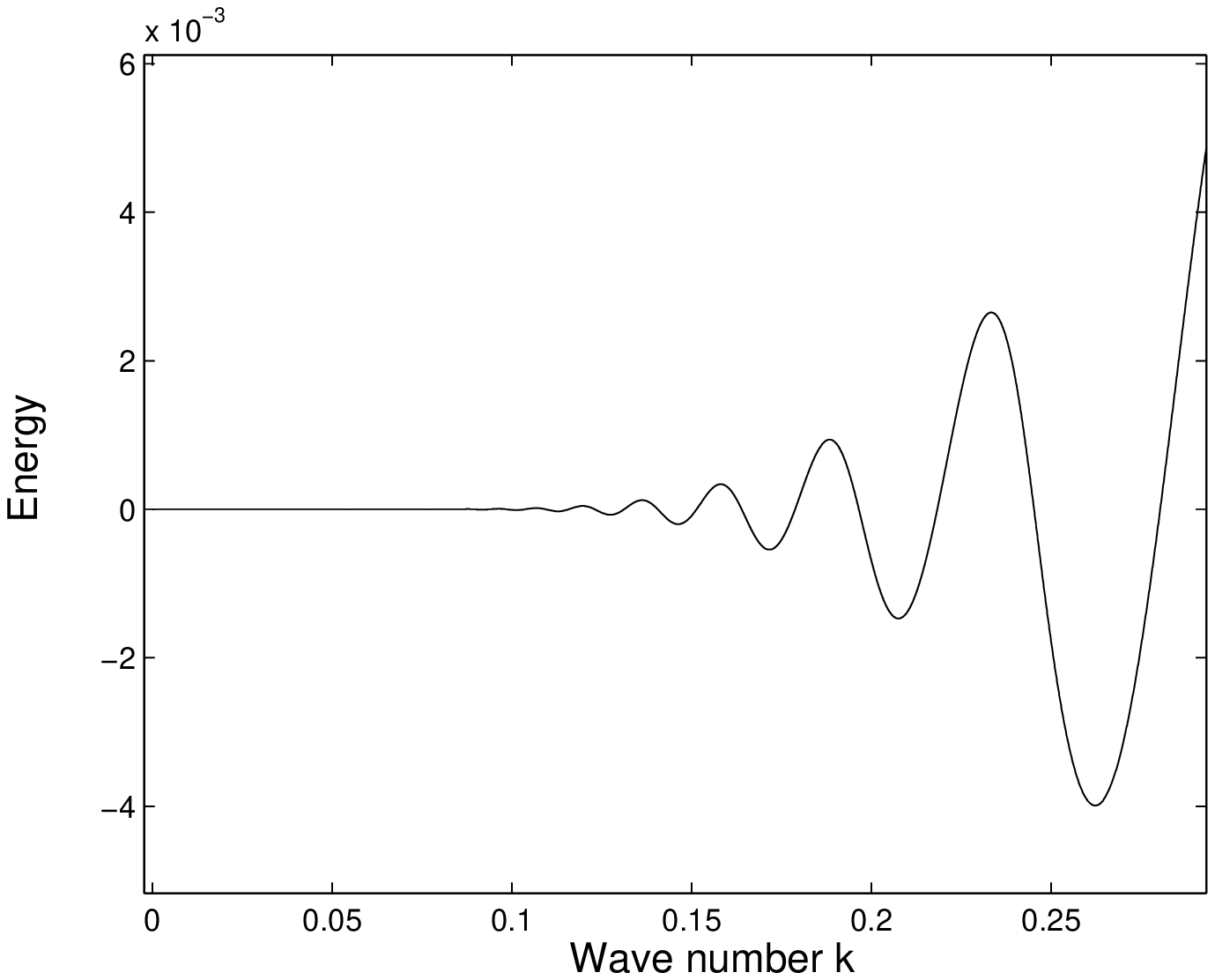}
\caption{Energy of the $\frac{2\pi}k$-periodic solutions as a function of $k$
for $q=1$, $c=1$ and $P=-1.9$, for low values of $k$.}
\label{fig_energy_3}
\end{figure}

In this case the energy is a monotone function of wavenumber and
the convergence $k\to0$ is rapid.
Keeping $q=1$ and $c=1$ but decreasing $P$ to $P=-1$ begins to show oscillations indicative of
a Shilnikov-type bifurcation as shown in Figure \ref{fig_energy_2}.
Decreasing $P$ further to $P=-1.9$ shows more dramatically the Shilnikov-type 
oscillations, as shown in Figure \ref{fig_energy_3}.
In Figure \ref{fig_energy_3} a sequence of bifurcations occurs along the branch.
Each point on the energy-wavenumber
diagram where $E'(k)=0$ corresponds to a saddle-centre
bifurcation of Floquet multipliers.  There are always two Floquet multipliers at
$+1$ due to the fact that (\ref{steady_kdv}) is autonomous.  When $E'(k)=0$
an additional pair of Floquet multipliers coalesces at $+1$.  
Each one of these saddle-centre bifurcations of the branch of periodic orbits
leads to a secondary homoclinic bifurcation \cite{bd05}.  So, in addition to the
limiting homoclinic orbit that we are principally interested in,
there is a countable number of other orbits generated along the branch, which are 
homoclinic to the branch of periodic orbits.  Although there is an infinite number
of bifurcations along the branch our numerical results show that the Maslov index of the limiting
homoclinic orbit is finite.

\subsection{Computing the Maslov index as a function of $\lambda$}

\begin{figure}
  \centering
\includegraphics[width=10cm]{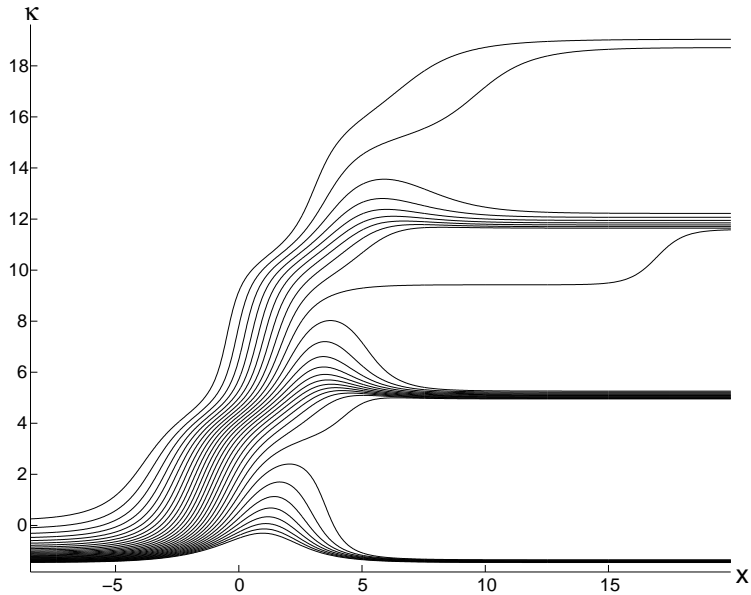}
\caption{$\kappa$ as a function of $x$ for the following values
of $\lambda$:$-2,-1.9,\ldots,0.9$.
$\kappa$ is $\lambda$-growing. The parameter values are $P=\frac{13}{6}$,
$c=1$ and $q=1$.}
\label{kappap1exact}
\end{figure}

First consider the case $P=\frac{13}{6}$, $q=1$ and $c=1$ where the unimodal solitary wave
is given explicitly.  The lifts $\kappa(x)$ of the Maslov angle for this system
are plotted as a function of $x$ in Figure \ref{kappap1exact} for various values of $\lambda$.
\begin{figure}
\centering
\includegraphics[width=10cm]{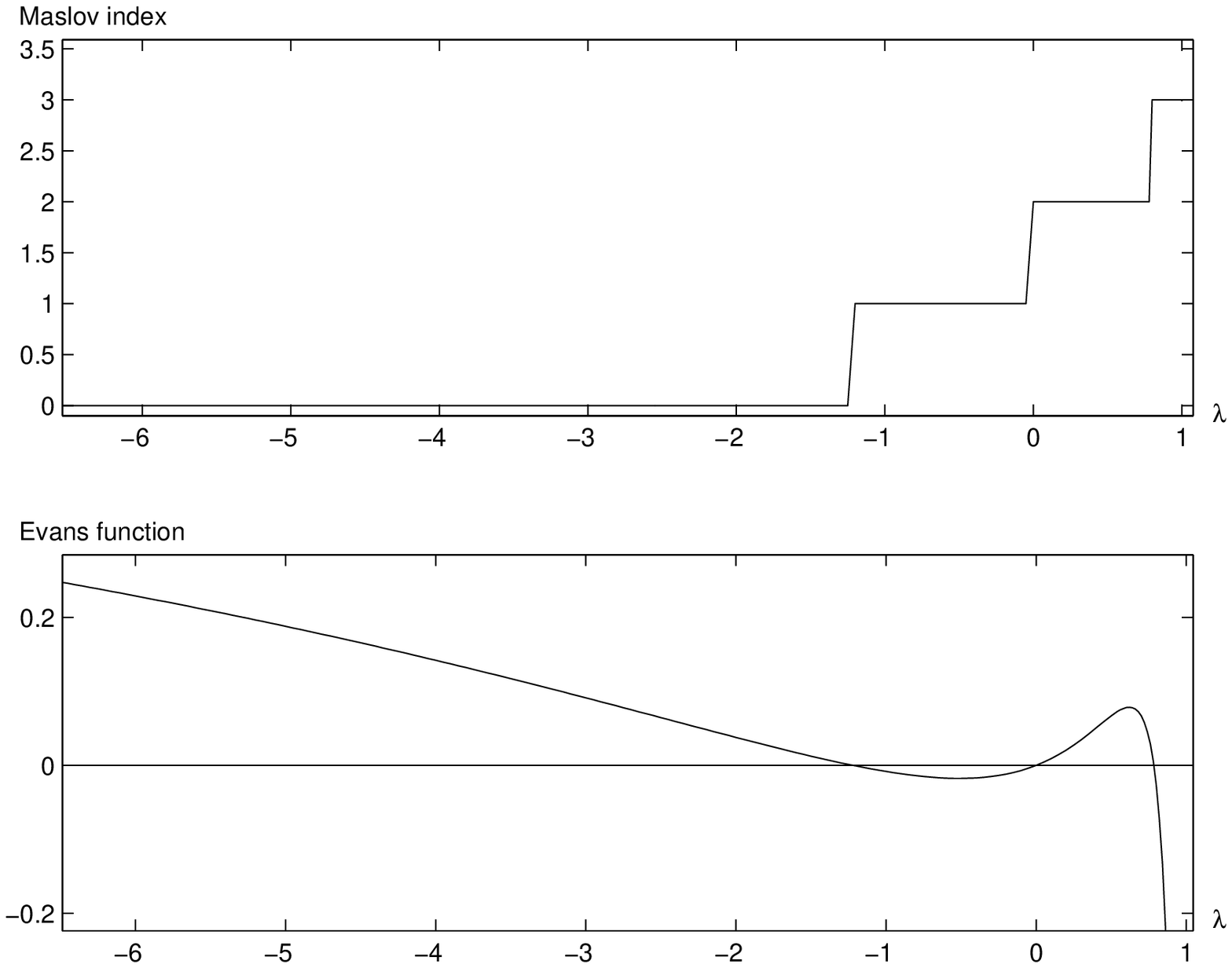}
\caption{Evans function and Maslov index as a function of $\lambda$
for the explicit unimodal solitary wave solution when $P=\frac{13}{6}$,
$c=1$ and $q=1$.  In this case, the Maslov index is $\lambda$-growing.}
\label{p1exact}
\end{figure}
In Figure \ref{p1exact}, the corresponding  Maslov indices 
have been plotted as a function of $\lambda$.  The Evans function shows that
$\mathscr{L}$ has exactly three eigenvalues in this case.  Denote these eigenvalues by
\[
\lambda_1 < \lambda_2=0 < \lambda_3 \,.
\]
The qualitative behaviour of the Maslov index in this case
is similar to the example on $\R^2$ in \S\ref{sec-example-r2}.  The values of the Maslov are
shown in the table below.  The Maslov index in this case is computed using
the Maslov angle, and this Maslov index is denoted by $\textsf{Maslov}(\kappa,\lambda)$.
\vspace{.15cm}

\begin{center}
\begin{tabular}{|c|c|c|c|c|}
\hline
\phantom{$\bigg|$}
$\lambda$ & $\lambda<\lambda_1$  & $\lambda_1<\lambda<\lambda_2$ & $\lambda_2<\lambda<\lambda_3$ & $\lambda>\lambda_3$ \\
\hline
\phantom{$\bigg|$}$\textsf{Maslov}(\kappa,\lambda)$ & $0$ & $1$ & $2$ & $3$ \\
\hline
\end{tabular}
\end{center} 
\vspace{.25cm}

Note that the Maslov index in each region predicts the number of eigenvalues of 
$\mathscr{L}$ in each $\lambda$ interval.
\vspace{.15cm}

\begin{center}
\begin{tabular}{|c|c|c|c|c|}
\hline
\phantom{$\bigg|$}
$\lambda$ region & $\lambda <\lambda_1$  & $\lambda<\lambda_2$ &  $\lambda<\lambda_3$ & $\lambda<\lambda^{\rm edge}$ \\
\hline
\phantom{$\bigg|$}$\#\ \textsf{Eigs}(\mathscr{L})$ & $0$ & $1$ & $2$ & $3$ \\
\hline
\end{tabular}
\end{center} 
\vspace{.25cm}

It is immediate from this table that 
\[
\textsf{Maslov}^{\rm homoclinic} = \lim_{\lambda\to 0^+}\textsf{Maslov}(\kappa,\lambda) = 2\,.
\]
The operator $\mathscr{L}$ has exactly one negative eigenvalue in this case. 
Our calculations indicate
that this is the case for all the unimodal homoclinic orbits.  It is easy to show analytically
that the Maslov index of a unimodal homoclinic orbit is greater than or equal to 2.
An elementary proof in given in Appendix \ref{app-one-neg-eig}.
This result has implications for stability of the solitary waves as solutions of
KdV5 and it is discussed in \S\ref{stability-kdv5}.

\begin{figure}
\centering
\includegraphics[width=10cm]{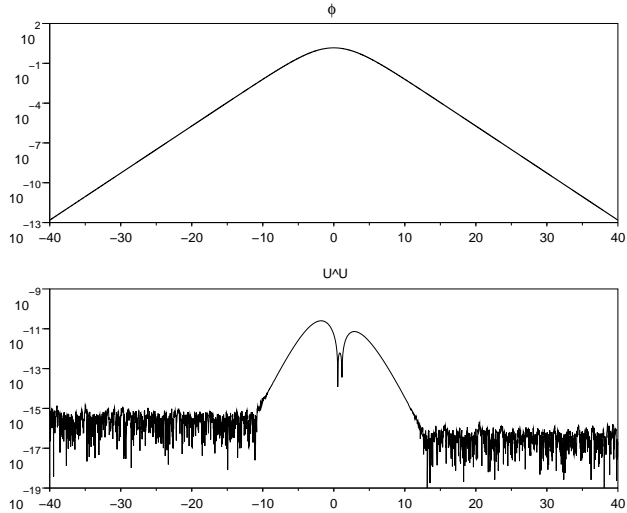}
\caption{ 
The decomposability of a $2$-form is equivalent to $I_1=0$ where
$I_1$ is defined in (\ref{def-I1I2}).  The upper figure shows the logarithm of 
the difference between  $a(x)$ and its limit $a_\infty$, and the lower
figure shows the logarithm of the value of $I_1$
when $q=1$, $P=\frac{13}{6}$, $\Delta x=.01$ and $\lambda=-10$
The fast and irregular oscillations are associated with round-off  errors.}
\label{rk4err}
\end{figure}

To test how accurately the Lagrangian Grassmannian is preserved
by the numerical scheme, the values of 
\begin{equation}\label{def-I1I2}
I_1 = {\bf U}\wedge{\bf U}\quad\mbox{and}\quad I_2 = \symp\wedge{\bf U}\,,
\end{equation}
are computed as a function of $x$.  In these calculations the standard explicit
fourth-order Runge-Kutta algorithm is used.  The value of $I_1$
is shown in Figure \ref{rk4err} and shows that the error is of order of the
machine accuracy, except for a small region around zero, but the error there is
still exceptionally small.
Concerning $I_2$, it is in fact exactly preserved, even numerically: 
if $I_2$ is below machine accuracy, then it remains so.

Using Proposition \ref{asymptotic} and the proof in Appendix
\ref{KdV_asymptotic}, it follows that 
the KdV5 system satisfies hypothesis
\ref{hyp_asymptotic} and therefore $D(\lambda)\to1$ and the Maslov index tends to $0$
as $\lambda\to-\infty$.

\section{Spectrum of $\mathscr{L}$ and the stability of solitary waves for KdV5}
\label{stability-kdv5}
\setcounter{equation}{0}

One of the intriguing properties of the Maslov index is the connection between
the number of eigenvalues in subsets of the $\lambda$ space, the Maslov index
and the stability of solitary waves.  We know that the Maslov index
counts the eigenvalues  of $\mathscr{L}$ and it turns out that
this number plays a role in previous stability results.

For KdV5 the connection between stability and  the Maslov index
is not obvious.  For unimodal solitary-wave solutions
of KdV5, \textsc{Kodama \& Pelinovsky}~\cite{kodama} have studied this connection
and they show the following result.
Suppose the following integral exists
\[
N(c,P)=\int_{\R} \widehat\phi(x,c,P)^2\rd x\,,
\]
and is a differentiable function of $c$, and define
\[
r=
\left\{\begin{array}{ll}
  0 &\text{ if }  \frac{\partial N}{\partial c}(1,P)<0\\

  1 &\text{ if }   \frac{\partial N}{\partial c}(1,P)>0\\
\end{array}\right.\,.
\]
The functional $N(c,P)$ is sometimes called the momentum of the
solitary wave.  If $q=1$, it is proved by \textsc{Lewandosky}~\cite{lewandosky} that
$r=1$ for all admissable $P$.

In \cite{kodama} it is argued (see proposition 3.8 there)
that a unimodal solitary wave is stable if $r=+1$ and $\#\mathscr{L}^-=1$,
where $\#\mathscr{L}^-$ is the number of negative eigenvalues of $\mathscr{L}$.

This observation is consistent with the theory
of \cite{bd-sima} where a instability results for a class of unimodal
solitary waves were obtained.

\subsection{Stability of two-pulse solitary waves}

More refined results on stability of two-pulse solitary waves were obtained
by \textsc{Chugunova \& Pelinovsky}~\cite{chug-peli-two-pulse}.
Suppose that ${\bf L}$ 
has only simple eigenvalues except a 
double eigenvalue at $0$ and suppose $\frac{\partial N}{\partial c}(1,P)\neq 0$.
In \cite{chug-peli-two-pulse} it is proved for symmetric solitary
wave solutions that
\[
 N_{unst} = \#\mathscr{L}^- - r - N^-_{imag}\,,
\]
where
$N_{unst}$ is the number of eigenvalues with 
strictly positive real part of ${\bf L}$ and 
$N^-_{imag}$ is the number of pure imaginary eigenvalues of ${\bf L}$ with 
negative Krein signature.  
Using the formula $\#\mathscr{L}^-=\textsf{Maslov}^{\rm homoclinic} - 1$
gives
\[
 N_{unst} = \textsf{Maslov}^{\rm homoclinic} - 1 - r - N^-_{imag}\,.
\]
It is immediate from this formula that
if $N^-_{imag}=0$ and $\textsf{Maslov}^{\rm homoclinic}=2$ then the basic state is stable
if $\frac{dN}{dc}>0$ and unstable if $\frac{dN}{dc}<0$.

Using the classification of \textsc{Buffoni, Champneys \& Toland}~\cite{bct}, a two-pulse
solitary wave has the classification ${\bf 2}(\ell)$ where $\ell$ is a natural number.
In \cite{cdb2008}, it is found numerically that two-pulse solutions have the following
formula for the Maslov index
\begin{equation}
\label{maslov-twopulse}
\textsf{Maslov}^{\rm homoclinic} = 
\left\{\begin{array}{ll}
  3 &\text{ if $\ell$ is even}\\
  4 &\text{ if $\ell$ is odd}\\
\end{array}\right.\,.
\end{equation}
We can make several observations using this formula for the Maslov index
of two-pulse homoclinic orbits.  Suppose $r=+1$, then
\[
N_{unst} = \textsf{Maslov}^{\rm homoclinic} - 2 - N_{imag}^- = 
\left\{\begin{array}{ll}
  1 - N_{imag}^- &\text{ if $\ell$ is even}\\
  2 - N_{imag}^- &\text{ if $\ell$ is odd}\\
\end{array}\right.
\,,
\]
We have the immediate observation that a necessary condition for a 
$2-$pulse homoclinic orbit to be stable is $N^-_{imag}>0$. From the
parity of $N^-_{imag}$, we have $N_{unst}=1$ and $N^-_{imag}=0$ when $l$ is even and the solitary
wave is unstable.
When $l$ is odd, the parity of $N^-_{imag}$ is not sufficient to determine
the stability of the solitary wave.

\textsc{Buryak \& Champneys}~\cite{bc-kdv-paper} used a completely
different method to study stability and they found that $2-$pulse solitary 
waves are stable if $\ell$ is odd and unstable if $\ell$ is even, 
assuming that $r=1$. This is consistent with the
value we found for the Maslov index.
When $l$ is odd, the stability of
the solitary wave is equivalent to $N_{imag}^-=2$.

The Maslov index does not give any information about the purely 
imaginary eigenvalues, and so to determine their number, 
a calculation of the spectral problem is necessary.  Some results
on this are reported by \textsc{Chardard}~\cite{chardard-thesis}.  
There it is found that  when $l$ is odd,
there are eigenvalues with non-zero imaginary part, 
but they appear to have very small
real parts.  
Further results are necessary to be certain about the 
spectral stability of $2-$pulse solitary waves when $l$ is odd.

One way to check whether the real part of a complex eigenvalue
is nonzero is to use the formula (for the case $q=1$)
\begin{equation}\label{realpart-lambda}
{\rm Re}(\lambda) = -\frac{1}{\|u\|^2}\int_{-\infty}^{+\infty}\widehat\phi_x|u(x)|^2\,\rd x\,.
\end{equation}
Here $\lambda$ is the complex eigenvalue associated with the stability exponent
and $u$ is the associated eigenmode:
\begin{equation}\label{spec-eqn}
\frac{d\ }{dx}\left( u_{xxxx} -Pu_{xx}+cu - 2\widehat\phi\,u\right) = \lambda u\,.
\end{equation}
The formula (\ref{realpart-lambda}) is derived by multiplying (\ref{spec-eqn})
by the complex conjugate of $u(\cdot)$ and integrating over $\R$:
$$\lambda \|u\|^2=\int_\R |u|^2\widehat{\phi_x}+\int_\R \widehat{\phi} u_x \overline{u}
=\int_\R |u|^2\widehat{\phi_x}
+\int_\R \widehat{\phi} \frac 1 2 (u_x \overline{u}-u \overline{u_x})
=\int_\R |u|^2\widehat{\phi_x}
+\i\int_\R \widehat{\phi} \mathrm{Im}(u_x \overline{u})$$

In spite of the simplicity of  the formula (\ref{realpart-lambda}) there is not
much that one can say in general.  
If $\widehat\phi(x)$  is an even function then $\widehat\phi_x$ is an odd function.
Then it is immediate that $|u(x)|^2$ even implies that ${\rm Re}(\lambda)=0$.
However, this is a highly special case.

\section{A model PDE for long-wave short-wave resonance}
\label{sec8}
\setcounter{equation}{0}

In this section the Maslov index is computed for a class of solitary
waves which arise in 
a model PDE for long-wave short-wave (LW-SW) resonance (cf.\ 
\textsc{Kawahara et al.}~\cite{ksk75},
\textsc{Ma}~\cite{ma79}, \textsc{Benilov \& Burtsev}~\cite{bb83},
 \textsc{Latifi \& Leon}~\cite{ll91}).  The LW-SW equations 
are a coupled system with one equation of nonlinear Schr\"odinger type
and the other of KdV type.  A typical form is
\begin{equation}
\label{eq:lwsw}
\begin{array}{rcl}
E_t &=& \ri (E_{xx} + \rho E - \nu E)\\
\rho_t &=& \partial_x(\rho_{xx} - c\rho + 3\rho^2 +|E|^2)\,,\\
\end{array}
\end{equation}
were $\rho(x,t)$ is real valued, $E(x,t)$ is complex valued, and
$c,\nu$ are considered to be positive real parameters.
In real coordinates,
$E = u+\ri v$ and $\rho = w$, the above equations can be written:
\begin{equation}
  \label{eq:real-lwsw}
\begin{array}{rcl}
u_t &=& -v_{xx} - vw +\nu v\\
v_t &=& u_{xx}  + uw - \nu u\\
w_t &=& w_{xxx}-cw_x + 6ww_x+2uu_x+2vv_x\,.
\end{array}  
\end{equation}
This system can be expressed as a Hamiltonian system in the time
direction.  However, we will not emphasize this property since it is
the spatial Hamiltonian structure that is associated with the Maslov index
(see Appendix \ref{app-LW-SW-Ham} for the temporal Hamiltonian structure). 
Solitary waves satisfy the steady equations
\begin{equation}\label{53}
\begin{array}{rcl}
-2u_{xx}  - 2uw + 2\nu u &=& 0 \\
-2v_{xx} - 2vw +2\nu v  &=&  0\\
-w_{xx}+cw - 3w^2-u^2-v^2 &=& \textsf{constant}\,,
\end{array}
\end{equation}
where the signs and coefficients are modified to ensure
that they are the Euler-Lagrange equation associated with the
Hamiltonian function $H(Z)$ in Appendix \ref{app-LW-SW-Ham}.
Exact solutions of this problem are known \cite{ma79}; for example,
\begin{equation}\label{Er-sol-1}
u(x) = A\,\text{sech}(\sqrt{\nu}\,x)\,,\quad v(x)=0\quad\mbox{and}\quad
w(x) = 2\nu\,\text{sech}^2(\sqrt{\nu}\,x)\,,\quad
\end{equation}
with $\textsf{constant}=0$ and $A^2 = 2\nu\,(c-4\nu)$, and
the existence condition $c-4\nu>0$.

To study the Maslov index of these solutions, linearize the steady
equations about the basic solitary wave and introduce a spectral parameter:
${\bf L}Z = \lambda Z$ with ${\bf L} = D^2H(\widehat Z)$.
Written out, this equation is
\begin{equation}\label{5}
\begin{array}{rcl}
-2u_{xx} - 2\widehat w u-2\widehat u w + 2\nu u &=& \lambda u \\
-2v_{xx} - 2\widehat w v-2\widehat v w + 2\nu v &=& \lambda v \\
-w_{xx} +cw - 6\widehat w w - 2\widehat u u - 2\widehat v v &=& 
\lambda w 
\end{array}
\end{equation}

When $\widehat v=0$ this system decouples into a second order equation
for $v$, and a fourth order coupled system for $u,w$,
\begin{equation}\label{5a}
\begin{array}{rcl}
-2u_{xx} - 2\widehat w u-2\widehat u w + 2\nu u &=& \lambda u \\
-w_{xx} +cw - 6\widehat w w - 2\widehat u u &=& 
\lambda w 
\end{array}
\end{equation}
The decoupled equation for $v$ is then
\begin{equation}\label{5b}
-2v_{xx} - 2\widehat w v + 2\nu v = \lambda v 
\end{equation}
This latter system can be analyzed completely and the result in given in
Appendix \ref{app-2times2}.

The fourth-order system for $u,w$ (\ref{5a}) can be written as a standard Hamiltonian ODE
in the form (\ref{h1}) with $n=2$ by taking
\[
{\bf u}(x,\lambda)=\left(
\begin{array}{c}
u \\
w \\
2u_x \\
w_x \\
\end{array}
\right)\,,\quad
{\bf B}(x,\lambda)=\left(
\begin{array}{cccc}
\lambda-2\nu+2\hat{w}(x) & 2\hat{u}(x) & 0 & 0 \\
2\hat{u} & \lambda-c+6\hat{w}(x)& 0 & 0 \\
0 & 0 & \frac{1}{2} & 0 \\
0 & 0 & 0 & 1 \\
\end{array}
\right)
\]
The essential spectrum for this equation is
\[
\sigma_{ess} =  \{\ \lambda \in \R\ :\ \lambda\geq 2\nu\ \mbox{and}\quad \lambda\geq c \}
\,.
\]
Adding the condition that $c>4\nu$, the system at infinity is hyperbolic for all
$\lambda\in\R$ such that $\lambda<2\nu$.

\begin{figure} 
\centering
  \includegraphics[width=14 cm]{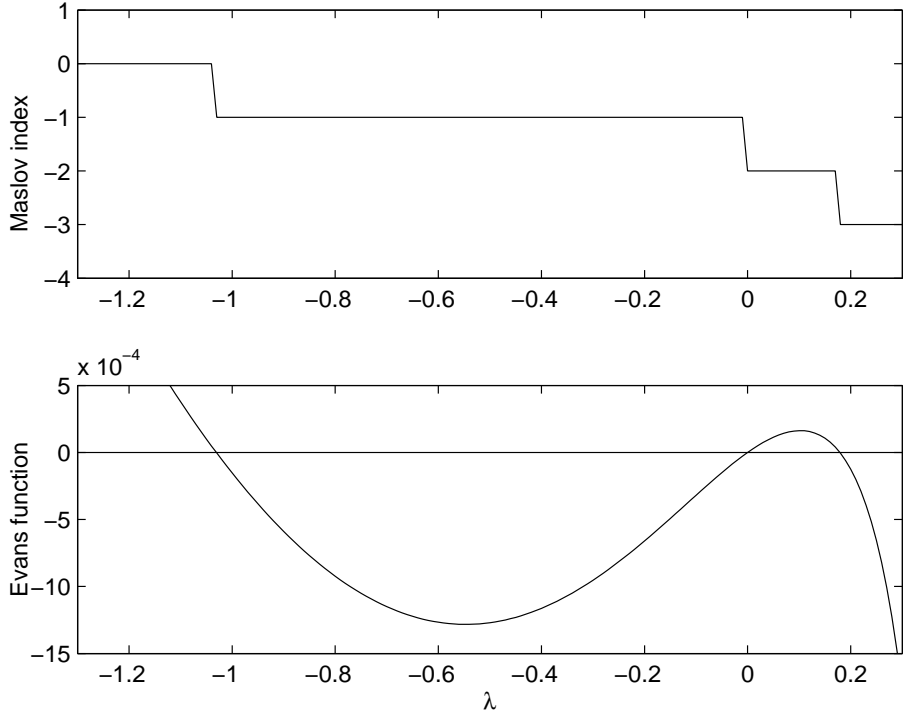}
\caption{Longwave-Shortwave problem for the following parameters $c=1$, $\nu=0.2$. 
Top: Maslov index as a function of $\lambda$. Bottom: Evans 
function as a function of $\lambda$. The Maslov index is $\lambda$-decreasing.}
\label{fig-lwsw-1}
\end{figure}

The Maslov index of the $4\times 4$ subsystem
is computed for the case $c=1$ and $\nu=0.2$ and the results
are shown, along with the associated Evans function, 
in Figure \ref{fig-lwsw-1} and tabulated
in the table below, where $\lambda_1<\lambda_2=0<\lambda_3$ are the three roots of
the Evans function.
\vspace{.25cm}

\begin{center}
\begin{tabular}{|c|c|c|c|c|}
\hline
\phantom{$\bigg|$}
$\lambda$ & $\lambda <\lambda_1$  & $\lambda_1<\lambda<\lambda_2$ & $\lambda_2<\lambda<\lambda_3$ & $\lambda>\lambda_3$ \\
\hline
\phantom{$\bigg|$}$\textsf{Maslov}({\bf U}^+,{\rm E}_\infty^s)$ & $0$ & $-1$ & $-2$ & $-3$ \\
\hline
\end{tabular}
\end{center} 

\section{A non-monotone Maslov index.}
\label{sec-10}
\setcounter{equation}{0}

In the LW-SW system and the KdV5 equation, the Maslov index is a monotone function of
$\lambda$. (Note however that it is not a monotone function of $x$.)  Here we show an
example where the Maslov index is not a monotone function of $\lambda$.  It is a slight
modification of the LW-SW resonance equations.  In this case the correlation between
the number of roots of the Evans function and the value of the Maslov index is no
longer apparent.  Look at the eigenvalue problem
\begin{equation}\label{non-monotone-L}
{\bf L}\begin{pmatrix}
  u\\
  w
\end{pmatrix}= \lambda\begin{pmatrix}
  u\\
  w
\end{pmatrix}\,,\quad \mbox{with}\quad
{\bf L}\begin{pmatrix}
  u\\
  w
\end{pmatrix} :=
\begin{pmatrix}
-2u_{xx} - 2\widehat w(x) u+2\widehat u(x) w + 2\nu u\\
 w_{xx} -cw + 6\widehat w(x) w + 2\widehat u(x) u
\end{pmatrix}\,,
\end{equation}
with $c>4\nu>0$,
\[\widehat u(x) = A\text{sech}(\sqrt{\nu}\,x) \quad\mbox{and}\quad
\widehat w(x) = 2\nu\,\text{sech}^2(\sqrt{\nu}\,x)
\]
with $A^2 = 2\nu(c-4\nu)$, and the requirement $c>4\nu>0$.

The spectral problem associated to this operator can be expressed in the form
(\ref{h1}) with $n=2$,
\[
{\bf u}(x)=\left(
\begin{array}{c}
u \\
w \\
2u_x \\
w_x \\
\end{array}
\right)\,,\quad\mbox{and}\quad
{\bf B}(x,\lambda) = \left(
\begin{array}{cccc}
\lambda - 2\nu + 2\hat{w}(x) & 2\hat{u}(x) & 0 & 0 \\
2\hat{u}(x) & -\lambda-c+6\hat{w}(x)& 0 & 0 \\
0 & 0 & \frac{1}{2} & 0 \\
0 & 0 & 0 & 1 \\
\end{array}
\right)\,.
\]
The essential spectrum of ${\bf L}$ consists of
\[
] -\infty,-c] \cup\ [2\nu,+\infty[.
\]
The essential spectrum is unbounded from above and below, hence a \emph{Morse index}
cannot be defined for ${\bf L}$.  However, we will still be able to compute a Maslov index.
The key property that leads to non-monotonicity is that the matrix
$\partial_\lambda{\bf B}(x,\lambda)$ is not semi-definite: the matrix
$\partial_\lambda{\bf B}(x,\lambda)$ has eigenvalues $\{0,0,-1,+1\}$ and so is not semi-definite.

\begin{figure} 
\begin{center}
\includegraphics[width=.80\columnwidth]{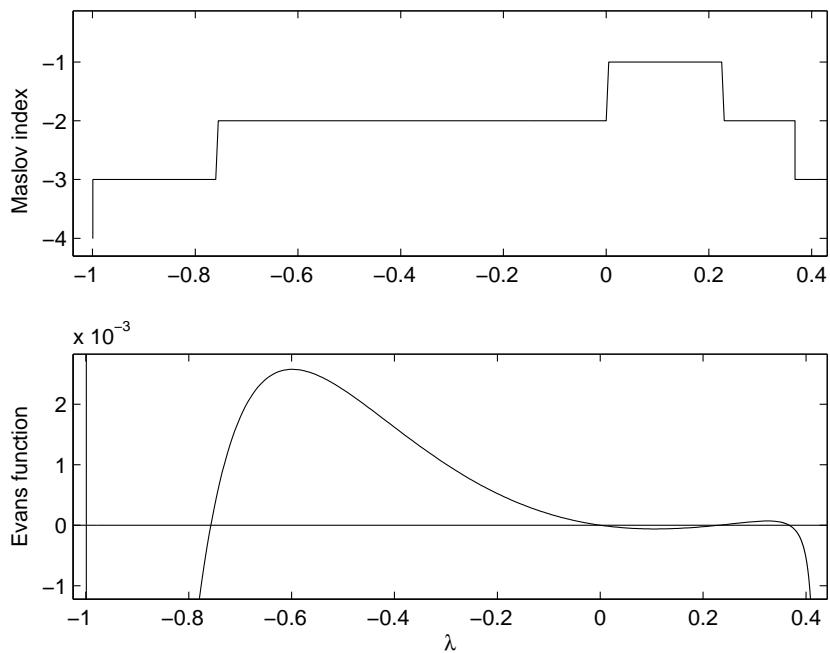}
\caption{Plot of the Maslov index for the non-monotone example
(\ref{non-monotone-L}) for the case $c=1$ and $\nu=0.21$.  Upper
figure shows the Maslov index and the lower figure the Evans function.} 
\label{fig-lwsw-3}
\end{center}
\end{figure}

Results for the case $c=1$ and $\nu=0.21$ are tabulated below and
shown in Figure \ref{fig-lwsw-3}.
In this case there are $5$ eigenvalues,
but there is no longer a correlation between the Maslov index and the
number of eigenvalues in a subset of $\lambda$.
\vspace{.25cm}

\begin{center}
\begin{tabular}{|c|ccccccccc|}
\hline
$\lambda$ & $-c$  & & & & $0$ & & & & $2\nu$\\
\hline
$D(\lambda)$ & $+\infty$& $+$ & $-$ & $+$ & $0$ & $-$ &
$+$ & $-$ & $-\infty$\\
\hline
$\textsf{Maslov}(\lambda)$ & & $-4$ & $-3$ & $-2$ & & $-1$ & $-2$ & $-3$ & \\
\hline
\end{tabular}
\end{center} 

\section{Concluding remarks}
\label{sec-conclusion}
\setcounter{equation}{0}

We have only just scratched the surface of the implications of the
Maslov index for homoclinic orbits and solitary waves.  Other
important questions for Hamiltonian systems on four-dimensional phase space
are (a) the connection between transversality of the homoclinic orbit
and the parity of the Maslov index, (b) the jump of the Maslov index
at bifurcations, (c) whether the angles $\kappa_1$ and $\kappa_2$ in the
decomposition $\kappa=\kappa_1+\kappa_2$ contain other useful information,
(d) the role of purely imaginary eigenvalues in the stability
of solitary waves in KdV5, and (e) the Maslov index of multi-pulse homoclinic orbits.
The latter question is addressed in the paper \cite{cdb2008}.

The extension of the Maslov index of homoclinic orbits to phase space with dimension
greater than four is straightforward in principle but there are some
differences in detail.  
First, the dimension of the basic manifolds  ($\bwedge^n(\R^{2n-1})$, $\rp^{C_{2n}^n-1}$, ${\rm G}_n(\R^{2n})$, $\Lambda(n)$, $\Lambda^1(n)$) jumps a lot
when $n$ goes from
$2$ to $3$.
\vspace{.15cm}

\begin{center}

\begin{tabular}{|c|c|c|c|c|c|}
\hline
\phantom{$\bigg|$}
\textsf{Manifold} & $\bwedge^3(\R^6)$ & $\rp^{19}$ & ${\rm G}_3(\R^6)$ & $\Lambda(3)$ & $\Lambda^1(3)$ \\
\hline
\phantom{$\bigg|$}\textsf{Dimension} & $20$ & $19$ & $9$ & $6$ & $5$ \\
\hline
\end{tabular}
\end{center} 
\vspace{.25cm}

In this table, $\rp^{19}$ is the projectification of $\bwedge^3(\R^6)$.  The biggest
change in the numerics is the difficulty due to the 
jump in dimension of the Lagrangian Grassmannian.  Whereas
it is $3-$dimensional in the case $n=2$, it jumps to double that dimension when
$n=3$.  
The details of the theory and numerics for $n\geq 3$ are given in
Part 2 \cite{cdb-part2}.

\newpage

\begin{center}
\hrule height.15 cm
\vspace{.2cm}
--- {\Large\bf Appendix} ---
\vspace{.2cm}
\hrule height.15cm
\end{center}
\vspace{.25cm}

\appendix

\renewcommand{\theequation}{A-\arabic{equation}}
\section{Kernel of ${\bf A}^{(2)}$ on $\bwedge^2(\R^4)$}
\label{A2-kernel}
\setcounter{equation}{0}

Let ${\bf A}$ be an arbitrary $4\times 4$ matrix with entries $a_{ij}$.
Then, with respect to the standard basis (\ref{E-basis}) 
on $\bwedge^2(\R^4)$ the induced matrix is
\begin{equation}\label{2.10}
{\bf A}^{(2)}=
\left[\begin{matrix}
a_{11} + a_{22} & a_{23} & a_{24} & -a_{13} & -a_{14} & 0 \\
a_{32}  & a_{11} + a_{33} & a_{34} & a_{12} & 0 & -a_{14} \\
 a_{42} & a_{43} & a_{11}+a_{44} & 0 & a_{12} & a_{13} \\
-a_{31} & a_{21} & 0 & a_{22}+a_{33} & a_{34} & -a_{24} \\
-a_{41} & 0 & a_{21} & a_{43} & a_{22}+a_{44} & a_{23} \\
0 & -a_{41} & a_{31} & -a_{42} & a_{32} & a_{33}+a_{44} 
\end{matrix}\right]\,,
\end{equation}
A constructive proof is given in \S2 of \cite{allen}.
\vspace{.15cm}

\begin{prop}
The induced matrix
${\bf A}^{(2)}$ satisfies 
${\bf A}^{(2)}\symp=0$, where $\symp$ is defined in (\ref{omega-def}),
if and only if ${\bf JA}$ is symmetric.
\end{prop}

\noindent{\bf Proof.} 
An explicit calculation gives
\[
{\bf A}^{(2)}\symp = 
\begin{pmatrix} a_{23} -a_{14} \\ a_{11}+a_{33} \\ a_{43}+a_{12} \\ a_{21}+a_{34}  \\ a_{22}+a_{44} \\ -a_{41} +a_{32}
\end{pmatrix}\,.
\]
On the other hand
\[
{\bf JA} = \left[\begin{matrix}
-a_{31} & - a_{32} & - a_{33} & -a_{34} \\
-a_{41} & - a_{42} & - a_{43} & -a_{44} \\
a_{11} &  a_{12} &  a_{31} & a_{14} \\
a_{21} &  a_{22} &  a_{23} & a_{24} \end{matrix}\right]\,.
\]
In order for ${\bf JA}$ to be symmetric we require
\[
a_{41} = a_{32} \,,\quad
a_{11} = -a_{33} \,,\quad
a_{12} = -a_{43} \,,\quad
a_{21} = -a_{34} \,,\quad
a_{22} = -a_{44} \,,\quad
a_{23} = a_{14}\,.
\]
These conditions are satisfied if and only if ${\bf A}^{(2)}\symp=0$.
$\hfill\blacksquare$

\renewcommand{\theequation}{B-\arabic{equation}}
\section{The connection between ${\bf A}^{(2)}$, ${\bf J}^{(2)}$ and ${\bf B}^{(2)}$}
\label{app_J2B2}
\setcounter{equation}{0}

In this appendix the role of the induced symplectic operator ${\bf J}^{(2)}$ on
$\bw$ is explored.
Using the standard formula for the induced matrix (\ref{2.10}),
the induced form of the symplectic operator is
\[
{\bf J}^{(2)} = \begin{pmatrix}
\hfill0 & \phantom{-}0 & -1      & \hfill1 & \phantom{-}0 & \hfill0 \\
\hfill0 & \phantom{-}0 & \hfill0 & \hfill0 & \phantom{-}0 & \hfill0 \\
\hfill1 & \phantom{-}0 & \hfill0 & \hfill0 & \phantom{-}0 & -1 \\
-1      & \phantom{-}0 & \hfill0 & \hfill0 & \phantom{-}0 & \hfill1 \\
\hfill0 & \phantom{-}0 & \hfill0 & \hfill0 & \phantom{-}0 & \hfill0 \\
\hfill0 & \phantom{-}0 & \hfill1 & -1      & \phantom{-}0 & \hfill0 \end{pmatrix} \,.
\]
Note that ${\bf J}^{(2)}\symp = 0$ and $\symp^T{\bf J}^{(2)} = 0$ using
$\symp$ defined in (\ref{omega-def}).  In fact, the kernel of ${\bf J}^{(2)}$ 
is four dimensional ${\rm Kernel}({\bf J}^{(2)}) = {\rm span}\{
{\bf E}_2,{\bf E}_5,{\bf E}_1+{\bf E}_6,{\bf E}_3+{\bf E}_4\}$,
where ${\bf E}_1,\ldots,{\bf E}_6$ is the standard basis for $\bwedge^2(\R^4)$.

Now suppose ${\bf B}$ is a symmetric matrix and
\[
{\bf A} = {\bf J}^{-1}{\bf B} = -{\bf JB}\,.
\]
then
\[
{\bf A}^{(2)} =
\left[\begin{matrix}
b_{13}+b_{24} & b_{34} & b_{44} & -b_{33} & -b_{34} & 0 \\
-b_{12}  & 0 & -b_{14} & b_{23} & 0 & -b_{34} \\
-b_{22} & -b_{23} & b_{13}-b_{24} & 0 & b_{23} & b_{33} \\
b_{11} & b_{14} & 0 & b_{24}-b_{13} & -b_{14} & -b_{44} \\
b_{12} & 0 & b_{14} & -b_{23} & 0 & b_{34} \\
0 & b_{12} & -b_{11} & b_{22} & -b_{12} & -b_{13}-b_{24} 
\end{matrix}\right]\,.
\]
The induced matrix ${\bf A}^{(2)}$ does not equal the product of the
induced matrices for ${\bf J}$ and ${\bf B}$ but it has the following form
\[
{\bf A}^{(2)} = -{\bf J}^{(2)}{\bf B}^{(2)} + {\bf S}\,,
\]
where ${\bf S}$ is the skew-symmetric matrix
\[
{\bf S} = \left[\begin{matrix}
0 & b_{12} & -b_{11} & b_{22} & -b_{12} & -b_{24}-b_{13} \\
-b_{12} & 0 & - b_{14} & b_{23} & 0 & -b_{34} \\
b_{11} & b_{14} & 0 & b_{24}-b_{13} & -b_{14} & -b_{44} \\
-b_{22} & -b_{23} & b_{13}-b_{24} & 0 & b_{23} & b_{33} \\
b_{12} & 0 & b_{14} & -b_{23} & 0 & b_{34} \\
b_{13}+b_{24} & b_{34} & b_{44} & -b_{33} & -b_{34} & 0 
\end{matrix}\right]\,.
\]
The skew-symmetric matrix ${\bf S}$ has the properties
${\bf S}\symp = 0$, $\symp^T{\bf S} = 0$, and ${\bf B}^{(2)}$
has the properties
\[
{\bf B}^{(2)}\symp\neq 0\quad\mbox{but}\quad {\bf B}^{(2)}\symp
\quad\mbox{is in the kernel of}\quad {\bf J}^{(2)}\,.
\]
Hence the property ${\bf A}^{(2)}\symp=0$ is recovered.  Moreover,
since $\symp^T{\bf J}^{(2)}=0$, we also have $\symp^T{\bf A}^{(2)}=0$.

\renewcommand{\theequation}{C-\arabic{equation}}
\section{Hyperbolic subspaces and Lagrangian planes}
\label{montaldi-lagrangian}
\setcounter{equation}{0}

\begin{figure}
\includegraphics[width=14cm]{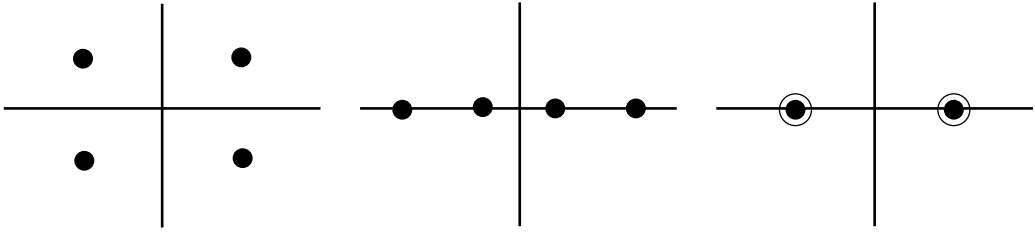}
\caption{ 
The three cases of hyperbolic spectra for constant coefficient
Hamiltonian systems on $\R^4$.}
\label{fig_hl}
\end{figure}
For a linear constant-coefficient Hamiltonian system on $\R^4$, 
${\bf u}_x={\bf A}{\bf u}$, there
are three cases where the spectrum of ${\bf A}$ is strictly hyperbolic and they are
shown in Figure \ref{fig_hl}.  The unstable subspace in each case is
a Lagrangian plane.  This observation is a special case of a result of
\textsc{Montaldi}~\cite{montaldi}.  Here we sketch a proof based on
eigenvectors.

Take the first case.  The unstable eigenvalues are of the form
$\mu = \nu \pm \ri \tau$ with $\nu>0$ and $\tau>0$, with eigenvectors
$\xi_1\pm\ri\xi_2$ where
\[
{\bf B}_\infty(\xi_1+\ri \xi_2) = (\nu+\ri \tau){\bf J} (\xi_1+\ri \xi_2)\,,
\]
or
\[
{\bf B}_\infty\xi_1 = \nu{\bf J}\xi_1 - \tau{\bf J}\xi_2\,,\quad
{\bf B}_\infty\xi_2 = \nu{\bf J}\xi_2 + \tau{\bf J}\xi_1 \,.
\]
Taking the inner product of the first equation with $\xi_2$ and
the second with $\xi_1$ and subtracting,
\[
0 = \xi_2^T{\bf B}_\infty\xi_1 - \xi_1^T{\bf B}_\infty\xi_2 =\nu\xi_2^T{\bf J}\xi_1 - \nu\xi_1^T{\bf J}\xi_2 =
2\nu \xi_2^T{\bf J}\xi_1 \,.
\]
Now since $\nu>0$ it follows that $\xi_2^T{\bf J}\xi_1=0$ and so
${\rm span}\{\xi_1,\xi_2\}$ is a Lagrangian subspace.

A similar calculation verifies the other two cases, with the third
requiring the introduction of generalized eigenspaces.

Now, use the fact that Lagrangian subspaces are invariant for the
$x-$dependent system, as shown in \S\ref{subsec-u-lagr},
to conclude that hyperbolic sets are Lagrangian manifolds.

\renewcommand{\theequation}{D-\arabic{equation}}
\section{The induced inner product on $\bwedge^k(\R^{2n})$}
\label{app-ip}
\setcounter{equation}{0}

In this appendix the equivalence between the induced inner product
$\lbk\cdot,\cdot\rbk_k$ on $\cwedgeknn$, which has dimension $d$,
and the standard inner product on $\R^d$ is established.  Here the
general case of $2n-$dimensional phase space is considered, which will
be required in Part 2 \cite{cdb-part2}.

With the standard orthonormal basis for $\R^{2n}$, 
$\{{\bf e}_1,\ldots,{\bf e}_{2n}\}$,
the nonzero and distinct members of the set 
\begin{equation}\label{2.3}
\{{\bf e}_{i_1}\wedge\cdots\wedge{\bf e}_{i_k}\ :\ i_1,\ldots,i_k
=1,\ldots, 2n\,\}
\end{equation}
form a basis for the vector space $\cwedgeknn$,
with exactly $d=\frac{n!}{(n-k)!k!}$ distinct elements.

Choose an ordering such as a standard lexical ordering and label
the nonzero distinct elements in the set (\ref{2.3})
by ${\bf E}_1,\ldots,{\bf E}_d$.
Then, any element ${\bf U}\in\cwedgeknn$ can be represented as
${\bf U} = \sum_{j=1}^d U_j\, {\bf E}_j$.
The inner product $\langle\cdot,\cdot\rangle$ on $\R^{2n}$
induces an inner product on each vector space $\cwedgeknn$ as
follows.  Let
\[
{\bf U} = {\bf u}_1\wedge\cdots\wedge{\bf u}_k
\quad{\rm and}\quad{\bf V} =
{\bf v}_1\wedge\cdots\wedge{\bf v}_k\,,\quad
{\bf u}_i,{\bf v}_j\in\R^{2n}\,,\quad\forall\ i,j=1,\ldots,k\,,
\]
be any decomposable $k$-forms.  A $k-$form is decomposable if
it can be written as a pure form: a wedge product between $k$
linearly independent vectors in $\R^{2n}$.  The inner product of
${\bf U}$ and ${\bf V}$ is defined by
\begin{equation}\label{ip-det-def}
\lbk {\bf U} , {\bf V} \rbk_k := {\rm det}\left[
\begin{matrix}
\langle {\bf u}_1,{\bf v}_1\rangle & \cdots & \langle{\bf u}_1,
{\bf v}_k\rangle\\
\vdots & \ddots & \vdots \\
\langle {\bf u}_k,{\bf v}_1\rangle & \cdots & \langle{\bf u}_k,
{\bf v}_k\rangle
\end{matrix}\right]\,,\quad
{\bf U},{\bf V}\in\cwedgeknn\,.
\end{equation}
Since every element in $\cwedgeknn$ is a sum of decomposable elements,
this definition extends by (bi)-linearity to any $k$-form.  Using the
orthonormality of the induced basis
\[
\lbk{\bf E}_i,{\bf E}_j\rbk_k = \left\{ \quad \begin{matrix} 1 & \mbox{if $i=j$}\\
 0 & \mbox{if $i\neq j$}\end{matrix}\right.\,,
\]
the inner product between two elements
${\bf U} = \sum_{i=1}^d U_i{\bf E}_i$ and ${\bf V} = \sum_{j=1}^d V_j{\bf E}_j$ is
\[
\begin{array}{rcl}
\lbk {\bf U},{\bf V}\rbk_k &=&  
\lbkbig\sum_{i=1}^d U_i{\bf E}_i,\sum_{k=1}^d V_j{\bf E}_j\rbkbig_k
= \sum_{i=1}^d\sum_{j=1}^d U_iV_j\lbk {\bf E}_i,{\bf E}_j\rbk_k\\[2mm]
&=& \sum_{i=1}^d U_iV_i := \langle {\bf U},{\bf V}\rangle_d\,,
\end{array}
\]
yielding the equivalent representation
\begin{equation}\label{ip-equiv}
\lbk {\bf U},{\bf V}\rbk_k = \langle {\bf U},{\bf V}\rangle_d\,,\quad
{\bf U},{\bf V}\in\cwedgeknn\,.
\end{equation}

\renewcommand{\theequation}{E-\arabic{equation}}
\section{Plus and minus subspaces in $\bwedge^2(\R^4)$}
\label{app-plus-minus}
\setcounter{equation}{0}

Consider $\R^4$ with the standard symplectic and volume forms.
Let $V={\rm span}\{\xi_1,\xi_2\}$ be a two-dimensional
\emph{oriented} subspace of $\R^4$. 

An oriented subspace is defined as follows. 
A subspace is an equivalence class of bases; that is,
${\rm span}\{\xi_1,\xi_2\}$ and ${\rm span}\{\eta_1,\eta_2\}$ 
represent the same subspace if and only if there is an invertible
$2\times 2$ matrix ${\bf m}$ such that $[\xi_1|\xi_2] =[\eta_1|\eta_2]{\bf m}$.
A subspace is oriented if ${\bf m}$ is restricted to have positive determinant.
The oriented subspace ${\rm span}\{\xi_1,\xi_2\}$ has one of three types
\[
\begin{array}{rcl}
\textsf{plus subspace} &\quad&
 \textsf{if}\ \omega\wedge\xi_1\wedge\xi_2 >0 \\[2mm]
\textsf{Lagrangian subspace} &\quad&
 \textsf{if}\ \omega\wedge\xi_1\wedge\xi_2 =0 \\[2mm]
\textsf{minus subspace} &\quad&
 \textsf{if}\ \omega\wedge\xi_1\wedge\xi_2 <0,.
\end{array}
\]

\begin{prop}
The sign of $\omega\wedge{\bf U}$ is an invariant of (\ref{E2}).
\end{prop}

Use Proposition \ref{prop-5-im} to conclude that
\[
\symp\wedge{\bf U}(x) = \textsf{constant}\,,
\]
along solutions of (\ref{E2}).

There is an interesting connection between Krein signature and the above classification
of oriented subspaces.  Krein signature is a sign which is associated with purely
imaginary eigenvalues (in the linearization about an equilibrium) or Floquet multipliers
(in the linearization about a periodic orbit).  Consider the case of a simple purely
imaginary eigenvalue $\ri \nu$, $\nu>0$.  
Its complex eigenvector $\zeta=\xi_1+\ri\xi_2$ satisfies
\[
{\bf B}\zeta = \ri\nu {\bf J}\zeta\,.
\]
The Krein signature is defined as the sign of
\[
\ri\langle{\bf J}\overline{\zeta},\zeta\rangle = 2\langle{\bf J}\xi_2,\xi_1\rangle\,.
\]
Now use the identity
$\langle  {\bf J}\xi_2 , {\bf v}\rangle = 
\lbk \symp,\xi_2\wedge{\bf v}\rbk_2 \,,\quad \forall\ {\bf v}
\in\R^{2n}\,$ to obtain

\[
\ri\langle{\bf J}\overline{\zeta},\zeta\rangle\vol =
2\langle \mathbf J\xi_2,\xi_1\rangle\vol = \lbk \symp, \xi_2\wedge\xi_1\rbk_2\vol
= -\lbk \symp, \xi_1\wedge\xi_2\rbk_2\vol = \symp\wedge\xi_1\wedge\xi_2\,.
\]
This observation also emphasizes the fact that a choice of orientation underlies the
definition of Krein signature.

\renewcommand{\theequation}{F-\arabic{equation}}
\section{A spectral problem with ${\rm sech}^2$ potential}
\label{app_spectral}
\setcounter{equation}{0}

This appendix establishes the basic properties of the ODE
eigenvalue problem
\begin{equation}\label{a1}
\phi_{xx} + 12\,{\rm sech}^2x\,\phi = \kappa\,\phi\,,
\end{equation}
in the set $\mathscr{K}:=\{\kappa\in\R\ |\ \kappa>0\}$.
The solutions of this ODE can be determined explicitly.  The
eigenvalues are $\kappa=1,4,9$.
For all $\kappa\in\mathscr{K}\setminus\{1,4,9\}$, the two functions
\begin{equation}\label{phi-pm}
\phi^{\pm}(x,\kappa) = {\rm e}^{\pm\sqrt{\kappa}\,x}\left(
\pm a_0 + a_1\,{\rm tanh}(x) \pm a_2\,{\rm tanh}^2(x) + {\rm tanh}^3(x)
 \right)\,,
\end{equation}
are linearly independent, where
\[
a_0 = \frac{\sqrt{\kappa}}{15}(4-\kappa)\,,\quad
a_1 = \frac{1}{5}(2\kappa-3)\,,\quad a_2 = -\sqrt{\kappa}\,.
\]

The eigenvalues can be verified by explicit calculation.
That $\kappa=1,4,9$ are the only eigenvalues in $\mathscr{K}$, and that
$\phi^{\pm}$ are linearly independent is verified by computing the
Wronskian
\[
W(x,\kappa) = {\rm det}\left[\begin{matrix}
\phi^+ & \phi^- \\ \phi^+_x & \phi^-_x \end{matrix}\right]\,.
\]
It is easily verified that $W_x=0$ and so $W(x,\kappa)$ is independent
of $x$.  Evaluate at $x=0$
\[
W(0,\kappa) =
{\rm det}\left[\begin{matrix}
a_0 & -a_0 \\
a_1 + a_0\sqrt{\kappa} & a_1 + a_0\sqrt{\kappa} \end{matrix}\right]
= 2 a_0(a_1 + a_0\sqrt{\kappa})\,.
\]
Substituting for $a_0$ and $a_1$,
\[
W(0,\kappa) = \frac{2}{225}\sqrt{\kappa}(\kappa-1)(\kappa-4)(\kappa-9)\,.
\]
Hence $\phi^{\pm}$ are linearly independent for all
$\kappa\in\mathscr{K}\setminus\{1,4,9\}$.

The eigenfunctions are 
\[
\begin{array}{rcl}
\phi(x,1) &=& {\rm sech}(x)( 4 - 5\,{\rm sech}^2(x))\,,\quad\mbox{when}
\quad\kappa=1\\[2mm]
\phi(x,4) &=& {\rm tanh}(x)\,{\rm sech}^2(x)\,,\quad\mbox{when}
\quad\kappa=4\\[2mm]
\phi(x,9) &=& {\rm sech}^3(x)\,,\quad\mbox{when}
\quad\kappa=9\,.
\end{array}
\]
modulo an arbitrary multiplicative constant.

\renewcommand{\theequation}{G-\arabic{equation}}
\section{Attractivity of the Lagrangian Grassmannian $\Lambda(2)$}
\label{sec_attractivity}
\setcounter{equation}{0}

One of the advantages of subtracting off the growth rate at infinity 
in the equations on $\bwedge^2(\R^4)$, as in
(\ref{key}), is that the Lagrangian Grassmannian
becomes an attracting invariant manifold.  When $\Lambda(2)$ is
attractive, one has
greater freedom in choosing the numerical integration scheme.

To prove attractivity, consider the integration of the $2-$form representing
the unstable subspace ${\bf U}^+(x,\lambda)$
\[
\frac{d\ }{dx}{\bf U}^+ = {\bf A}^{(2)}(x,\lambda){\bf U}^+\quad
{\bf U}^+\in\bwedge^2(\R^4)\quad -L<x<+L\,.
\]
Introduce the transformation
\[
{\bf U}^+(x,\lambda) = {\rm e}^{\sigma_+(\lambda)x}\,
\widehat{\bf U}^+(x,\lambda)
\]
where $\sigma_+(\lambda)$ is the sum of the eigenvalues of
${\bf A}_\infty(\lambda)$ with positive real part.
Then $\widehat{\bf U}^+$ satisfies
\begin{equation}\label{e-1}
\frac{d\ }{dx}\widehat{\bf U}^+ = [{\bf A}^{(2)}(x,\lambda)
-\sigma_+(\lambda){\bf I}]
\widehat{\bf U}^+\quad
\quad -L<x<+L
\end{equation}

The Lagrangian Grassmannian is the set
\[
 \widehat{\bf U}^+\wedge \widehat{\bf U}^+=0 \quad\mbox{and}\quad
\symp\wedge \widehat{\bf U}^+=0\,.
\]
When evaluated on the differential equation (\ref{e-1}) these invariants
satisfy
\[
\begin{array}{rcl}
\frac{d\ }{dx} \widehat{\bf U}^+\wedge \widehat{\bf U}^+ &=& \frac{d\ }{dx}\widehat{\bf U}^+\wedge \widehat{\bf U}^+
+ \widehat{\bf U}^+\wedge \frac{d\ }{dx}\widehat{\bf U}^+\\[2mm]
&=& {\bf A}^{(2)}\widehat{\bf U}^+\wedge \widehat{\bf U}^+ +
\widehat{\bf U}^+\wedge {\bf A}^{(2)}\widehat{\bf U}^+
-2\sigma_+\,\widehat{\bf U}^+\wedge\widehat{\bf U}^+\\[2mm]
&=& {\rm Trace}({\bf A})\,\widehat{\bf U}^+\wedge \widehat{\bf U}^+ 
-2\sigma_+\,\widehat{\bf U}^+\wedge\widehat{\bf U}^+\\[2mm]
&=& -2\sigma_+\, \widehat{\bf U}^+\wedge \widehat{\bf U}^+\,,
\end{array}
\]
since ${\rm Trace}({\bf A})=0$.
A similar calculation with $\symp\wedge \widehat{\bf U}^+$ yields
\[
\frac{d\ }{dx}\symp\wedge \widehat{\bf U}^+ = \symp \wedge ({\bf A}^{(2)}-\sigma_+{\bf I}){\bf U}^+
= - {\bf A}^{(2)}\symp\wedge{\bf U}^+ - \sigma_+\,\symp\wedge \widehat{\bf U}^+ = -\sigma_+\,\symp\wedge \widehat{\bf U}^+\,,
\]
using the fact that ${\bf A}^{(2)}\symp=0$, which is proved in Appendix \ref{A2-kernel}.,
This proves that
\[
\widehat{\bf U}^+\wedge \widehat{\bf U}^+(x) = {\rm e}^{-2\sigma_+x}\widehat{\bf U}^+\wedge \widehat{\bf U}^+\bigg|_{x=-L}\quad\mbox{and}\quad
\symp\wedge \widehat{\bf U}^+(x) = {\rm e}^{-\sigma_+x}\symp\wedge \widehat{\bf U}^+\bigg|_{x=-L}\,,\quad\mbox{for}\quad x>-L\,.
\]
The eigenvalue $\sigma_+$ is real and positive.  
Hence when integrating the unstable subspace ${\bf U}^+$ along
the Lagrangian Grassmannian, both $\widehat{\bf U}^+\wedge \widehat{\bf U}^+$
and $\symp\wedge \widehat{\bf U}^+(x)$ are 
exponentially attracted to the zero set.  Therefore a special
integrator is not required for maintaining the constraints;
a standard Runge-Kutta algorithm is quite satisfactory.

\renewcommand{\theequation}{H-\arabic{equation}}
\section{The existence of at least one negative eigenvalue}
\label{app-one-neg-eig}
\setcounter{equation}{0}

Consider the linear operator
\begin{equation}\label{L1}
\mathscr{L}\phi := \phi_{xxxx} - P\phi_{xx} + a(x) \phi\,,
\end{equation}
introduced in (\ref{4d-example}) with $a(x) = c - (q+1)\widehat\phi(x)^q$
and $\widehat\phi(x)$ satisfying (\ref{steady_kdv}).
Assume
\begin{equation}\label{hyp-L2}
P+2c\geq0\quad\mbox{and}\quad 0 < c\leq 1\quad\mbox{or}\quad P>0\quad\mbox{and}\quad
c>0\,.
\end{equation}
The essential spectrum for this
problem is non-negative.  Here it is proved that $\mathscr{L}$ has at least one
negative eigenvalue in the point spectrum. 

Multiply (\ref{steady_kdv}) by the basic state $\widehat\phi(x)$,
\[
\begin{array}{rcl}
\widehat\phi^{q+2} &=& c\widehat\phi^2 - P\,\widehat\phi\widehat\phi_{xx} + \widehat\phi\widehat\phi_{xxxx}\\
&=& c\widehat\phi^2 - (P+2c)\,\widehat\phi\widehat\phi_{xx} +2c\widehat\phi\widehat\phi_{xx} + \widehat\phi\widehat\phi_{xxxx}\\
&=& c(\widehat\phi+\widehat\phi_{xx})^2 - (P+2c)\,\widehat\phi\widehat\phi_{xx} - c\widehat\phi_{xx}^2 + \widehat\phi\widehat\phi_{xxxx}\,.
\end{array}
\]
Hence integrating, using the fact that $\widehat\phi$ and its derivatives decay
exponentially as $x\to\pm\infty$, and the hypotheses (\ref{hyp-L2})
\begin{equation}\label{q2-positive}
\int_{-\infty}^{\infty} \widehat\phi^{q+2}\rd x =
\int_{-\infty}^{\infty} c(\widehat\phi+\widehat\phi_{xx})^2\rd x + (P+2c)
\int_{-\infty}^{\infty} \widehat\phi_x^2\rd x +(1-c)\int_{-\infty}^{+\infty}\widehat\phi_{xx}^2\,\rd x > 0\,.
\end{equation}
or if $P>0$ and $c>0$,
\begin{equation}\label{q2-formula}
\int_{-\infty}^{\infty} \widehat\phi^{q+2}\rd x =
\int_{-\infty}^{\infty} (c\widehat\phi^2 + P\widehat\phi_x^2  + \widehat\phi_{xx}^2)\,\rd x
>0\,.
\end{equation}

To prove that (\ref{L1}) has a negative eigenvalue, we
will show that the quadratic form $\langle u,\mathscr{L}u\rangle$
is negative when $u=\widehat\phi$ where
$\langle u,\mathscr{L}u\rangle := \int_{-\infty}^{\infty} 
u\,\mathscr{L} u\,\rd x$.  Now
\[
\begin{array}{rcl}
\langle u,\mathscr{L}u\rangle\big|_{u=\widehat\phi} &=& 
\int_{-\infty}^{\infty} (\widehat\phi(\widehat\phi_{xxxx}-P\widehat\phi_{xx}+a(x)\widehat\phi))\,\rd x\\[2mm]
&=&
\int_{-\infty}^{\infty}(\widehat\phi_{xx}^2+P\,\widehat\phi_x^2+c\widehat\phi^2)\,\rd x 
- (q+1)\int_{-\infty}^{\infty}\widehat\phi^{q+2}\rd x \\[2mm]
&=& - q\,\int_{-\infty}^{\infty}\widehat\phi^{q+2}\rd x\,,
\end{array}
\]
using (\ref{q2-formula}) in the last line.  It follows from (\ref{q2-positive}) 
or (\ref{q2-formula}) that
$\langle\widehat\phi,\mathscr{L}\widehat\phi\rangle <0$.

\renewcommand{\theequation}{I-\arabic{equation}}
\section{Check of hypothesis \ref{hyp_asymptotic} for the KdV5 system}
\label{KdV_asymptotic}
\setcounter{equation}{0}

In this appendix, the details are given of the proof
that ${\bf A}(x,\lambda)={\bf J}^{-1}{\bf B}(x,\lambda)$ for KdV5,
with ${\bf B}(x,\lambda)$ defined in (\ref{u-B-def}),
satisfies Hypothesis \ref{hyp_asymptotic}.

First set
\[
s=\frac{1}{ (1-\lambda)^{\frac 1 4}}\,.
\]
When $\lambda$ is large and negative, $s$ is a small parameter.
This parameter will be used to obtain series expansions of 
the eigenvalues and eigenvectors.

The characteristic polynomial of $\mathbf A_{\infty}(\lambda)$ is
\[
0={\rm det}[X{\bf I} - {\bf A}_\infty(\lambda)] = X^4-P X^2+\frac 1 {s^4}\,.
\]
This polynomial is a biquadratic and for $s$ small 
it has four complex roots, one in each quadrant. 
Let $\theta(s)$ be the eigenvalue in the right-upper quadrant.
Its Taylor expansion is: 
\[
\theta(s)= \frac{1}{s\sqrt{2}}  
\left( 1-\frac 1 4\,\ri P{s}^{2}-\frac 1 {32} {P}^{2}{s}^{4}
-{\frac 1 {128}}\,\ri{P}^{3}{s}^{6}+ O\left(s^8\right)
 \right) \,.
\]
The other eigenvalues are $\overline{\theta(s)}$,
$-\theta(s)$,$-\overline{\theta(s)}$.
The eigenvector associated with $\theta(s)$ is :
\[
{\bf v}(s)=\frac 1 s\begin{pmatrix}
 -{\frac {1}{{s}^{4}\theta}}\\
\theta\\
1\\
{\theta}^{2}  
\end{pmatrix}\,.
\]
A Taylor expansion of this eigenvector is:
\[
{\bf v}(s)=s{\bf v}_1+\mathrm i s {\bf v}_2+ O\left(s^6\right)\,.
\]
with
\[
{\bf v}_1(s) =\begin{pmatrix}
- \left( 1/2\,\sqrt {2}+\frac 1 8\,\sqrt {2}P{s}^{2}-{\frac {1}{64}}\,\sqrt {2}{P}^{2}{s}^{4}+{\frac {1}{256}}\,\sqrt {2}{P}^{3}{s}^{6} \right) {s}^{-3}\\
{\frac {1}{256}}\,\sqrt {2}{P}^{3}{s}^{5}-{\frac {1}{64}}\,\sqrt {2}{P}^{2}{s}^{3}+\frac 1 8\,\sqrt {2}Ps+{\frac {\sqrt {2}}{2s}}\\
1\\
\frac 1 2 P
\end{pmatrix}
\]
and
\[
{\bf v}_2(s) =\begin{pmatrix}
 - \left( -1/2\,\sqrt {2}+\frac 1 8\,\sqrt {2}P{s}^{2}+{\frac {1}{64}}\,\sqrt {2}{P}^{2}{s}^{4}+{\frac {1}{256}}\,\sqrt {2}{P}^{3}{s}^{6} \right) {s}^{-3}\\
-{\frac {1}{256}}\,\sqrt {2}{P}^{3}{s}^{5}-{\frac {1}{64}}\sqrt {2}{P}^{2}{s}^{3}-\frac 1 8\sqrt {2}Ps+{\frac {\sqrt {2}}{2s}}\\
0\\
{s}^{-2}-\frac 1 8{s}^{2}{P}^{2} 
\end{pmatrix}\,.
\]
$(\Re v(s),\Im v(s))$ is a basis of the unstable space.
Let $V_{unst}$ the matrix whose columns are $\Re v(s)$ 
and $\Im v(s))$.

The eigenvector associated to $-\theta(s)$ is:
\[
{\bf w}(s)=s\begin{pmatrix}
 {\frac {1}{{s}^{4}\theta}}\\
-\theta\\
1\\
{\theta}^{2}  
\end{pmatrix}
\]
$(\Re {\bf w}(s),\Im {\bf w}(s))$ is a basis of the unstable space.
Let $U_{st}$ the matrix whose columns are $\Re {\bf w}(s)$ 
and $\Im {\bf w}(s))$.

The matrix $\begin{pmatrix}
  V_{unst} | U_{st}
\end{pmatrix}$ is not  a symplectic matrix but
$V(s)=\begin{pmatrix}
  V_{unst} |V_{st} 
\end{pmatrix}$, with $V_{st}=-U_{st} (V_{unst}^TJU_{st})^{-1}$, is.
Besides, we have:

$$V_{st}=\begin{pmatrix} 
\frac 1 {2s}&{\frac {-\frac 1 4\,{s}^{2}P-\frac 1 {32}\,{P}^{3}{s}^{6}}{s}}\\
0&{\frac {\frac 1 2\,{s}^{2}+\frac 1 {16}\,{s}^{6}{P}^{2}}{s}}\\
{\frac {\frac 1 4\,\sqrt {2}{s}^{3}-\frac 1 {16}\,\sqrt {2}P{s}^{5}}{s}}&{\frac {-\frac 1 4\,\sqrt {2}{s}^{3}-\frac 1 {16}\,\sqrt {2}P{s}^{5}}{s}}\\
 \frac{ -\frac 1 4\,\sqrt {2}s+\frac 1 {16}\,\sqrt {2}P{s}^{3}-{\frac {3}{128}}\,\sqrt {2}{P}^{2}{s}^{5}} s & \frac{ -\frac 1 4\,\sqrt {2}s-\frac 1 {16}\,\sqrt {2}P{s}^{3}-{\frac {3}{128}}\,\sqrt {2}{P}^{2}{s}^{5} } s
\end{pmatrix}+O(s^6)
$$
We also have $V^{-1}=-J( ^TV) J$ since $V$ is a symplectic matrix.

Let 
\[
{\bf B}=\begin{pmatrix}
   0&0&1&0\\
0&0&0&0\\
0&0&0&0\\
0&0&0&0
\end{pmatrix}\,.
\]

We are now able to evaluate:
$V^{-1}BV$:
$V^{-1}BV=O\left(s^2\right)$
but also
\[
V^{-1}BV
\begin{pmatrix}
  1 & 0\\
  0 & 1\\
  0 & 0\\
  0 & 0
\end{pmatrix}
=VBV_{unst}=  O\left(s^2\right)\,.
\]
Therefore, $(V^{-1}BV)^{(2)} e_1=  O\left(s^2\right)$.
Hence, as $R(x,\lambda)=A(x,\lambda)-A_{\infty}(\lambda)=(1-a(x))B$ and as 
$|1-a(x)|\leq C_1 {\rm e}^{-C_2 |x|}$,  this proves that
Hypothesis \ref{hyp_asymptotic} of 
Proposition \ref{asymptotic} is satisfied.

\renewcommand{\theequation}{J-\arabic{equation}}
\section{Hamiltonian evolution equation for LW-SW equations}
\label{app-LW-SW-Ham}
\setcounter{equation}{0}

The LW-SW equations (\ref{eq:real-lwsw}) can be expressed in Hamiltonian form as follows.
Let
\[
{\bf K} = \left[\begin{matrix}
0 & \fr & 0 \\ -\fr & 0 & 0 \\ 0 & 0 & -\partial_x 
\end{matrix}\right]\,,
\]
\[
H(Z) = \int\limits_{-\infty}^{+\infty}
\left( u_x^2+v_x^2+\frac{1}{2} w_x^2 - w(w^2+u^2+v^2)  
+ \frac{1}{2} cw^2 +\nu(u^2+v^2) \right) \d x,
\]
with $Z=(u,v,w)$.  Then the system becomes
\[
Z_t = {\bf K}\nabla H(Z)\,,
\]
since, with respect to an $L_2(\R)$ inner product,
\[
\nabla H(Z) = \begin{pmatrix} H_u\\H_v\\H_w\end{pmatrix} =
\begin{pmatrix}  -2u_{xx}  - 2uw + 2\nu u \\ 
-2v_{xx} - 2vw +2\nu v \\
-w_{xx}+cw - 3w^2-u^2-v^2
\end{pmatrix}
\]

\renewcommand{\theequation}{K-\arabic{equation}}
\section{The reduced eigenvalue problem associated with LW-SW equations}
\label{app-2times2}
\setcounter{equation}{0}

The two-dimensional ODE (\ref{5b}) that arises in the reduced problem for LW-SW resonance
can be written in the form
\[
v_{xx} + 2\nu\,{\rm sech}^2(\sqrt{\nu}x)\, v = \left(\nu -\frac{1}{2}\lambda\right) v\,.
\]
ODEs of this type can be solved explicitly as noted in Appendix \ref{app_spectral}.
The essential spectrum is the semi-infinite interval
$\sigma_{ess}(L)=[2\nu,+\infty)$.  Now suppose that $\lambda < 2\nu$.
Then the system at infinity is hyperbolic and one can explicitly construct the
solutions $(v^+,v^-)$ which give the solutions for the stable and unstable subspace
\[
v^{\pm}(x;\lambda)=e^{\pm\mu\sqrt{\nu} x}(\mp\mu+\tanh(\sqrt{\nu} x))\,,\quad
\mu=\sqrt{1-\frac{\lambda}{2\nu}}\,.
\]
The Evans function can be obtained from
\[
D(\lambda) = \left.{\rm det}\left[ \begin{matrix}
v^+(x;\lambda) & v^-(x;\lambda) \\
v_x^+(x;\lambda) & v_x^-(x;\lambda) \end{matrix}\right]\right|_{x=0} = 2\mu\sqrt{\nu}(\mu^2-1)=
-\frac{\lambda}{\sqrt{\nu}}\sqrt{ 1 - \frac{\lambda}{2\nu}}\,.
\]
The Maslov index is 
\[
\textsf{Maslov}(\lambda)=\begin{cases}
1\quad \text{ if } \lambda<0 \\
0\quad \text{ if } 0< \lambda< 2\nu \\
\end{cases}\,.
\]



\begin{thebibliography}{}

\bibitem{agj}
\textsc{J.W. Alexander, R. Gardner \& C.K.R.T. Jones}.
{\it A topological invariant arising in the stability analysis of 
traveling waves}, J. Reine Angew. Math. {\bf 410} 167--212 (1990).

\bibitem{allen}
\textsc{L. Allen \& T.J. Bridges}.
{\it Numerical exterior algebra and the compound matrix method},
Numerische Mathematik {\bf 92}, 197--232 (2002).

\bibitem{arnold1}
\textsc{V.I. Arnol'd}.
{\it On a characteristic class entering into conditions of quantization},
Funct. Anal. Appl. {\bf 1} 1--14 (1967).

\bibitem{arnold2}
\textsc{V.I. Arnol'd}.
{\it The Sturm theorems and symplectic geometry},
 Funct. Anal. Appl. {\bf 19} 251--259 (1986).

\bibitem{bb83}
\textsc{E.S. Benilov \& S.P. Burtsev}.
{\it To the integrability of the equations describing the langmuir-wave
ion-acoustic wave interaction}, Phys. Lett. A {\bf 98} 256--258 (1983).

\bibitem{bose}
\textsc{A. Bose \& C.K.R.T. Jones}.
{\it Stability of the in-phase travelling wave solution in a pair of 
coupled nerve fibers},
Indiana Univ. Math. J. \textbf{44}, 189--220 (1995). 

\bibitem{bd-sima}
\textsc{T.J. Bridges \& G. Derks}.
{\it Linear instability of solitary wave solutions of the Kawahara equation and
its generalizations}, SIAM J. Math. Anal. {\bf 33} 1356--1378 (2002).

\bibitem{bd-analyticity}
\textsc{T.J. Bridges \& G. Derks}.
{\it Constructing the symplectic Evans matrix using maximally analytic
individual vectors}, Proc. Roy. Soc. Edin. A {\bf 133} 505--526 (2003).

\bibitem{bdg02}
\textsc{T.J. Bridges, G. Derks \& G. Gottwald}.
{\it Stability and instability of solitary waves of the fifth-order KdV
equation: a numerical framework},
Physica D {\bf 172}, 190--216 (2002).

\bibitem{bd05}
\textsc{T.J. Bridges \& N.M. Donaldson}.
{\it Degenerate periodic orbits and homoclinic torus bifurcation},
Phys. Rev. Lett. {\bf 95}(10) 104301 (2005).

\bibitem{bct}
\textsc{B. Buffoni, A.R. Champneys \& J.F. Toland}.
{\it Bifurcation and coalescence of a plethora of multi-modal homoclinic orbits
in a Hamiltonian system},
J. Dyn. Diff. Eqns. {\bf 8} 221--281 (1996).

\bibitem{bc-kdv-paper}
\textsc{A.V. Buryak \& A.R. Champneys}.
{\it On the stability of solitary wave solutions of the fifth-order KdV
equation}, Phys. Lett. A {\bf 233} 58--62 (1997).

\bibitem{champ}
\textsc{A.R. Champneys}.
{\it Homoclinic orbits in reversible systems and their applications in
mechanics, fluids and optics},
Physica D {\bf 112}, 158--186 (1999).

\bibitem{F.Chardard2007}
\textsc{F. Chardard}.
{\it Maslov index for solitary waves obtained as a limit
 of the Maslov index for periodic waves},
C. R. Acad. Sci. Paris, Ser. I \textbf{345} 689--694 (2007).

\bibitem{chardard-thesis}
\textsc{F. Chardard}.
{\it Stabilit\'e des ondes solitaires},
Th\`ese de Doctorat de L'\'Ecole Normale Sup\'erieure de Cachan (2009).

\bibitem{F.Chardard2006}
\textsc {F. Chardard, F. Dias \& T.J. Bridges}. 
{\it Fast computation of the Maslov index for hyperbolic linear systems
with periodic coefficients}, 
J. Phys. A: Math. Gen. \textbf{39} 14545--14557 (2006).

\bibitem{cdb2008}
\textsc {F. Chardard, F. Dias \& T.J. Bridges}.
{\it On the Maslov index of multi-pulse homoclinic orbits},
Preprint (2008).

\bibitem{cdb-part2}
\textsc {F. Chardard, F. Dias \& T.J. Bridges}.
{\it Computing the Maslov index of solitary waves. Part 2. Hamiltonian
systems on a $2n-$dimensional phase space}, Preprint (2008).

\bibitem{chen_hu}
\textsc{C.-N. Chen \& X. Hu}.
{\it Maslov index for homoclinic orbits of {H}amiltonian systems}, 
Ann. Inst. H. Poincar\'e Anal. Non Lin\'eaire
{\bf 24} 589--603 (2007). 

\bibitem{chug-peli-two-pulse}
\textsc{M. Chugunova \& D. Pelinovsky}.
{\it Two-pulse solutions in the fifth-order {K}d{V} equation: 
rigorous theory and numerical approximations},
Discrete Contin. Dyn. Syst. Ser. B \textbf{8} 773--800 (2007).

\bibitem{crl}
\textsc{S.C. Creagh, J.M. Robbins \& R.G. Littlejohn}.
{\it Geometrical properties of Maslov indices in the semiclassical trace
formula for the density of states}, Phys. Rev. A {\bf 42} 1907--1922 (1990).

\bibitem{diasiooss}
\textsc{F. Dias \& G. Iooss}.
{\it Water-waves as a spatial dynamical system},
Handbook of Mathematical Fluid Dynamics {\bf 2},
Elsevier Science: Amsterdam (2003).

\bibitem{diaskuz}
\textsc{F. Dias \& E.A. Kuznetsov}.
{\it On the non-linear stability of solitary wave solutions of the fifth-order
{K}orteweg-de {V}ries equation},
Physics Letters A {\bf 263}, 98--104 (1999).

\bibitem{dmv}
\textsc{F. Dias, D. Menasce \& J.-M. Vanden-Broeck}.
{\it Numerical study of capillary-gravity solitary waves},
Eur. J. Mech. B/Fluids {\bf 15} 17--36 (1996).

\bibitem{duistermaat}
\textsc{J.J. Duistermaat}.
{\it On the Morse index in variational calculus},
Adv. in Math. {\bf 21} 173--195 (1976).

\bibitem{gutzwiller}
\textsc{M.C. Gutzwiller}.
\textit{Chaos in Classical and Quantum Mechanics},
Springer-Verlag: New York (1990).

\bibitem{jones}
\textsc{C.K.R.T. Jones}.
{\it Instability of standing waves for non-linear Schr\"odinger-type equations},
Ergodic Theory and Dynamical Systems {\bf 8*}, 119--138 (1988).

\bibitem{kawahara}
\textsc{T. Kawahara}.
{\it Oscillatory solitary waves in dispersive media},
J. Phys. Soc. Japan {\bf 33} 260--264 (1972).

\bibitem{ksk75}
\textsc{T. Kawahara, N. Sugimoto \& T. Kakutani}.
{\it Nonlinear interaction between short and long capillary-gravity
waves}, J. Phys. Soc. Japan {\bf 39}  1379--1386 (1975).

\bibitem{kodama}
\textsc{Y. Kodama \& D. Pelinovsky}.
{\it Spectral stability and time evolution of $N$-solitons in the {K}d{V} hierarchy},
J. Phys. A: Math. Gen. {\bf 38}, 6129--6140 (2005).

\bibitem{klr}
\textsc{A. Kushner, V. Lychagin \& V. Rubtsov}.
{\it Contact Geometry and Nonlinear Differential Equations},
Cambridge University Press (2007).

\bibitem{ll91}
\textsc{A. Latifi \& J. Leon}.
{\it On the interaction of Langmuir waves with acoustic waves in plasmas},
Phys. Lett. A {\bf 152} 171--177 (1991).

\bibitem{lewandosky}
\textsc{S. Lewandosky}.
{\it Stability of solitary waves of a fifth-order water wave model}, 
Physica D {\bf 227} 162--172 (2007).

\bibitem{lr87}
\textsc{R.G. Littlejohn \& J.M. Robbins}.
\textit{New way to compute Maslov indices},
Phys. Rev. A {\bf 36} 2953--2961 (1987).

\bibitem{ma79}
\textsc{Y.C. Ma}.
{\it On the multi-soliton solutions of some nonlinear
evolution equations},
Stud. Appl. Math. {\bf 60} 73--82 (1979).

\bibitem{sympl}
\textsc{D. McDuff \& D. Salamon}.
{\it Introduction to symplectic topology},
Oxford Mathematical Monographs,
Oxford University Press: New York (1995). 

\bibitem{montaldi}
\textsc{J. Montaldi}.
{\it A note on the geometry of linear Hamiltonian systems
of signature $0$ in $\R^4$}, Diff. Geom. Appl. {\bf 25} 344--350 (2007). 

\bibitem{m-g03}
\textsc{P. Muratore-Ginanneschi}.
{\it Path integration over closed loops and Gutzwiller's trace formula},
Phys. Rep. {\bf 383} 299--397 (2003).

\bibitem{pb03}
\textsc{M. Pletyukhov \& M. Brack}.
{\it On the canonically invariant calculation of Maslov indices},
J. Phys. A {\bf 36} 9449--9469 (2003).

\bibitem{rs93}
\textsc{J.W. Robbin \& D.A. Salamon}.
{\it The Maslov index for paths}, Topology {\bf 32} 827--844 (1993).

\bibitem{rs95}
\textsc{J.W. Robbin \& D.A. Salamon}.
{\it The spectral flow and the Maslov index},
Bull. London Math. Soc. {\bf 27} 1--33 (1995).

\bibitem{robbins91}
\textsc{J. Robbins}.
{\it Maslov indices in the Gutzwiller trace formula},
Nonlinearity {\bf 4} 343--363 (1991).

\bibitem{robbins92}
\textsc{J. Robbins}.
{\it Winding number formula for Maslov indices},
Chaos {\bf 2} 145--147 (1992).

\bibitem{vf92}
\textsc{A. Vanderbauwhede \& B. Fiedler}.
{\it Homoclinic period blow-up in reversible and conservative systems},
ZAMP {\bf 43} 292--318 (1992).

\bibitem{pw92}
\textsc{M.I. Weinstein \& R.L. Pego}.
{\it Eigenvalues, and instabilities of solitary waves}, 
Phil. Trans. Royal Soc. London A \textbf{340} 47--94 (1992).

\bibitem{wong}
\textsc{Y.C. Wong}.
{\it Differential geometry of Grassmann manifolds},
Proc. Nat. Acad. Sci. USA {\bf 57} 589--594 (1967).

\end{thebibliography}
\end{document}